\documentclass[draft,11pt]{article}

\usepackage{amsfonts,amsmath,amsthm,amscd,amssymb,latexsym,cite,verbatim,texdraw,floatflt,caption2,pb-diagram}

   \usepackage[T2A]{fontenc} % for russian and ukrainian language
    \usepackage[cp1251]{inputenc}
    \usepackage[ukrainian]{babel}

\usepackage[mathscr]{eucal}
\newtheorem{theorem}{Теорема}[section]
\newtheorem{lemma}{Лема}[section]

\newtheorem{corollary}{Наслідок}[section]
\newtheorem{proposition}{Твердження}[section]
\theoremstyle{definition}

\renewcommand{\abstract}{\textbf{Анотація. }\medskip}

\numberwithin{equation}{section}

\begin{document}

%%%%%%%%%%%%%%%%%%%%%%%%%%%%%%%%%%%%%%%%%%%%%%%%%%%%%%%%%%%%%%%%%%%%%%%%%%%%%
%%%%%%%%%%%%%%%%%%%%%%%%%%%%%%%%%%%%%%%%%%%%%(\ref {3})%%%%%%%%%%%%%%%%%%%%%%
%%%%%%%%%%%%%%%%%%%%%%%%%%%%%%%%%%%%%%%%%%%%%%%%%%%%%%%%%%%%%%%%%%%%%%%%%%%%%

\title{Задачі теорії наближень в абстрактних лінійних просторах}

\author{ А.\,С.~Сердюк, А.\,Л.~Шидліч}

\date{ }
\maketitle

{ \hfill \begin{minipage}[t]{6.5cm}
{\it

Присвячується світлій пам'яті

Олександра Івановича Степанця...}

 \end{minipage}}

\bigskip

\begin{center}

{\begin{minipage}[t]{13cm}
{  \small

В даній оглядовій роботі наведено результати, які охоплюють дослідження
актуальних проблем теорії апроксимації в абстрактних лінійних просторах. Такі дослідження набули активного розвитку, починаючи з 2000-х рокiв, на базі ідей та підходів, започаткованих в роботах О.\,І.~Степанця.  Зокрема, в огляді містяться результати, які стосуються  найкращих, найкращих $n$-членних наближень  та поперечників  деяких функціональних компактів у просторах ${\mathcal S}^p$, а також сформульовано прямі та обернені теореми наближення у цих просторах.}

 \end{minipage}}

\end{center}

%%%%%%%%%%%%%%%%%%%%%%%%%%%%%%%%%%%%%%%%%%%%%%%%%%%%%%%%%%%%%%%%%%%%%%%%%%%%%%%%%%%%%%%%%%%%

\section{Вступ}\label{Serdyuk_Shidlich_Intro}

Результати наукових досліджень, які будуть висвітлені у даній оглядовій роботі,
виникли внаслідок пошуку О.\,І.~Степанцем, його учнями та послідовниками нових підходів до  задач
теорії наближення функцій багатьох змінних і, зокрема,
періодичних функцій. В цій теорії існує багато проблем і одними з визначальних, напевно, є такі:
вибір апроксимативних агрегатів, вибір класів функцій та апроксимаційних характеристик. В той час, як в
одновимірному випадку вигляд найпростішого агрегату наближення визначається природним порядком натурального ряду, в
багатовимірному  випадку, тобто, коли задано абстрактну множину ${\mathscr X}$ -- банахів простір функцій $f({\bf x})=f(x_{1} \ldots,x_{d})$,
${\bf x}\in {\mathbb R}^d$, $d$ змінних, вибір найпростіших агрегатів є дещо проблематичним.
Перші труднощі тут починаються з того, що саме слід вважати  аналогом частинної суми для кратного ряду
  \begin{equation}\label{B.1}
  \sum_{{\bf k}\in {\mathbb Z}^d} c_{\bf k},\quad {\bf k}=(k_{1}, \ldots k_{d}),
  \end{equation}
 де ${\mathbb Z}^d$~--- цілочисельна решітка в ${\mathbb R}^{d}$. Природним є розгляд  ``прямокутних''\ сум
 і відповідних їм апроксимативних агрегатів -- у періодичному випадку тригонометричних поліномів вигляду
 \begin{equation}\label{B.2}
 \sum_{k_1=-n_1}^{n_1} \cdots \sum_{k_d=-n_d}^{n_d} c_{k_{1}, \ldots, k_{d}} {\mathrm e}^{{\mathrm i}(k_{1}t_{1}+\cdots+k_{d}t_{d})}.
  \end{equation}
Проте  частинні суми кратного ряду можна означати багатьма іншими способами, зокрема, наприклад, у такий спосіб.
Нехай $\{G_{\alpha}\}$~-- сім'я  обмежених областей взаємно неперетинних в  $ {\mathbb R}^d$, які залежать
від параметра $\alpha$, $\alpha\in {\mathbb N}$, і такі, що будь-який вектор ${\bf n}\in {\mathbb Z}^d$
належить усім областям $G_{\alpha}$ при достатньо великих значеннях $\alpha$.
 Тоді вираз
 $%\[
 \sum_{{\bf k}\in G_{\alpha}} c_{\bf k}
 $ %\]
називають частинною сумою ряду (\ref{B.1}), яка відповідає  області $G_{\alpha}$. За аналогією з цим вводяться і відповідні частинні суми тригонометричних рядів:
 \begin{equation}\label{B.3}
 \sum_{{\bf k}\in G_{\alpha}}c_{\bf k}{\mathrm e}^{{\mathrm i}{\bf k} {\bf x}} = \sum_{{\bf k}\in G_{\alpha}} c_{k_1 \ldots k_d} {\mathrm e}^{{\mathrm i}(k_1x_1+\cdots+k_d x_d)}.
 \end{equation}

Досить швидко виявилось, що у випадку наближення функцій з відомих класів Соболєва $W_p^r({\mathbb R}^d)$ замість прямокутних
сум вигляду (\ref{B.2}) доцільніше  застосовувати суми (\ref{B.3}), які побудовані за областями, що
визначаються деякими   гіперболічними поверхнями.
Такі області вперше були введені  К.\,І.~Бабенком в \cite{Babenko_1960_2, Babenko_1960_5} і отримали назву
гіперболічних хрестів. Їх поява дала істотний поштовх у розвитку сучасної теорії наближення функцій багатьох змінних.
В цьому напрямку отримано велику кількість важливих та цікавих результатів, з яким можна ознайомитись, наприклад, з робіт
\cite{Temlyakov_1986, Temlyakov_B1993, Romanyuk_2012, Dung_Temlyakov_Ullrich_2018}.

Слід зазначити, що більшість результатів, які стосуються наближення функцій з використанням
гіперболічних хрестів, у просторах $L_p({\mathbb R}^d)$ мають порядковий характер, а точні рівності  отримуються лише у
гільбертових просторах (при $p=2$). Спроби використання гіперболічних хрестів, а також їх модифікацій -- ступінчастих гіперболічних хрестів
при наближенні функцій з класів, відмінних від соболєвських,
взагалі кажучи, бажаних результатів майже не дають.  У зв'язку з цим, природно, виникає припущення,
що для кожного конкретного класу $\mathfrak N $ (або ж деякої сім'ї таких класів) потрібно підбирати
відповідну йому сім'ю областей $G_{\alpha}$, яка визначається його параметрами.

Іншою причиною, яка ускладнює отримання точних результатів  по наближенню функцій багатьох змінних
є історично сформована практика розглядати задачі саме у просторах $L_p({\mathbb R}^d)$. У періодичному випадку норма в цих просторах  означається рівністю
 \begin{equation}\label{Lp_norm}
 \|  f\|_{_{\scriptstyle L_p ({\mathbb  T}^d)}} =
 \bigg((2\pi)^{-d}\int%\limits
 _{{\mathbb T}^d} |f({\bf x})|^p {\mathrm d}{\bf x}\bigg)^\frac{1}{p},\quad  {\mathbb T}^d=[0,2\pi)^d.
 \end{equation}
і характеризує величину середнього значення $p$-го  степеня модуля заданої функції.

При  $p=2$  добре відомою є рівність Парсеваля
 \[
 \| f\|_{_{\scriptstyle L_2 ({\mathbb  T}^d)}} = \bigg(\sum\limits_{{\bf k}\in {\mathbb Z}^d } |c_{\bf k}|^2 \bigg)^\frac{1}{2},
 \]
де  $c_{\bf k} = c_{k_1, \ldots, k_d}$ -- коефіцієнти Фур'є функції
$f$. Тобто, у цьому випадку  норма функції $f$ повністю характеризує всю множину
$\{ c_{\bf k}\}_{{\bf k}\in {  {\mathbb Z}^d}}$ (при інших значеннях $p$, зрозуміло, подібні рівності можливі лише у
тривіальних випадках). Тому є доцільною спроба введення норм функцій за допомогою
величин, пов'язаних саме з їх коефіцієнтами Фур'є. Такий підхід розглядався, зокрема, у роботах \cite{Sterlin_1972, Kahan_M1976} та ін., але найбільш ретельно був розвинутий, починаючи з 2000-х років, у циклі робіт О.\,І.~Степанця та його послідовників
\cite{Stepanets_Preprint2001, Stepanets_UMZh2001_3, Stepanets_UMZh2001_8, Stepanets_M2002_2, Stepanets_UMZh2003_10, Stepanets_UMZh2004_10, Stepanets_UMZh2005_4, Stepanets_UMZh2006_1, Stepanets_B2006, Stepanets_Serdyuk_UMZh2002, Stepanets_Rukasov_UMZh2003_2, Stepanets_Rukasov_UMZh2003_5, Stepanets_Shydlich_UMZh2003, Stepanets_Shidlich_Pr_2007, Abdullayev_Ozkartepe_Savchuk_Shidlich_2019, Abdullayev_Chaichenko_Shidlich_2020, Abdullayev_Serdyuk_Shidlich_2020, Vakarchuk_2004, Vakarchuk_Shchitov_2006, Voicexivskij_2003, Voicexivskij_UMZh2003, Rukasov_UMZh2003, Savchuk_Shidlich_UMZh2007, Savchuk_Shidlich_2014, Serdyuk_2003, Timan_M2009, Shydlich_Zb2003, Shydlich_UMZh2004, Shidlich_UMZh2008, Shydlich_Zb_2008, Shydlich_UMZh2009, Shidlich_Zb2011, Shidlich_Zb2013, Shidlich_Zb2014, Shidlich_2016}.

Цей підхід, зокрема, дозволяє
розповсюджувати ідеї та методи теорії наближень  на абстрактні лінійні
простори, що в свою чергу, дає можливість дивитись на функції з загальних позицій аналізу та
дозволяє отримувати завершені  змістовні результати.

\section{ Наближення в просторах ${\mathcal S}^p_\varphi$}\label{Approx_Sp_phi}

\subsection{ Означення і деякі властивості просторів ${\mathcal S}^p_\varphi$.}
Нехай ${\mathscr X}$
-- деякий лінійний комплексний простір, $\varphi=\{\varphi_k\}_{k=1}^\infty$ -- фіксована зліченна
лінійно незалежна система в ньому, і нехай існує комплекснозначна функція $(x,y)$,
визначена  для кожної пари $x,y\in {\mathscr X}$, в якій хоча  б один із елементів
належить до $\varphi$, така, що виконуються  умови:

\quad 1) $(x,y) = \overline {(y,x)}$, де $\overline {z}$ ---
число, комплексно-спряжене  з $z$;

\quad 2) $(\lambda x_1 + \mu x_2,y) = \lambda (x_1,y) + \mu
(x_2,y),\quad \lambda ,\: \mu $ --- довільні числа;

\quad 3) $(\varphi_k,\varphi_l) = \left\{ \begin{matrix} 0,\quad k \neq
l;\cr 1,\quad k = l.\end{matrix} \right. $

\noindent Тобто, визначено скалярний добуток елементів простору
${\mathscr X}$ із елементами системи $\varphi$.

Кожному елементу  $x \in {\mathscr X}$ ставиться у відповідність
послідовність чисел
 $
 \widehat x_\varphi(k) = (x , \varphi_k)$, $k = 1,2,\ldots$ $(k \in {\mathbb N}),
 $
 i при даному фіксованому $p\in(0,\infty)$ розглядають простори ${\mathcal S}^p_\varphi = {\mathcal S}^p_\varphi({\mathscr X})$ всіх елементів
 $x \in  {\mathscr X}$ зі скінченною (квазі-)нормою
 \begin{equation}\label{Def_Sp}
 \| x\|_{_{\scriptstyle p}} := \| x\|_{_{\scriptstyle p,\varphi}} =
 \Big(\sum^\infty_{k=1}
 | \widehat x_\varphi(k)|^p \Big)^\frac{1}{p}.
 \end{equation}
При  цьому елементи $x,y\in {\mathscr X}$ вважаються тотожними в ${\mathcal S}^p_\varphi$,
якщо для будь-якого $k\in \mathbb N$ виконується рівність
$\widehat{x}_\varphi(k)= \widehat{y}_\varphi(k)$.

Зрозуміло, що при $ p=2$ простір $ {\mathcal S}_{\varphi }^2$ за умови його повноти є гільбертовим.
При всіх інших $p\in (0,\infty )$ простори $ {\mathcal S}^p_\varphi$ наслідують важливі властивості гільбертових
просторів -- рівність Парсеваля у вигляді співвідношення (\ref{Def_Sp}) і мінімальну властивість
частинних сум ряду Фур'є, яка формулюється в такий спосіб:

\begin{proposition}[\cite{Stepanets_Preprint2001, Stepanets_UMZh2001_3}]\label{Min_Fourier_Sums}
      Нехай $ f\in {{\mathcal S}^p_\varphi},$ $ p\in (0, \infty ),$
 \begin{equation}\label{b4}%$$
       S[f]= S[f]_{\varphi } = \sum _{k=1}^{\infty }\widehat f(k)\varphi _k
  \end{equation}
--  ряд Фур'є елемента $ f$  за системою  $ \varphi $  і
   \[
   S_n(f)=S_n(f)_{\varphi }=\sum _{k=1}^n\widehat f(k)\varphi _k, \ \   k\in {\mathbb N},
   \]
-- частинні суми цього ряду.  Серед  усіх поліномів вигляду
  $
   \Phi _n = \sum _{k=1}^n c_k\varphi _k
  $, $c_k\in {\mathbb C}$,
при даному $ n\in {\mathbb N}$ найменше відхиляється  від $ f$ частинна  сума $S_n(f)$, тобто,
$$
\inf\limits_{c_k}\|f- \Phi _n\|_{_{\scriptstyle p}} = \|f - S_n(f)\|_{_{\scriptstyle p}}.
$$
Крім того, виконується рівність
 \begin{equation}\label{b5}%$$
       \|f - S_n(f)\|_{_{\scriptstyle p}}^p=\|f\|_{_{\scriptstyle p}}^p - \sum _{k=1}^n|\widehat f(k)|^p=
       \sum _{k=n+1}^\infty |\widehat f(k)|^p.
  \end{equation}
\end{proposition}

При $ n\to \infty $ права частина в (\ref{b5}) прямує до нуля. Тобто,
для довільного  елемента $ f$ з ${\mathcal S}^p_\varphi$ його ряд Фур'є (\ref{b4}) збігається до
$ f,$   система $ \varphi $ є повною в ${\mathcal S}^p_\varphi$, а простір ${\mathcal S}^p_\varphi$ сепарабельний.

Звернемо увагу ще на одну властивість просторів ${\mathcal S}^p_\varphi:$ {\it якщо систему $ \varphi '=\{\varphi
'_k\}_{k=1}^{\infty }$ отримано із  системи $ \{\varphi _k\}_{k=1}^{\infty }$
шляхом будь-якої перестановки її членів, то справджуються рівності}
 \begin{equation}\label{b6}%$$
      {\mathcal S}^p_\varphi={\mathcal S}^p_{\varphi '}, \quad \mbox {і} \quad
      \|f\|_{_{\scriptstyle \varphi ,p}}=\|f\|_{_{\scriptstyle \varphi',p}}\ \ \ \forall f\in {\mathcal S}^p_\varphi.
   \end{equation}
Цей факт випливає з означення просторів ${\mathcal S}^p_\varphi$ та рівності (\ref{Def_Sp}).

Останнє зауваження дає можливість узагальнити твердження \ref{Min_Fourier_Sums} наступним чином.

\begin{proposition}[\cite{Stepanets_UMZh2006_1}]\label{Min_Fourier_Sums_2}
  Нехай $ \{g_{\alpha}\}$ -- сім'я обмежених підмножин, які залежать від параметра
  $\alpha\in {\mathbb N}$ і таких, що будь-яке число $ n\in {\mathbb N}$ належить усім множинам $ g_{\alpha}$
  з достатньо великими індексами $\alpha$.
  Нехай, далі, $ f\in {{\mathcal S}^p_\varphi}$, $ p\in (0, \infty ),$  і
   $$
    S_{g_{\alpha}}(f)=S_{g_{\alpha}}(f)_{\varphi }=\sum _{k\in g_{\alpha}}\widehat f(k)\varphi _k
    $$
-- частинна сума ряду $ S[f]_{\varphi },$ яка відповідає множині $ g_{\alpha}.$ Тоді серед
усіх   сум вигляду $\Phi _{g_{\alpha}}=\sum _{k\in g_{\alpha}}c_k\varphi _k$, $c_k\in {\mathbb C}$,
 найменше відхиляється  від $ f$ частинна сума  $ S_{g_{\alpha}}(f),$ тобто
$$
\inf\limits_{c_k}\|f-\Phi _{g_{\alpha}}\|_{_{\scriptstyle p}} = \|f- S_{g_{\alpha}}(f)\|_{_{\scriptstyle p}}.
$$
При цьому
$$
\|f-S_{g_{\alpha}}(f)\|_{_{\scriptstyle p}}^p=\|f\|_{_{\scriptstyle p}}^p - \sum _{k\in
g_{\alpha}}|\widehat f(k)|^p
$$
і
$$
\lim _{\alpha\to \infty }\|f-S_{g_{\alpha}}(f)\|_{_{\scriptstyle p}}=0.
$$
\end{proposition}

\subsection{Деякі реалізації та узагальнення.}\label{Examples}
Розглянемо декiлька прикладiв   реалiзацiй та узагальнень розглядуваних у п.~2.1 побудов (див., наприклад, \cite{Stepanets_UMZh2006_1}).

\vskip 2mm
\noindent{\bf \ref{Examples}.1. Простори ${\mathcal S}^p$.} Нехай ${\mathbb R}^d$ --  $d$-вимірний, $d\ge 1$, евклідів простір, ${\bf x}=(x_1,\ldots ,x_d)$ -- його елементи,
${\mathbb Z}^d$ -- цілочисельна решітка в ${\mathbb R}^d$, тобто, множина векторів ${\bf k}=(k_1,\ldots ,k_d)$ з
цілочисельними координатами, $({\bf x},{\bf y})=x_1y_1+\cdots
+x_dy_d$, $|{\bf x}|=\sqrt {({\bf x},{\bf x})}$ і, зокрема,
$({\bf k},{\bf x})=k_1x_1+\cdots +k_dx_d$,  $|{\bf k}|=\sqrt{k_1^2+\cdots +k_d^2}.$

Нехай, далі,  $L=L({\mathbb T}^d)$  -- множина всіх $2\pi$-періодичних за кожною зі
 змінних функцій $f({\bf x})=f(x_1,\cdots ,x_d)$, сумовних на кубі періодів ${\mathbb T}^d:=[0,2\pi)^d$.

 Якщо $f\in L$, то через $S[f]$  позначають  ряд Фур'є функції $f$ за тригонометричною системою
 $\{{\mathrm e}^{{\mathrm i}({\bf k},{\bf x})}\}_{{\bf k}\in {\mathbb Z}^d}$,    тобто
\begin{equation}\label{b8}%$$
    S[f]=\sum\limits_{{\bf k}\in {\mathbb Z}^d}\widehat{f} ({\bf k}){\mathrm e}^{{\mathrm i}({\bf k},{\bf x})}.
    \end{equation}
де
\begin{equation}\label{Fourier_Coeff}%$$
\widehat f({\bf k}):=(2\pi)^{-d}\int_{\mathbb T^d}f({\bf x}){\mathrm e}^{-{\mathrm i}({\bf k},{\bf x})}
 {\mathrm d}{\bf x},\quad{\bf k}\in\mathbb Z^d.
    \end{equation}
Якщо ототожнити  функції, еквівалентні відносно міри  Лебега, то за простір ${\mathscr X}$ можна
взяти простір $L({\mathbb T}^d)$, а за систему $\varphi$ -- тригонометричну систему
$\tau=\{\tau_s({\bf x})\}_{s=1}^\infty$, де
 \begin{equation}\label{b9}%$$
    \tau_s({\bf x})={\mathrm e}^{{\mathrm i} ({\bf k}_s,{\bf x})}, \quad {\bf k}_s\in {\mathbb  Z}^m,\ \ s=1,2,\ldots,
    %\eqno(1.7) $$
    \end{equation}
утворену із системи   $\{{\mathrm e}^{{\mathrm i}({\bf k},{\bf x})}\}_{{\bf k}\in {\mathbb Z}^d}$ шляхом довільної нумерації її елементів;
 скалярний добуток в такому випадку задається у відомий спосіб:
 \begin{equation}\label{b10}%$$
 (f,\tau_s)= (2\pi)^{-d}\int_{\mathbb T^d}f({\bf x})\overline\tau_s({\bf x}){\mathrm d}{\bf x}
 =\widehat{f}({\bf k}_s)=\widehat{f}_\tau(s).
     %\eqno(1.7) $$
    \end{equation}
Отримані при цьому множини  ${\mathcal S}^p_\tau$ згідно з  (\ref{b6})   не залежать від нумерації системи $\{{\mathrm e}^{{\mathrm i}({\bf k},{\bf x})}\}_{{\bf k}\in {\mathbb Z}^d}$ і надалі позначаються через  ${\mathcal S}^p$.

\vskip 2mm
\noindent{\bf \ref{Examples}.2. Простори $l_p$.} Виберемо тепер в ролі  ${\mathscr X}$ -- простір усіх послідовностей
$x=\{x_i\}_{i=1}^\infty$ комплексних чисел, у якому операції додавання та множення на скаляр визначаються в стандартний спосіб. У ролі  $\varphi$ --
 систему послідовностей $e = \{e_k\}_{k=1}^\infty,$ де $e_k =\{e_{ki}\}_{i=1}^\infty$ такі, що  $e_{kk}=1$ і $e_{ki}=0$ при $k\not=i$.

Скалярний добуток  елементів  $x\in {\mathscr X}$ на елементи системи $e$ визначимо співвідношеннями
 $$
 (x,e_k) = \widehat x_{e}(k) = x_k, \quad (e_k,x)=\overline{x_k},\quad k\in {\mathbb N}.
 $$
і при фіксованому  $p\in(0,\infty)$  розглянемо простори ${\mathcal S}^p_e({\mathscr X})$ всіх послідовностей
$x=\{x_i\}_{i=1}^\infty$ комплексних чисел зі скінченною (квазі-)нормою
 $$% \begin{equation}\label{Def_Sp}
 \| x\|_{_{\scriptstyle p,e}} = \Big(\sum^\infty_{k=1}  | \widehat x_e(k)|^p \Big)^\frac{1}{p}= \Big(\sum^\infty_{k=1}  |x_k|^p \Big)^\frac{1}{p}.
 \eqno(\ref{Def_Sp}')
 $$%\end{equation}
Очевидно, що ${\mathcal S}^p_e({\mathscr X})$ збігаються з відомими просторами послідовностей $l_p$.

\vskip 2mm
\noindent{\bf \ref{Examples}.3. Простори ${\mathcal S}^{p,\,\mu}_\varphi$} є деяким узагальненням просторів ${\mathcal S}^p_{\varphi}$. Вони були введені в роботі О.\,І.Степанця та В.\,І.~Рукасова \cite{Stepanets_Rukasov_UMZh2003_2} і будуються за тією ж схемою, що й останні, однак в цьому випадку функціонал вигляду
 $$
 \bigg(\sum^\infty_{k=1} {\big| \cdot \big|}^p
 \bigg)^\frac{1}{p}
 $$
у співвідношенні (\ref{Def_Sp}) слід замінити на функціонал з вагою $\mu$
  $$
 \bigg(\sum^\infty_{k=1} {\big| \cdot \big|}^p
 \mu^p_k \bigg)^\frac{1}{p},
 $$
де $\mu=\{\mu_k\}^\infty_{k=1}$ --- задана система невід'ємних
чисел, $\mu_k\geq 0$, $k\in  {\mathbb N}.$ При цьому якщо
$\mu_k\equiv 1,$ то ${\mathcal S}^{p,\,\mu}_\varphi ={\mathcal S}^p_{\varphi}$.

\vskip 2mm
\noindent{\bf  \ref{Examples}.4. Простори ${\mathcal S}^p_\Phi$} введено в 2003 році О.\,І.~Степанцем \cite{Stepanets_UMZh2003_10}.
При їх означенні використовуються подібна до наведених вище схема, яка полягає в наступному.
Нехай ${\mathscr X}$ та  ${\mathscr Y}$ -- деякі лінійні комплексні простори
векторів $x$ та $y$ відповідно. Припустимо, що на ${\mathscr X}$
задано лінійний оператор $\Phi$, який діє  в ${\mathscr Y}$, а на деякій
підмножині ${\mathscr Y}'\subset {\mathscr Y}$ визначено функціонал $f$. Нехай,
далі, $E(\Phi)$ -- множина значень оператора $\Phi$, і ${\mathscr X}'$ -- прообраз множини ${\mathscr Y}'\subset E(\Phi) $ при відображенні $\Phi$. В такому випадку на ${\mathscr X}'$ можна визначити
функціонал $f'$ за допомогою рівності
\begin{equation}\label{(f.1)}
 f'(x) =f (\Phi(x)), \qquad x\in {\mathscr X}'
\end{equation}
Якщо в ролі $f$ вибрати функціонал, що задає на ${\mathscr Y}' $
норму (або квазінорму), то рівність (\ref{(f.1)}) буде визначати
аналогічну величину на ${\mathscr X}'$.

Нехай $({\mathbb R}^d,d\mu)$, $d\ge 1$,
--- $d$-вимірний евклідів простір точок ${\bf x}=(x_1,\ldots,x_d)$, визначений на борелевій $\sigma$-алгебрі ${\mathcal B}$ зі скінченною $\sigma$-аддитивною неперервною мірою, $A$
--- $\mu$-вимірна підмножина з $({\mathbb R}^d,{\mathrm d}\mu)$,
$\mu$-міра якої дорівнює $a$, де $a$ --- або скінченне число, або
ж  $a=\infty$%:  $ \mbox{\rm mes}_\mu A=|A|_\mu=a,\  a\in (0,\infty]$
; ${\mathscr Y}={\mathscr Y}(A,{\mathrm d}\mu)$ --- множина всіх заданих на $A$ функцій $y=y({\bf
x})$, вимірних відносно міри ${\mathrm d}\mu$.

 При заданому $p\in (0,\infty]$ через $L_p(A,{\mathrm d}\mu)$ позначають
 підмножину  функций з ${\mathscr Y}(A,{\mathrm d}\mu)$, для яких є скінченною (квазі-)норма
\begin{equation}\label{(f.2)}
 \|y\|_{L_p(A,{\mathrm d}\mu)}= \left\{ \begin{matrix} \Big( \int_A
|\ y({\bf x})|^p\ {\mathrm d}\mu\Big)^{1/p}  , \quad\hfill & p \in
(0,\infty),\\
 \mathop{\mbox{ess\,sup}}\limits_{{\bf x}\in A} |\ y({\bf x})| ,\hfill & p = \infty.\end{matrix}\right.
 \end{equation}

 Нехай тепер ${\mathscr X}$ --- деякий лінійний простір векторів $x$, і $\Phi$ ---
 лінійний оператор, який діє з ${\mathscr X}$ в ${\mathscr Y}(A,{\mathrm d}\mu)$:
 $$ \Phi:{\mathscr X} \rightarrow {\mathscr Y}{(A,{\mathrm d}\mu)},\quad \Phi(x)
 \mathop{=}\limits^{\rm df}
 \widehat x,\quad x \in {\mathscr X},\quad \widehat x \in {\mathscr Y}{(A,{\mathrm d}\mu)}.$$
 При довільному фіксованому $p \in (0,\infty]$ покладають
\[
 {\mathcal S}_\Phi^p ={\mathcal S}_\Phi^p \: {({\mathscr X};{\mathscr Y})} =
\left\{ x \in {\mathscr X} : \|x\|_p = \|x\|_{p,\Phi} =\|\widehat x \|_{L_p(A,{\mathrm d}\mu)} <
\infty \right\}.
\]
Елементи $x_1, x_2\in {\mathscr X}$ вважають тотожними в ${\mathcal S}_\Phi^p$, якщо за
мірою ${\mathrm d}\mu$ майже скрізь \ \mbox{$\widehat x_1({\bf t}) = \widehat
x_2({\bf t})$.}

Таким чином, множина ${\mathcal S}_\Phi^p$ ---  множина всіх векторів $x{\in}
{\mathscr X}$, які є прообразами функцій з множини
$L_p{( A, {\mathrm d}\mu)}$  при \mbox{відображенні $\Phi$.}

Простори   ${\mathcal S}_\varphi^{p,\,\mu}$ (а отже, і ${\mathcal S}_\varphi^p$)
є частковими випадками просторів ${\mathcal S}_\Phi^p$. Дійсно, якщо в даному просторі ${\mathscr X}$ означити оператор
$\Phi$, який кожному $x \in {\mathscr X}$ ставить у відповідність
послідовність $y = \left\{y_k\right\}^\infty_{k=1}$, де $y_k=\widehat  x_\varphi(k)$; за
множину $({\mathbb R}^d,{\mathrm d}\mu)$ взяти простір ${\mathbb R}^1$ з
мірою ${\mathrm d}\mu$, носієм якої є множина $\mathbb Z^1$ цілочисельних
 точок $k$, в яких $\mu(k)\equiv \mu_k$; і  покласти $A =
\left\{ k \in {\mathbb Z}^1,\quad k \geq 1\:\right\}$, то в такому
випадку ${\mathscr Y}(A,{\mathrm d}\mu)$
--- множина всіх послідовностей $y$, для яких є скінченною величина
 $$
 \|y\|_{\scriptstyle_{L_p(A,{\mathrm d}\mu)}}\: = \bigg(\sum^\infty_{k=1}
 {|y_k|}^p \bigg)^\frac{1}{p},\quad p \in (0,\infty)
 $$

\subsection{ $\psi$-інтеграли та характеристичні послідовності.}\label{Charteristic}

\noindent{\bf \ref{Charteristic}.1. } У 2001 році О.\,І.~Степанець ввів до розгляду наступні об'єкти
наближення у просторах ${\mathcal S}^p_\varphi$, тобто,
підмножини елементів, які відповідають в класичній теорії апроксимації
поняттю класу функцій \cite{Stepanets_UMZh2001_3}, \cite[Гл.~11]{Stepanets_M2002_2}.

 Нехай $\psi = \{\psi_k \}_{k=1}^{\infty}$ -- довільна система
комплексних чисел. Якщо для данного елемента $f\in {\mathscr X}$,
ряд Фур'є якого має вигляд (\ref{b4}), існує елемент $F \in
{\mathscr X}$, для якого ряд Фур'є $S[F]_{\varphi}$ має вигляд
 \begin{equation}\label{b11}%$$
   S[F]_{\varphi} =\sum_{k=1}^\infty \psi_k \,\widehat f(k) \varphi_k,
 % \eqno (4.5)$$
 \end{equation}
     тобто, коли  \begin{equation}\label{b12}%$$
  \widehat F_\varphi(k) = \psi_k\,\widehat f(k),\quad k\in {\mathbb N}, %\eqno (4.6)$$
 \end{equation}
то елемент $F$ називається $\psi$-інтегралом елемента $f$. В такому
випадку записують $F=\mathcal{J}^{\psi}f$. Якщо $\mathfrak N$ -- деяка підмножина з ${\mathscr X}$, то через $\psi\mathfrak
N$ позначають множину $\psi$-інтегралів усіх елементів з $\mathfrak N$. Зокрема, $\psi {\mathcal S}^p_\varphi$ --
множина $\psi$-інтегралів всіх елементів, які належать даному
простору $ {\mathcal S}^p_\varphi$.

 Якщо $f$ і $F$ пов'язані співвідношенням (\ref{b11}) або (\ref{b12}), то  $f$ називають $\psi$-похідною елемента $F$ і позначають  $f = D^{\,\psi}F = F^{\,\psi}.$

Надалі обмежуємося випадком, коли система $\varphi$ задовольняє умову
  \begin{equation}\label{b13}%$$
  \lim_{ {k\rightarrow\infty}}|\psi_k| = 0. %\eqno(4.7)$$
  \end{equation}
Зрозуміло, що ця умова забезпечує вкладення $\psi {\mathcal S}^p_\varphi \subset {\mathcal S}^p_\varphi$, яке
має місце, зокрема,  за умови обмеженості множини чисел $|\psi_k|$, $k\in {\mathbb N}$.

Нехай
  \begin{equation}\label{b14}%$$
  U_{\varphi}^{p} =\,
   \left\{f\in {\mathcal S}^p_\varphi:\quad  \|f\|_{_{\scriptstyle p}} \leq 1  \right\} %\eqno (4.8) $$
    \end{equation}
-- одинична куля у даному просторі ${\mathcal S}^p_\varphi$ і $\psi U_{\varphi}^{p}$ --
множина $\psi$-інтегралів всіх елементів з $U_{\varphi}^{p}$. Саме множини $\psi U_{\varphi}^{p}$ є
основними об'єктами апроксимації в просторах ${\mathcal S}^p_\varphi$. Якщо
простір ${\mathcal S}^p_\varphi$ є повним, а
  \begin{equation}\label{b15}%$$
  \psi_k \neq 0 \quad \forall k \in {\mathbb N},
  %\eqno()$$
  \end{equation}
то внаслідок (\ref{b12}) та (\ref{b14})
  \begin{equation}\label{b16}%$$
      \psi\,U_\varphi^p = \bigg\{f \in {\mathcal S}^p_\varphi:\quad \sum^{\infty}_{k=1} \bigg |\frac{\widehat f(k)}{\psi_k}\bigg|^{p}
  \leq 1  \bigg\},%\eqno(4.8') $$
    \end{equation}
тобто, множина $\psi\,U_\varphi^p$ є $p$-еліпсоїдом в просторі ${\mathcal S}^p_\varphi$ з півосями, які дорівнюють $|\psi_k|$.

\vskip 2mm
\noindent{\bf \ref{Charteristic}.2. } Конструкцію агрегатів, які використовуються
 для наближення елементів $f \in {\mathcal S}^p_\varphi$, зручно визначати за допомогою спеціально підібраних
 характеристичних послідовностей $\varepsilon(\psi),\,g(\psi)$ і $\delta(\psi)$ системи $\psi$, які задаются в такий спосіб \cite{Stepanets_UMZh2001_3}, \cite[Гл.~11]{Stepanets_M2002_2}.

 Нехай $\psi = \{\psi_k \}_{k=1}^{\infty}$ -- довільна система комплексних чисел, які задовольняють умову
(\ref{b13}). Через $\varepsilon(\psi) = \{\varepsilon_1,\varepsilon_2,\ldots\} $ позначають
множину значень величин  $|\psi_k|$, впорядковану за їх спаданням, через $g(\psi) =\{ g_1,g_2,\ldots\} $ -- послідовність
множин
 $$
  g_n = g_n(\psi) =\left\{k \in {\mathbb N}:| \psi_k\,|\geq \varepsilon_n
 \right\}
 $$
 і через $\delta(\psi) = \delta_1,\delta_2,\ldots $ --  послідовність  чисел $\delta_n =|g_n|,$ де $|g_n|$--
кількість чисел $k \in {\mathbb N}$, які належать  множині $g_n$.  Через  $g_0 =g_0(\psi)$ позначають порожню
множину і вважають, що $\delta_{0}=0$.

Враховуючи умову (\ref{b13}), послідовності $\varepsilon(\psi)$ і $g(\psi)$ можна визначити такими співвідношеннями:
   \begin{equation}\label{b17i}%$$
\displaystyle{\begin{matrix}\varepsilon_1=\sup\limits_{k\in {\mathbb
N}}|\psi_k|,\   g_1=\{k\in {\mathbb N}:\
|\psi_k|=\varepsilon_1\},\quad \varepsilon_n=\sup\limits_{k\bar{\in} g_{n-1}}|\psi_k|,\cr\\ \
g_n=g_{n-1}\cup \{k\in {\mathbb N}: \ |\psi_k|=\varepsilon_n\},\quad n\in {\mathbb N}\setminus\{1\}.
\end{matrix}}%$$
  \end{equation}
За такого означення будь-яке число $n^{\ast} \in {\mathbb N}$ належить усім
множинам $g_n $ з достатньо великими номерами $n$ і
 \begin{equation}\label{b17}%$$
  \lim_{k\rightarrow\infty}\,\delta_k = \infty.
  \end{equation}

Зазначимо також, що якщо $\widetilde{\psi}=\{\widetilde{\psi}_k\}_{k=1}^\infty$ -- спадна
перестановка системи чисел $|\psi_k|$, $k=1,2,\ldots$, то має місце рівність
     \begin{equation}\label{b18}%$$
     \widetilde{\psi}_k=\varepsilon_n\quad \forall k\in
     (\delta_{n-1},\delta_n],\ \ n=1,2,\ldots .
     \end{equation}

\subsection{Найкращі наближення індивідуальних елементів множин $\psi {\mathcal S}^p_\varphi$.}\label{DI_TH}

 Нехай   $\psi=\{\psi_k\}_{k=1}^\infty$ -- довільна система комплексних чисел, підпорядкованих умові (\ref{b13}),
і  $\varepsilon(\psi),\,g(\psi)$ та $\delta(\psi)$ -- відповідні їй характеристичні послідовності.

Величину
    \begin{equation}\label{b21}%$$
          E_n(f)_{\psi,\,p} = \inf_{c_k\in {\mathbb C}} \Big\|\,f - \sum \limits_{k \in {g_{n-1}(\psi)}}c_k\,\varphi_k
    \Big\|_{_{\scriptstyle p}}
    \end{equation}
називають найкращим наближенням елемента $f\in {\mathcal S}^p_\varphi$ довільними   поліномами, побудованими по областях $g_{n-1}(\psi)$.

Наступне твердження встановлює зв'язок між найкращим  наближенням елемента $f$
і найкращими наближеннями його $\psi$-похідних.  Подібні твердження в теорії наближень прийнято називати прямими теоремами.

\begin{theorem}[{\cite[Гл.~11]{Stepanets_M2002_2}, \cite{Stepanets_UMZh2006_1}}]\label{Direct_Th} Нехай $f\in {\mathcal S}^p_\varphi$, $ p>0$ і система
$\psi =\{\psi_k\}_{k=1}^\infty$ підпорядкована умовам  (\ref{b13}) та (\ref{b15}).   Тоді ряд
 $$
  \sum\limits_
   {k=1}^{\infty}(\,\varepsilon^p_k - \varepsilon^p_{k-1})
   E^p_k(f)_{\psi,p}
 $$
 збігається і при  довільному $n\in {\mathbb N}$ справджується рівність
      \begin{equation}\label{b31}%$$
  E^p_n(f)_{\psi,p} = {\varepsilon}^{p}_{n} \,E^p_n(f^\psi)_{\psi,p} +
  \sum\limits^\infty_{k=n+1}\,(\varepsilon^p_k -
  \varepsilon^p_{k-1})\,E^p_k(f^\psi)_{\psi,p}\;, %\eqno (5.2)$$
  \end{equation}
  у якому величини $E_n(\cdot)_{\psi,p}$ визначаються рівністю
(\ref{b21}), а $\varepsilon_k$, $k=1,2,\ldots ,$  -- елементи
характеристичної послідовності $\varepsilon(\psi)$.
\end{theorem}

 Теорема \ref{In_Th} у певному розумінні є оберненою до попередньої: у ній за властивостями
 найкращого наближення елемента $f$ стверджується про існування у нього похідних і дається
 інформація про найкраще наближення цих похідних.

\begin{theorem}[{\cite[Гл.~11]{Stepanets_M2002_2}, \cite{Stepanets_UMZh2006_1}}]\label{In_Th}  Нехай $f\in {\mathcal S}^p_\varphi\cap \psi {\mathscr X}$, $p>0$, система $\psi =\{\psi_k\}^\infty_{k=1}$ підпорядкована умовам  (\ref{b13}) та (\ref{b15}) і
\begin{equation}\label{b32}%$$
  \lim\limits_{k\rightarrow\infty}\,\varepsilon^{-1}_k\,E_k(f)_{\psi,p}
  = 0.%\eqno (5.2)$$
  \end{equation}
 Тоді для того, щоб виконувалось включення $%\begin{equation}\label{b33}%$$
 f\in \psi {\mathcal S}^p_\varphi,$ %  \end{equation}
необхідно та достатньо, щоб збігався ряд
\begin{equation}\label{b34}%$$
  \sum\limits_{k=2}^{\infty}(\,\varepsilon^{-p}_k - \varepsilon^{-p}_{k-1})
   E^p_k(f)_{\psi,p} .%\eqno (5.2)$$
  \end{equation}
Якщо цей ряд збігається, то при довільному $n\in {\mathbb N}$ справджується рівність
\begin{equation}\label{b35}%$$
 E^p_n(f)_{\psi,p} = {\varepsilon}^{-p}_{n} \,E^p_n(f^\psi)_{\psi,p} +
  \sum\limits^\infty_{k=n+1}\,(\varepsilon^{-p}_k -
  \varepsilon^{-p}_{k-1})\,E^p_k(f)_{\psi,p}\;,%\eqno (5.2)$$
  \end{equation}
 у якому величини $E_n(\cdot)_{\psi,p}$ та $\varepsilon_k$ мають той же сенс, що і в теоремі \ref{Direct_Th}.
 \end{theorem}

%%%%%%%%%%%%%%%%%%%%%%%%%%%%%%%%%%%%%%%%%%%%%%%%%%%%%%%%%%%%%%%%%%%%%%%%%%%%%%%

\subsection{Найкращі наближення та базисні поперечники $q$-еліпсоїдів.}\label{Basis_Width}

\vskip 2mm
\noindent{\bf \ref{Basis_Width}.1. Означення найкращих наближень та базисних поперечників.}
 Нехай $f$ -- довільний елемент простору ${\mathcal S}^p_\varphi$ і
 $\gamma_n$, $n\in {\mathbb N}$, -- будь-який набір з $n$ різних натуральных чисел. Величину
    \begin{equation}\label{Best_Approx}%$$
          E_{\gamma_n}(f)_{p} = \inf_{c_k\in {\mathbb C}} \Big\|\,f - \sum \limits_{k \in {\gamma_n}}c_k\,\varphi_k
    \Big\|_{_{\scriptstyle p}}
    \end{equation}
називають найкращим наближенням елемента $f\in {\mathcal S}^p_\varphi$   $n$-членними поліномами, що відповідають набору $\gamma_n$.

Нехай, далі,
 $%\[
  S_{\gamma_n}(f)=S_{\gamma_n}(f)_{\varphi} =\sum_{k \in \gamma_n}\widehat{f}_{\varphi}(k)\varphi_k
 $%\]
-- сума Фур'є, яка відповідає набору $\gamma_n$, і
    \begin{equation}\label{Approx_Fourier_Sums}%$$
  {\mathscr E}_{\gamma_n}(f)_p\,=\|\,f - S_{\gamma_n}(f)
  \|_p
  \end{equation}
-- наближення елемента $f \in {\mathcal S}^p_\varphi $ сумою Фур'є, що відповідає набору $\gamma_n$.

 Якщо ${\mathfrak N}$ -- деяка підмножина простору ${\mathcal S}^p_\varphi$, то через
 $E_{\gamma_n}(\mathfrak N)_{p}$ та ${\mathscr E}_{\gamma_n}(\mathfrak N)_{p}$ позначають точні верхні межі
 величин  (\ref{Best_Approx}) та (\ref{Approx_Fourier_Sums}) по множині ${\mathfrak N}$, тобто,
    \begin{equation}\label{Best_Approx_Set}
          E_{\gamma_n}(\mathfrak N)_{p} = \sup\limits_{f\in\mathfrak  N} E_{\gamma_n}(f)_{p}\quad \mbox{\rm та}\quad
          {\mathscr E}_{\gamma_n}(\mathfrak N)_p = \sup\limits_{f\in\,\mathfrak N} {\mathscr E}_{\gamma_n}(f)_{p}.
     \end{equation}
Характеристики
\begin{equation}\label{b70}%  $$
   {\mathscr D}_n({\mathfrak N})_{p} = \inf\limits_{\gamma_n}
  E_{\gamma_n}({\mathfrak N})_{p}\quad \mbox{\rm  та}\quad
   {\mathscr D}_n^\perp({\mathfrak N})_{p} = \inf\limits_{\gamma_n}
  {\mathscr E}_{\gamma_n}({\mathfrak N})_{p}%   \eqno (4.22) $$
 \end{equation}
називають базисним та проєкційним поперечниками порядку  $n$ множини ${\mathfrak N}$
 в просторах ${\mathcal S}^p_\varphi$.

 Зазначимо, що у випадку наближення періодичних функцій  тригонометричними поліномами величинам ${\mathscr D}_n({\mathfrak N})_{p}$ відповідають    тригонометричні (базисні) поперечники, а величинам ${\mathscr D}_n^\perp({\mathfrak N})_{p}$  -- проєкційні (Фур'є) поперечники.

\vskip 2mm
\noindent{\bf \ref{Basis_Width}.2. Найкращі наближення та  поперечники $q$-еліпсоїдів в просторах ${\mathcal S}^p_\varphi$ при $0<q\le p$.}
Нехай $\psi = \{\psi_k \}_{k=1}^{\infty}$ -- довільна система
комплексних чисел, які задовольняють  умови (\ref{b13}) та (\ref{b15}) і $q$ -- довільне додатне число таке, що  $0<q\le p$.
У ролі множин ${\mathfrak N}$ у співвідношеннях (\ref{Best_Approx_Set}) та (\ref{b70}) будемо вибирати множини $\psi\,U_\varphi^q$
 $q$-еліпсоїдів в просторах ${\mathcal S}^p_\varphi$, які задаються рівністю (\ref{b16}) при $p=q$.

Оскільки (див., наприклад, \cite{Hardy_B1948}) для будь-якої невід'ємної послідовності $a=\{a_k\}_{k=1}^\infty$, $a_k\ge 0$,
  \begin{equation}\label{b47}%$$
 \bigg(\sum\limits_{k=1}^\infty a_k^p\bigg)^{1/p}\le
 \bigg(\sum\limits_{k=1}^\infty a_k^q\bigg)^{1/q}, \quad 0<q\le p,
 %\eqno (4.11) $$
 \end{equation}
то
 \begin{equation}\label{b48}%$$
 S^q_\varphi\subset {\mathcal S}^p_\varphi\quad \mbox{\rm та}\quad \psi U^q_\varphi\subset \psi U^p_\varphi, \quad 0<q\le p.
 %\eqno (4.11) $$
 \end{equation}

Для довільної системи комплексних чисел  $\psi=\{\psi_k\}^{\infty}_{k=1}$ та будь-якого набору
 $\gamma_n$ із $n$ різних натуральных чисел  через $\psi_{\gamma_n}=\{\psi_{\gamma_n}(k)\}_{k=1}^\infty$  позначимо послідовність чисел таку, що
 \begin{equation}\label{b71}%  $$
  \psi_{\gamma_n}(k)=\left\{\begin{matrix}0,\quad\hfill & k\in\gamma_n,\\ \psi_k,\quad\hfill & k\overline{\in} \gamma_n.\end{matrix}\right.
  \end{equation}

\begin{theorem}[\cite{Stepanets_UMZh2006_1}]\label{Best_approx_p<q}    Нехай  $\psi =
\{\psi_k\}^{\infty}_{k=1}$,  -- довільна система комплексних чисел, підпорядкована умовам
(\ref{b13}) та (\ref{b15}), і  $0<q\leq p$.
Тоді  для довільного набору  $\gamma_n$ із $n$ різних натуральных чисел, $n \in {\mathbb N}$, справджуються рівності
 \begin{equation}\label{b73}%  $$
   E_{\gamma_n}(\psi U_{\varphi}^{q})_{p}=
  {\mathscr E}_{\gamma_n}(\psi U_{\varphi}^{q})_{p} =\widetilde{\psi}_{\gamma_n}(1), %\eqno (4.26) $$
     \end{equation}
 де $\widetilde{\psi}_{\gamma_n}(1)$ -- перший член послідовності $\widetilde{\psi}_{\gamma_n}=\{\widetilde{\psi}_{\gamma_n}(k) \}^{\infty}_{k=1}$,
яка є спадною перестановкою послідовності  $\{|\psi_{\gamma_n}(k)|\}_{k=1}^\infty$.
\end{theorem}

Нехай $\widetilde{\psi} =\{\widetilde{\psi}_k\}^{ {\infty}} _{{k=1}}$ -- спадна перестановка послідовності $\{|\psi_k|\}^{
{\infty}}_{ {k=1}}$. Тоді, розглядаючи точні нижні межі обох частин рівності  (\ref{b73}) по всіх можливих наборах  $\gamma_n$,
неважко помітити, що точна нижня межа  правої частини (\ref{b73}) реалізується набором %$\gamma_n^*$ вигляду
 \begin{equation}\label{b80}%$$
 \gamma_n^*=\{i_k\in {\mathbb N}\ :\ |\psi_{i_k}|=\widetilde{\psi}_k,\quad k=1,2,\ldots,n\},
 %\eqno (4.11) $$
 \end{equation}
і при цьому $\widetilde{\psi}_{\gamma_n^*}(k)=\widetilde{\psi}_{n+k}$, $k=1,2,\ldots$.
 Тому внаслідок (\ref{b70})
 $$
 {\mathscr D}_n(\psi  U_\varphi^q)_p={\mathscr D}_n^\perp(\psi  U_\varphi^q)_p=\widetilde{\psi}_{\gamma_n^*}(1)=\widetilde{\psi}_{n+1}.
 $$
Отже, має місце таке твердження про точні значення поперечників.

\begin{theorem}[\cite{Stepanets_UMZh2006_1}]\label{Basis_Width_p<q}  Нехай  $\psi =
\{\psi_k\}^{\infty}_{k=1}$,  -- довільна система комплексних чисел, підпорядкована умовам (\ref{b13}) та (\ref{b15}),
і $0<q\leq p$. Тоді при кожному $n
\in {\mathbb N}$  справджуються рівності
 \begin{equation}\label{b81}%  $$
 {\mathscr D}_n(\psi
 U_\varphi^q)_p={\mathscr D}_n^\perp(\psi  U_\varphi^q)_p=\widetilde{\psi}_{n+1}, %\eqno (4.26) $$
     \end{equation}
 де  $\widetilde{\psi}=\{\widetilde{\psi}_k \}^{\infty}_{k=1}$ -- спадна перестановка  послідовності $\{|\psi_k|\}_{k=1}^\infty$.
 При цьому для набору чисел $\gamma_n^*$, означеного    рівністю (\ref{b80}), мають місце рівності
\[
 {\mathscr D}_n(\psi
 U_\varphi^q)_p={\mathscr D}_n^\perp(\psi  U_\varphi^q)_p=E_{\gamma_n^*}(\psi U_{\varphi}^{q})_{p}
 \]
  \begin{equation}\label{b82}%  $$
  =
  {\mathscr E}_{\gamma_n^*}(\psi U_{\varphi}^{q})_{p} =
\widetilde{\psi}_{\gamma_n^*}(1)=\widetilde{\psi}_{n+1}. %\eqno (4.26) $$
     \end{equation}
     \end{theorem}

\vskip 2mm
\noindent{\bf \ref{Basis_Width}.3. Найкращі наближення та базисні поперечники $q$-еліпсоїдів в просторах ${\mathcal S}^p_\varphi$ при $0<p<q$.}
Наведемо  аналоги теорем \ref{Best_approx_p<q} та \ref{Basis_Width_p<q}  у випадку, коли $q>p>0$.  Як і вище, припускаємо,
 що система чисел $\psi$ підпорядкована умові (\ref{b15}),  а також умові
  \begin{equation}\label{b54}%$$
       \|\psi\|_{l_{\frac{pq}{q-p}}}=\bigg(\sum\limits_{k=1}^\infty|\psi_k|^{\frac{pq}{q-p}}\bigg)^{\frac{q-p}{pq}} <\infty, \quad 0<p<q.
 %\eqno (4.11) $$
 \end{equation}
яка забезпечує вкладення
  \begin{equation}\label{b55}%$$
\psi U^q_\varphi\subset {\mathcal S}^p_\varphi, \quad 0<p<q.
 %\eqno (4.11) $$
 \end{equation}

\begin{theorem}[\cite{Stepanets_UMZh2006_1}]\label{Best_approx_p>q}  Нехай $0<p<q$, і $\psi=\{\psi_k\}_{k=1}^\infty$  --
система чисел, підпорядкована умовам  (\ref{b15}) та (\ref{b54}). Тоді для довільного набору  $\gamma_n$ із $n$ різних натуральных чисел, $n \in {\mathbb N}$,  справджуються рівності
 \begin{equation}\label{b84}%  $$
 E_{\gamma_n}(\psi U_\varphi^q)_{p} = {\mathscr E}_{\gamma_n}(\psi U_\varphi^q)_{p}
 =\bigg( \sum_{k=1}^{\infty} \,(\widetilde{\psi}_{\gamma_n}(k))^{\frac{p\,q}{q-p}}\bigg)^{\frac{q-p}{p\,q}},
 %   \eqno (4.22) $$
 \end{equation}
де $\widetilde{\psi}_{\gamma_n}=\{\widetilde{\psi}_{\gamma_n}(k) \}^{\infty}_{k=1}$ -- спадна перестановка послідовності
$\{|\psi_{\gamma_n}(k)|\}_{k=1}^\infty$.
\end{theorem}

Розглядаючи точні нижні межі обох частин рівності (\ref{b84}) по всіх можливих наборах  $\gamma_n$,
можна переконатись, що точна нижня межа  правої частини (\ref{b84}) реалізується набором  $\gamma^{\ast}_n$, який визначається
 співвідношенням (\ref{b80}).  Тому внаслідок (\ref{b70})
 $$
   {\mathscr D}_n(\psi\,U_\varphi^q)_{p} \!=\!{\mathscr D}_n^\perp(\psi  U_\varphi^q)_p=\!\!\bigg( \sum_{k=1}^{\infty} \,(\widetilde{\psi}_{\gamma^{\ast}_n}(k) )^{\frac{p\,q}
   {q-p}}\!\!\bigg)^{\!\!\frac{q-p}{p\,q}} \!\!\!\!=\!\! \bigg(\sum_{k=n+1}^{\infty}
   \widetilde{\psi}_k^{ \frac{p\,q}{q-p}}\bigg)^{\!\!\frac{q-p}{p\,q}}\!\!\!\!.
$$
Отже, справджується таке твердження про значення поперечників.

\begin{theorem}[\cite{Stepanets_UMZh2006_1}]\label{Basis_Width_p>q} Нехай  $0<p<q$  і $\psi=\{\psi_k\}_{k=1}^\infty$  --
система чисел, підпорядкована умовам  (\ref{b15}) та (\ref{b54}). Тоді при кожному
$n\in {\mathbb N}$  мають місце рівності
\begin{equation}\label{b92}% $$
  {\mathscr D}_n(\psi\,U_\varphi^q)_{p}={\mathscr D}_n^\perp(\psi  U_\varphi^q)_p=\bigg( \sum_{k=n+1}^{\infty} \widetilde{\psi}_k^{\frac{p\,q}{q-p}}
  \bigg)^{\frac{q-p}{p\,q}},
% \eqno(3.75)
 \end{equation}
 де  $\widetilde{\psi}=\{\widetilde{\psi}_k \}^{\infty}_{k=1}$ -- спадна перестановка   послідовності $\{|\psi_k|\}_{k=1}^\infty$.
При цьому для набору чисел $\gamma_n^*%=\{i_k\}_{k=1}^n
 $, який визначається  рівністю (\ref{b80}), справджуються рівності
$$
 {\mathscr D}_n(\psi
 U_\varphi^q)_p={\mathscr D}_n^\perp(\psi  U_\varphi^q)_p=E_{\gamma_n^*}(\psi U_{\varphi}^{q})_{p}=
  {\mathscr E}_{\gamma_n^*}(\psi U_{\varphi}^{q})_{p}
   $$
    \begin{equation}\label{b93}
=\bigg( \sum_{k=1}^{\infty} \,(\widetilde{\psi}_{\gamma^{\ast}_n}(k))^{\frac{p\,q}
   {q-p}}\bigg)^{\frac{q-p}{p\,q}}=\bigg( \sum_{k=n+1}^{\infty} \widetilde{\psi}_k^{\frac{p\,q}{q-p}}
  \bigg)^{\frac{q-p}{p\,q}}. %\eqno (4.26) $$
     \end{equation}
     \end{theorem}

Звернемо увагу на те, що послідовність $\widetilde{\psi}_k$,
 $k\in {\mathbb N}$,  в загальному випадку є ступінчастою. Тому внаслідок (\ref{b81}) такий самий характер має і величина
  ${\mathscr D}_n(\psi U_\varphi^q)_{p}$ при $p\geq q > 0$. Якщо ж $p<q$,  то згідно з (\ref{b92}) ця величина строго спадає з ростом параметра
  $n$.

Зазначимо, що  інтегральні аналоги теорем \ref{Best_approx_p<q}--\ref{Basis_Width_p>q} встановлено в роботах
\cite{Stepanets_UMZh2003_10, Stepanets_Shidlich_Pr_2007, Stepanets_Shidlich_JAT_2010, Stepanets_Shidlich_IZV_2010, Shidlich_Chaichenko_LM_2014}.
У роботах \cite{Shidlich_Chaichenko_lp_2014}, \cite{Shidlich_Chaichenko_Orlicz_lM_2015} та \cite{Chaichenko_Shidlich_2018}
твердження теорем \ref{Best_approx_p<q}--\ref{Basis_Width_p>q}  було поширено відповідно на простори зі змінним показником підсумовування $l_{\bf p}$,  простори Орлича $l_{M}$  та модулярні простори
Мусєляка-Орлича  $l_{\bf M}$.

\subsection{Наближення поліномами, побудованими по областях $g_n(\psi)$, та колмогоровські поперечники $p$-еліпсоїдів.}\label{Kolmogorov_Width}

\noindent{\bf \ref{Kolmogorov_Width}.1. Наближення поліномами, побудованими по областях $g_n(\psi)$.}
 Нехай  $\psi=\{\psi_k\}_{k=1}^\infty$ -- довільна система комплексних чисел, які задовольняють умову (\ref{b13}).
Розглянемо окремо випадок, коли апроксимуючі поліноми визначаються областями  $g_n(\psi)$, побудованими по даній системі комплексних чисел
$\psi$ згідно з формулами (\ref{b17i}). Для довільного елемента $f\in \psi {\mathcal S}^p_\varphi$ позначимо
     \begin{equation}\label{b19}%$$
   S_n(f)_{\varphi,\psi}:= S_{g_n(\psi)}(f) = \sum_{k \in
  g_n(\psi)}\,\widehat f(k)\,\varphi_k,\quad
  S_0(f)_{\varphi,\psi} = \theta, %\eqno (4.9)$$
       \end{equation}
    де   $g_n(\psi)$, $n=1,2,\ldots$, -- елементи
послідовності $g(\psi)$, $\theta$ -- нульовий елемент
простору ${\mathcal S}^p_\varphi$, і
  \begin{equation}\label{b20i}%$$
   {\mathscr E}_n(f)_{\psi,\,p}\,= \|\,f - S_{n-1}(f)_{
   {\varphi,\,\psi}}\|_{_{\scriptstyle p}}.
  % \eqno (4.10)$$
    \end{equation}
Якщо  ${\mathfrak N}$ -- деяка підмножина з  $\psi {\mathcal S}^p_\varphi$, то покладаємо
     \begin{equation}\label{b20}%$$
    {\mathscr E}_n({\mathfrak N})
   _{\psi,\,p}\,= \sup\limits_{f \in {\mathfrak N}}\,{\mathscr E}_n(f)_{
   {\psi,\,p}}.
  % \eqno (4.10)$$
    \end{equation}
та
  \begin{equation}\label{b22}%$$
      E_n({\mathfrak N})_ {\psi,\,p} = \sup\limits_{f \in  {\mathfrak N}}\,E_n(f)_{\psi,\,p}.
    \end{equation}

Враховуючи рівності (\ref{b18}) та прийняті позначення, наводимо такі твердження
-- наслідки з теорем \ref{Best_approx_p<q} та~\ref{Best_approx_p>q}.

\begin{corollary}[\cite{Stepanets_UMZh2006_1}]\label{Best_approx_g_n_p<q}    Нехай  $\psi =
\{\psi_k\}^{\infty}_{k=1}$,  -- довільна система комплексних чисел, підпорядкована умовам
(\ref{b13}) та (\ref{b15}), і   $0<q\leq p$.
Тоді  для кожного $n \in {\mathbb N}$  справджуються рівності
   \begin{equation}\label{b23}%$$
   E_n(\psi U_\varphi^p)_{\psi,p}=  {\mathscr E}_n(\psi \,U_\varphi^p)_{\psi,\,p}=\varepsilon_n,
      %\eqno (4.11) $$
    \end{equation}
    де $\varepsilon_n$ -- $n$-й член характеристичної
    послідовності $\varepsilon(\psi)$.
\end{corollary}

\begin{corollary}[\cite{Stepanets_UMZh2006_1}]\label{Best_approx_g_n_p>q} Нехай  $\psi=\{\psi_k\}_{k=1}^\infty$  --
система комплексних чисел, підпорядкована умовам  (\ref{b15}) та (\ref{b54}), і $0<p<q$. Тоді  для кожного $n \in {\mathbb N}$  справджуються рівності
  \begin{equation}\label{b57}%$$
  E_n(\psi U^q_\varphi)_{\psi,p} ={\mathscr E}_n(\psi U^q_\varphi)_{\psi,p} =
 \bigg(\sum\limits_{k=\delta_{n-1}+1}^\infty \widetilde{\psi}_k^{\frac{pq}{q-p}}\bigg)^{\frac{q-p}{pq}}\,,
  %\eqno (4.11) $$
 \end{equation}
  де $\widetilde{\psi}=\{\widetilde{\psi}_k\}_{k=1}^\infty$ -- спадна перестановка  послідовності $\{|\psi_k|\}_{k=1}^\infty$, а
 $\delta_n$ -- члени характеристичної послідовності $\delta(\psi)$.
 \end{corollary}

\vskip 2mm
\noindent{\bf \ref{Kolmogorov_Width}.2. Колмогоровські поперечники $p$-еліпсоїдів.} Нехай $Y$ -- лінійний нормований простір, $\mathfrak M$ -- центрально-симетрична множина в ньому і $\mathcal{F}_n$ -- множина всіх підпросторів $F_n$ розмірності
$n \in {\mathbb N}$ простору $Y$. Величину
$$
  d_n(\mathfrak M;Y) = \inf\limits_{F_n \in {\mathcal{F}}_{n}}\,
  \sup \limits_{x \in \mathfrak M}\,\inf\limits_{u \in F_n}\,
  \|x - u\|_{ {Y}}
$$
називають поперечником за  Колмогоровим множини ${\mathfrak M}$ у просторі $Y$.

\begin{theorem}[\cite{Stepanets_UMZh2006_1}]\label{Kolmogorov_Width_p=q} Нехай $\psi=\{\psi_k\}$, $k=1,2,\ldots$,
-- система комплексних чисел, підпорядкована умовам (\ref{b13}) та (\ref{b15}). Тоді при довільних
$p\in [1,\infty)$ та $n\in {\mathbb N}$ справджуються рівності
   \[%
   d_{\delta_{n-1}}(\psi U_\varphi^p)=  d_{\delta_{n-1}+1}(\psi U_\varphi^p)=\ldots
   \]
   \begin{equation}\label{b27}%$$
   =
     d_{\delta_{n}-1}(\psi U_\varphi^p)=E_n(\psi U_\varphi^p)_{\psi,p}=\varepsilon_n,
      %\eqno (4.11) $$
   \end{equation}
у яких  $\delta_s$ та $\varepsilon_s$, $s=1,2,\ldots,$ -- елементи характеристичних
послідовностей $\delta(\psi)$ та $\varepsilon(\psi)$ системи $\psi$, а $\delta_0=0$.
     \end{theorem}

Зазначимо, що у скінченно вимірних просторах $l_p^d$ твердження, аналогічне до теореми \ref{Kolmogorov_Width_p=q},
випливає з теореми 2.1 глави VI монографії А.~Пінкуса \cite{Pinkus_1985}. У роботах \cite{Shidlich_Chaichenko_lp_2014}, \cite{Shidlich_Chaichenko_Orlicz_lM_2015} та \cite{Chaichenko_Shidlich_2018}
твердження теореми \ref{Kolmogorov_Width_p=q}   поширено відповідно на простори зі змінним показником підсумовування $l_{\bf p}$,  простори Орлича $l_{M}$  та модулярні простори Мусєляка-Орлича  $l_{\bf M}$.

%%%%%%%%%%%%%%%%%%%%%%%%%%%%%%%%%%%%%%%%%%%%%%%%%%%%%%%%%%%%%%%%%%%%%%%%%%%%%%%

\subsection{Найкращі $n$-членні наближення.} \label{Best_n-term}

\noindent{\bf \ref{Best_n-term}.1. Найкращі $n$-членні наближення  $q$-еліпсоїдів в просторах ${\mathcal S}^p_\varphi$ при $0<q\le p$.}
 Нехай $n\in {\mathbb N}$, $\gamma_n$ -- довільний
набір з $n$ натуральних чисел і
     \begin{equation}\label{b94}%  $$
P_{\gamma_n}=\sum\limits_{k\in
\gamma_n}\alpha_k\varphi_k, %\eqno (4.26) $$
     \end{equation}
     де $\alpha_k$ -- комплексні числа.    Величину
    \begin{equation}\label{b95}%
         e_n(f)_p=e_n(f)_{\varphi,p}=\inf\limits_{\alpha_k,\,\gamma_n}\|f-P_{\gamma_n}\|_p.
    \end{equation}
називають  найкращим $n$-членним наближенням елемента $f\in {\mathcal S}^p_\varphi$
в просторі ${\mathcal S}^p_\varphi$. Якщо ${\mathfrak N}$ -- деяка підмножина з ${\mathcal S}^p_\varphi$, то покладаємо
     \begin{equation}\label{b961}%  $$
         e_n({\mathfrak N})_p=\sup\limits_{f\in {\mathfrak N}} e_n(f)_p.
     \end{equation}

Величини, аналогічні до величин (\ref{b95}), вперше введені  С.\,Б.~Стєчкіним \cite{Stechkin_1955_DAN},
і їх властивості досліджувались в теорії наближень  періодичних функцій багатьма авторами (див., наприклад, \cite{Romanyuk_2012, Dung_Temlyakov_Ullrich_2018, Temlyakov_B2011_Greedy, Temlyakov_B2015, Stepanets_M2002_2} та ін.). Варто зазначити, що раніше Е.~Шмідт \cite{Schmidt_1906} розглядав величину найкращого білінійного наближення, яка в ідейному плані є близькою до величин вигляду (\ref{b95}).

В даному підрозділі визначаються величини вигляду (\ref{b961}) у випадку, коли  ${\mathfrak N}$ є
$q$-еліпсоїдами $\psi\,U_\varphi^q$ (див. означення (\ref{b16})). Як і вище, вважаємо, що $\psi = \{\psi_k \}_{k=1}^{\infty}$ -- довільна система
комплексних чисел, які задовольняють  умови (\ref{b13}) та (\ref{b15}). В такому випадку, як вже зазначалося, при $0<q\le p$
виконується вкладення  $\psi U^q_\varphi\subset {\mathcal S}^p_\varphi$, а отже,
величина (\ref{b961}) має зміст.

\begin{theorem}[{\cite[Гл.~11]{Stepanets_M2002_2}, \cite{Stepanets_UMZh2006_1}}]\label{Best_n_term_approx_p>q}
         Нехай $\psi=\{\psi_k\}_{k=1}^\infty$  -- система чисел, підпорядкована умовам (\ref{b13}) та (\ref{b15}),
         і  $0<q\le p$.
         Тоді при довільному $n\in {\mathbb N}$  справджується  рівність
\begin{equation}\label{b97}%$$
         e_n^p(\psi U^q_\varphi)_{p} =\sup\limits_{s>n}
         (s-n) \bigg(\sum\limits_{k=1}^s \widetilde{\psi}_k^{-q}\bigg)^{-\frac pq},
  %\eqno (4.26) $$
 \end{equation}
де $\widetilde{\psi}=\{\widetilde{\psi}_k\}_{k=1}^\infty$ -- спадна перестановка послідовності чисел
$\{|\psi_k|\}_{k=1}^\infty$. Точна верхня межа в правій частині (\ref{b97}) досягається при деякому скінченному значенні $s^*$.
\end{theorem}

У випадку, коли всі числа послідовності $\psi$  дорівнюють одиниці, тобто, коли $\psi U^q_\varphi=U^q_\varphi$, $0<q\le p$, має місце таке твердження.

\begin{theorem}[{\cite[Гл.~11]{Stepanets_M2002_2}, \cite{Stepanets_UMZh2006_1}}]\label{Best_n_term_approx_psi_1}
         Нехай  $0<q\le p$.      Тоді при довільному $n\in {\mathbb N}$  справджується  рівність
  \begin{equation}\label{b125}%$$
         e_n^p( U^q_\varphi)_{p} =
         \sup\limits_{s>n} \frac{s-n}{s^{1/q}}.
 \end{equation}
При  $p=q$ точна верхня межа в правій частині (\ref{b125}) дорівнює одиниці. Якщо ж
$q<p$, то вона досягається в одній із точок $\Big[\frac {n}{1-q/p}\Big]$ або ж $\Big[\frac {n}{1-q/p}\Big]+1$, де $[c]$ --- ціла частина числа $c$.
\end{theorem}

\vskip 2mm
\noindent{\bf \ref{Best_n-term}.2. Найкращі $n$-членні наближення  $q$-еліпсоїдів в просторах ${\mathcal S}^p_\varphi$ при $0<p<q$.}
 В цьому підрозділі наведено точні значення величини $e_n(\psi U^q_\varphi)_p$ за умови, що  $0<p<q$. Як і вище, припускаємо, що система  чисел $\psi$ підпорядкована умові  (\ref{b15}), а також  умові (\ref{b54}),
яка гарантує вкладення $\psi U^q_\varphi\subset S^p_\varphi$.

\begin{theorem}[{\cite[Гл.~11]{Stepanets_M2002_2}, \cite{Stepanets_UMZh2006_1}}]\label{Best_n_term_approx_p<q}
      Нехай $0<p<q$ і  $\psi=\{\psi_k\}_{k=1}^\infty$  --       система чисел, підпорядкована умовам  (\ref{b15}) та (\ref{b54}).
      Тоді при довільному $n\in {\mathbb N}$  справджується рівність
  \begin{equation}\label{b127}%$$
  e_n(\psi U^q_\varphi)_{p} =
  \bigg((s^*-n)^\frac q{q-p} \bigg(\sum\limits_{k=1}^{s^*} \widetilde{\psi}_k^{-q}
  \bigg)^ {-\frac p{q-p}}+
 \sum\limits_{k=s^*+1}^\infty  \widetilde{\psi}_k^\frac
 {pq}{q-p}\bigg)^\frac{q-p}{pq},
  %\eqno (4.26) $$
 \end{equation}
де $\widetilde{\psi}=\{\widetilde{\psi}_k\}_{k=1}^\infty$ -- спадна перестановка послідовності чисел
$\{|\psi_k|\}_{k=1}^\infty$, а  число $s^*$ вибране з умови
  \begin{equation}\label{b128}%$$
 \widetilde{\psi}_{s^*}^{-q}\le \frac 1{s^*-n}\sum\limits_{k=1}^{s^*}\widetilde{\psi}_k^{-q}
 <\widetilde{\psi}_{s^*+1}^{-q}.
  %\eqno (4.26) $$
 \end{equation}
 Таке число $s^*$ існує і єдине.
 \end{theorem}

Зазначимо, що  аналоги теорем \ref{Best_n_term_approx_p>q}, \ref{Best_n_term_approx_psi_1}
та \ref{Best_n_term_approx_p<q} у випадку апроксимації інтегралами заданого рангу встановлено в роботах \cite{Stepanets_UMZh2003_10, Stepanets_Shidlich_Pr_2007, Stepanets_Shidlich_IZV_2010}.  У просторах $l_p^d$ скінченних послідовностей аналогічні твердження при всіх $0<p,\,q\le \infty$ отримано в роботі \cite{Fuchang_Gao_2010}.

У роботі \cite{Rukasov_UMZh2003} твердження теорем  \ref{Best_n_term_approx_p>q} та \ref{Best_n_term_approx_p<q} розповсюджено
на простори  ${\mathcal S}^{p,\mu}_{\varphi}$, а в
\cite{Shidlich_Chaichenko_Orlicz_lM_2015}
твердження теореми \ref{Best_n_term_approx_p>q} дещо розповсюджено  на простори Орлича $l_{M}$.

%%%%%%%%%%%%%%%%%%%%%%%%%%%%%%%%%%%%%%%%%%%%%%%%%%%%%%%%%%%%%%%%%%%%%%%%%%%%%%%%%%%%

\subsection{Порядкові оцінки найкращих $n$-членних наближень та  поперечників  $q$-еліпсоїдів в просторах ${\mathcal S}^p_\varphi$.}\label{Order_Estimates}

Аналіз  теорем \ref{Best_approx_p<q}--\ref{Basis_Width_p>q}, \ref{Kolmogorov_Width_p=q} та \ref{Best_n_term_approx_p>q}--\ref{Best_n_term_approx_p<q}, показує, що  точні значення апроксимативних характеристик у теоремах \ref{Best_approx_p<q}, \ref{Basis_Width_p<q} та \ref{Kolmogorov_Width_p=q},
виражаються у термінах величин, для яких явно прослідковується їх  швидкість прямування до нуля при $n\to \infty$.
Вирази, у термінах яких визначені точні значення апроксимативних характеристик у теоремах
\ref{Basis_Width_p>q}, \ref{Best_n_term_approx_p>q} та \ref{Best_n_term_approx_p<q}, потребували додаткових досліджень.
Такі дослідження були здійснені, зокрема, у роботах \cite{Stepanets_Shidlich_Pr_2007, Stepanets_Shidlich_NK2007, Stepanets_Shidlich_IZV_2010, Shidlich_2016}. При цьому ефективним виявився розвинений О.\,І.~Степанцем та його учнями апарат дослідження,
 який базується на наведеній нижче класифікації опуклих  функцій \cite[Гл.3]{Stepanets_M2002_1}.

\vskip 2mm
\noindent{\bf \ref{Order_Estimates}.1. Класифікація  Степанця опуклих функцій.}
Нехай  ${\mathfrak M}$ -- множина всіх опуклих донизу функцій
$\psi(t)$, неперервного аргументу $t\in [1,\infty)$, які
задовольняють умову $\mathop{\lim}\limits_{t\rightarrow\infty}\psi(t)=0$:
 \[
 {\mathfrak {M}}=\{\psi (t):\ \psi(t)>0, \ \psi (t_1)-2\psi
 ((t_1+t_2)/2)
 \]
 \[
 +\psi (t_2)\ge 0\ \ \forall t_1,t_2\in [1,\infty ),
 \lim _{t\to \infty }\psi (t)=0\}.
 \]
Множина ${\mathfrak {M}}$ досить неоднорідна за швидкістю прямування
до нуля при $t\to \infty $ її елементів: функції $\psi (t)$ можуть
спадати як дуже повільно, так і дуже швидко. Тому виникає
необхідність розбиття множини ${\mathfrak {M}}$ на підмножини, що
об'єднують функції $\psi \in {\mathfrak {M}}$, які в певному сенсі
мають однакову швидкість  прямування до нуля.

В ролі характеристики, за допомогою якої можна здійснити  таке розбиття, О.\,І.~Степанець обрав
пару функцій $\eta (t)=\eta (\psi ;t)$ і  $\mu (t)=\mu (\psi ;t)$, що означаються в такий спосіб.
Нехай  $\psi \in {\mathfrak {M}}$, тоді через $\eta
(t)=\eta (\psi ;t)$ позначають функцію, яка пов'язана з $\psi $
рівністю
  \begin{equation}\label{Class_1}%$$
\psi (\eta (t)) = \frac {1}2\psi (t), \ \ t\ge 1.
 \end{equation} %\eqno(1.4)$$
Внаслідок строгої монотонності функції $\psi$, характеристика $\eta (t)$ для всіх
$t\ge 1 $ з (\ref{Class_1}) визначається однозначно:
$
\eta (t) =\eta (\psi ;t) = \psi ^{-1} (\frac 1{2}\psi (t)).
 $
Функція $\mu (t)$ задається рівністю
$$
\mu (t) = \mu (\psi ;t)=\frac t{\eta (t) - t}. %\eqno (1.3)
 $$
 В залежності від поведінки функції $\mu$  розрізняють такі підмножини множини ${\mathfrak M}$:
 $$
 {\mathfrak M}_0=\{\psi\in {\mathfrak M}\ :\
 \ \ 0<\mu(\psi;t)\le K\ \ \forall t\ge 1\ \ \}\footnote{$K$, $K_1$, $K_2,\ldots$ -- деякі додатні сталі, що не
залежать від параметра $t$.},
 $$
 $$
 {\mathfrak M}_\infty=\{\psi\in {\mathfrak M}\ :\
 \ \ 0<K\le \mu(\psi;t)<\infty\ \ \forall t\ge 1\ \ \},
$$
$$
  {\mathfrak M}_C{=}{\mathfrak M}_0\cap  {\mathfrak M}_\infty{=}\{\psi\in {\mathfrak M}:\, 0<K_1 \le  \mu(\psi;t) \le  K_2 \ \forall t\ge 1\}.
$$

Через $B$ позначимо множину всіх монотонно спадних до нуля при $t\to\infty$ функцій $\psi(t)$, $t\ge 1$, які задовольняють так звану $\Delta_2$-умову:
\begin{equation}\label{5.m9}
% \frac{\psi(t)}{\psi(2t)}\le K.
{\psi(t)}\le K {\psi(2t)}.
\end{equation}
Як показано в \cite[Гл.3, \S 3.16]{Stepanets_M2002_1} має місце рівність
\begin{equation}\label{5.m10}
B\cap {\mathfrak M}={\mathfrak M}_0.
\end{equation}

Зазначимо, що природними представниками множин ${\mathfrak {M}}_C$ є функції  $t^{-r}$ при $r > 0$, а також функції
$t^{-r}\ln^{\varepsilon }(t{+}a)$ при довільних $\varepsilon \in {\mathbb R}$, додатних $r$ і $a$, для яких
$a\ge e^{3\varepsilon/r}-1$ та ін. До множини ${\mathfrak {M}}_0$ належать також функції
$\ln ^{-r}(t+a)$ при довільних  додатних  $r$ і $a$.

Через ${\mathfrak M}^{+}_\infty$ позначають підмножину   всіх функцій
$\psi\in {\mathfrak M}$, для яких  $\mu(\psi;t)$ монотонно і
необмежено зростає при $t\to \infty$:
 $$
{\mathfrak M}^{+}_\infty=\{\psi\in {\mathfrak M}\ :\
 \ \ \mu(\psi;t)\uparrow \infty\ \ \}.
 $$
З цієї множини виділяють такі підмножини:
 \[
{\mathfrak M}\,'_\infty=\{\psi\in {\mathfrak M}_\infty^+\ :\quad
\alpha(\psi;t)\downarrow 0,\quad
 {\psi(t)}/{|\psi'(t)|}\uparrow \infty\ \},
 \]
 де
 \[
  \alpha(\psi;t)=\frac{\psi(t)}{t|\psi'(t)|},\ \
 \psi'(t):=\psi'(t+),
 \]
 \[
{\mathfrak M}^c_{\infty}=\{\psi\in {\mathfrak M}_\infty^+\ :\quad
\alpha(\psi;t)\downarrow 0,\
 0<K_1<{\psi(t)}/{|\psi'(t)|}<K_2 \}
 \]
 і
 \[
{\mathfrak M}''_\infty=\{\psi\in {\mathfrak M}_\infty^+\ :\quad
 {\psi(t)}/{|\psi'(t)|}\downarrow 0 \ \}.
  \]
Природними представниками множин ${\mathfrak
 M}\,'_\infty$ та ${\mathfrak M}^c_{\infty}$ є  функції ${\rm exp}(-\lambda t^s)$, $\lambda>0$,
при $s\in (0,1)$ та  $s=1$ відповідно. До множини ${\mathfrak M}''_\infty$ належать функції
 $\exp (-\lambda (t+a)^r)$  при $\lambda >0$, $r>1$ і ${a\ge((r-1)/(r\lambda))^{1/r}-1}$.

\vskip 2mm
\noindent{\bf \ref{Order_Estimates}.2. Порядкові оцінки найкращих $n$-членних наближень та  поперечників  $q$-еліпсоїдів в просторах ${\mathcal S}^p_\varphi$.}
 Порядкові оцінки найкращих $n$-членних наближень $q$-еліпсоїдів в просторах ${\mathcal S}^p_\varphi$ містяться в такому твердженні.

\begin{theorem}[{\cite{Stepanets_Shidlich_Pr_2007, Stepanets_Shidlich_NK2007, Shidlich_2016}}]\label{Order_Best_n_term_approx}
    Нехай   $0<p,\,q<\infty$,   система чисел $\psi=\{\psi_k\}_{k=1}^\infty$ при всіх $k\in {\mathbb N}$ задовольняє рівність $\widetilde{\psi}_k=\psi_1(k)$, де  $\widetilde{\psi}=\{\widetilde{\psi}_k\}_{k=1}^\infty$ -- спадна перестановка послідовності чисел
$\{|\psi_k|\}_{k=1}^\infty$, а $\psi_1$ -- деяка додатна функція.

    $1)$ Якщо  функція $\psi_1^p$ належить множині $B$, а при $0<p<q$,    крім цього, при всіх $t$, більших
    деякого числа $t_0$,   є опуклою  та задовольняє умову
     \begin{equation}\label{1.2c121}%$$%
         t|\psi'_1(t)|/{\psi_1(t)}\ge K_0>\beta,
     \end{equation}
    де $\psi'_1(t):=\psi'_1(t+)$, $\beta=d(1/p-1/q)$, то
     \[
       e_n(\psi\,U_\varphi^q)_{p}\asymp  {\psi_1(n+1)}{\,n^{\frac 1p\,-\frac 1q}}.\footnote{Для додатних послідовностей $a(n)$ та $b(n)$ вираз "$a(n)\asymp b(n)$"\  означає, що  існують такі сталі $K_1,K_2>0$, що при всіх $n\in\mathbb{N}$ виконуються нерівності $\ a(n)\le K_2 b(n)$
       і ${a(n)\ge K_1 b(n)}$.}
     \]
     $2)$  Якщо функція  $\psi_1^p \in {\mathfrak M}'_\infty$, то
     \[
       e_n(\psi U_\varphi^q)_{p}\asymp %\frac
       {\psi_1(n+1)}{(\eta(\psi_1,n)-n)^{\frac 1p-\frac 1q}}.
     \]
     $3)$ Якщо функція $\psi_1^p$ належить множині ${\mathfrak M}^c_\infty$ або ${\mathfrak M}\,''_\infty$, то
      \[
        e_n(\psi U_\varphi^q)_{p}\asymp  {\psi_1(n+1)},\quad n\to\infty.
     \]
 \end{theorem}

Наведемо   порядкові оцінки поперечників ${\mathscr D}_n(\psi\,U_\varphi^q)_{p}$ та ${\mathscr D}_n^\perp(\psi\,U_\varphi^q)_{p}$ при $0<p<q<\infty$.

\begin{theorem}[{\cite{Stepanets_Shidlich_Pr_2007, Shidlich_2016}}]\label{Order_Basis_Width}
    Нехай   $0<p<q<\infty$,   система чисел $\psi=\{\psi_k\}_{k=1}^\infty$ при всіх $k\in {\mathbb N}$ задовольняє рівність $\widetilde{\psi}_k=\psi_1(k)$,
     де  $\widetilde{\psi}=\{\widetilde{\psi}_k\}_{k=1}^\infty$ -- спадна перестановка послідовності чисел $\{|\psi_k|\}_{k=1}^\infty$, а $\psi_1$ -- деяка додатна  функція.

    $1)$ Якщо  функція $\psi_1^p$ належить множині $B$, при всіх $t$, більших
    деякого числа $t_0$,   є опуклою  та задовольняє умову (\ref{1.2c121}) з  $\beta=d(1/p-1/q)$, то
     \[
       {\mathscr D}_n(\psi\,U_\varphi^q)_{p}\asymp
       {\mathscr D}_n^\perp(\psi\,U_\varphi^q)_{p}\asymp  {\psi_1(n+1)}{\,n^{\frac 1p\,-\frac 1q}}.
     \]
     $2)$  Якщо функція  $\psi_1^p \in {\mathfrak M}'_\infty$, то
     \[
       {\mathscr D}_n(\psi\,U_\varphi^q)_{p}\asymp
       {\mathscr D}_n^\perp(\psi\,U_\varphi^q)_{p}\asymp   {\psi_1(n+1)}{(\eta(\psi_1,n)-n)^{\frac 1p-\frac 1q}}.
     \]
     $3)$ Якщо функція $\psi_1^p$ належить множині ${\mathfrak M}^c_\infty$ або ${\mathfrak M}\,''_\infty$, то
      \[
        {\mathscr D}_n(\psi\,U_\varphi^q)_{p}\asymp
       {\mathscr D}_n^\perp(\psi\,U_\varphi^q)_{p}\asymp  {\psi_1(n+1)}.
     \]
 \end{theorem}

Порівнюючи порядкові оцінки для величин $e_n(\psi
U_\varphi^q)_{p}$ та ${\mathscr D}_n(\psi\,U_\varphi^q)_{p}$ бачимо, що у випадку, коли $0<q\le
 p$, а послідовність $\psi=\{\psi_k\}_{k=1}^\infty$ така, що при всіх натуральних $k$
 виконується  рівність\ \
$\widetilde{\psi}_k=\psi_1(k)$, де функція $\psi_1$ задовольняє одну з умов 1) чи 2) теореми \ref{Order_Best_n_term_approx},
мають місце рівності
 $$
 \lim\limits_{n\to\infty}\frac{e_n(\psi
U_\varphi^q)_{p}}{{\mathscr
 D}_n(\psi\,U_\varphi^q)_{p}}= \lim\limits_{n\to\infty}\frac{e_n(\psi
U_\varphi^q)_{p}}{{\mathscr
 D}_n^\perp(\psi\,U_\varphi^q)_{p}}=0.
 $$
Якщо ж $0<p<q$ і  $\psi_1$ задовольняє одну з  умов 1) чи 2), або ж якщо $0<p,\,q<\infty$ і
$\psi_1$ належить до однієї з множин  ${\mathfrak M}^c_\infty$ чи ${\mathfrak M}\,''_\infty$, то
 $$
 {e_n(\psi
U_\varphi^q)_{p}}\asymp{{\mathscr
 D}_n(\psi\,U_\varphi^q)_{p}}\asymp{{\mathscr
 D}_n^\perp(\psi\,U_\varphi^q)_{p}}.
 $$

%%%%%%%%%%%%%%%%%%%%%%%%%%%%%%%%%%%%%%%%%%%%%%%%%%%%%%%%%%%%%%%%%%%%%%%%%%%%%%%%%%%%%%%%%%%%%%%%%%%%%%
\subsection{Найкращі  $n$-членні наближення з обмеженнями}\label{Best_Approx_Restrictions}

\noindent{\bf \ref{Best_Approx_Restrictions}.1. Найкращі  $n$ -членні наближення з обмеженнями $q$-еліпсо\-їдів в просторах
 ${\mathcal S}^p_\varphi$ при $0<q\le p$.}
Нехай $\Gamma_n$  -- множина всіх наборів $\gamma_n$ з $n$ різних натуральних чисел.
В такому випадку величину $e_n(f)_{p}$, означена рівністю
(\ref{b95}), можна записати у  вигляді
  $$
  e_n(f)_p = \inf\limits_{\gamma_n \in \Gamma_n}
   E_{\gamma_n}(f)_{p}.
   $$
Поряд з  $e_n(f)_p$  можна розглядати і величини
 \begin{equation}\label{b165}%$$
     e_n(f;\Gamma^{\,'}_{n})_p = \inf\limits_{\gamma_n \in\, \Gamma_{n}^{\,'}}E_{\gamma_n}(f)_p,
  \end{equation}
де $\Gamma_{n}^{\,'}$ -- деяка   підмножина з $\Gamma_n$.  У зв'язку з цим  величину $e_n(f)_p$ зручно назвати
абсолютним найкращим $n$-членним наближенням, а величину $e_n(f;\Gamma^{\,'}_n)_p$ -- найкращим $n$-членним
наближенням з обмеженнями, маючи на увазі, що тут термін ``обмеження''  стосується вибору підмножини  $\Gamma^{\,'}_n$.

В ролі $\Gamma^{\,'}_n$ розглянемо дві підмножини $\Gamma^{(1)}_n$ та $\Gamma^{(2)}_n$. Через $\Gamma^{(1)}_n$
позначають множину наборів
$$
  \gamma^{(1)}_n = \{ j\,n+1,\;j\,n+2,\ldots (j+1)n
  \},\;j=0,1,\ldots;
$$ а через $\Gamma^{(2)}_n$ -- множину наборів
  $$
  \gamma^{(2)}_n = \{ j+1,\;j+2,\ldots j+n \},\;j=0,1,\ldots .
  $$
Зрозуміло, що завжди  $\Gamma^{(1)}_n \subset \Gamma^{(2)}_n \subset \Gamma_n$  і
 \begin{equation}\label{b166}%$$
  e_n(f)_p \leq e_n(f;\Gamma^{(2)}_n)_p \leq
  e_n(f;\Gamma^{(1)}_n)_p.
 %\eqno (4.26) $$
  \end{equation}
Тому якщо   $\mathfrak N$ -- деяка підмножина з $ {\mathcal S}^p_\varphi$ і
 \begin{equation}\label{b167}%$$
  e_n(\mathfrak N;\Gamma^{\,'}_n)_p = \sup\limits_{f \in\, \mathfrak
  N}e_n(f;\Gamma^{\,'}_n)_p\,,
 %\eqno (4.26) $$
  \end{equation}
то мають місце нерівності
 \begin{equation}\label{b168}%$$
  e_n(\mathfrak N)_p \leq e_n(\mathfrak N;\Gamma^{(2)}_n)_p \leq
  e_n(\mathfrak N;\Gamma^{(1)}_n)_p.
   %\eqno (4.26) $$
  \end{equation}
Як і раніше, в ролі множин $\mathfrak N$ вибираємо множини   $\psi\,U_\varphi^q$, які задаються рівністю (\ref{b16}) при $p=q$.

\begin{theorem}[\cite{Stepanets_Rukasov_UMZh2003_5, Stepanets_UMZh2005_4, Stepanets_UMZh2006_1}]\label{Best_Approx_Restr_q<p}
  Нехай $0<q\le p$ і $\psi=\{\psi_k\}_{k=1}^\infty$ --
 система комплексних чисел, для яких послідовність $|\psi_k|$, $k=1,2,\ldots,$ монотонно прямує до нуля.
 Тоді  при довільному $n\in {\mathbb N}$ виконуються рівності
 \[
  e_n(\psi \,U_\varphi^q;\Gamma^{(1)})_p = e_n(\psi\,U_
  \varphi^q;\Gamma^{(2)})_{p} = \frac{(s^{\ast}-1)^{1/p}}{ \Big( \sum\limits_{k=1}^{s^{\ast}} |\,{\psi}_{(k-1)n+1}|^
  {-q} \Big)^{{1/q}}},
 \]
 де $s^{\ast}$ -- деяке натуральне число, для якого
 \[
 \sup\limits_{s>1} \frac{(s-1)^{1/p}}{ \Big(\sum\limits_{k=1}^{s}|{\psi}
   _{(k-1)n+1}|^{-q} \Big)^{{1/q}}} = \frac{ (s^{\ast}-1)^{1/p}}{
   \Big(\sum\limits_{k=1}^{s^{\ast}} |{\psi}_{(k-1)n+1}|^
  {-q} \Big)^{{1/q}}}.
 \]
 Таке число $s^{\ast}$ завжди  існує.
 \end{theorem}

%%%%%%%%%%%%%%%%%%%%%%%%%%%%%%%%%%%%%%%%%%%%%%%%%%%%%%%%%%%%%%%%%%%%%%%%%%%%%
\vskip 2mm
\noindent{\bf \ref{Best_Approx_Restrictions}.2. Найкращі  $n$ -членні наближення з обмеженнями $q$-еліпсоїдів в просторах
 ${\mathcal S}^p_\varphi$ при $0<p<q$.} У випадку, коли  $0<p<q$ точні значення величин $e_n(\psi \,U_\varphi^q;\Gamma^{(1)})_p $ містяться в такому твердженні.

\begin{theorem}[\cite{Stepanets_UMZh2005_4, Stepanets_UMZh2006_1}]\label{Best_Approx_Restr_q>p} Нехай  $0<p<q$
і  $\psi=\{\psi_k\}_{k=1}^\infty$ -- система  таких комплексних чисел, що виконуються умови (\ref{b15})
та (\ref{b54}) і  послідовність $\{|\psi_k|\}_{k=1}^\infty$ не зростаючи прямує до нуля.
Тоді при довільному $n\in {\mathbb N}$
\[
  e_n(\psi\,U_\varphi^{q}; \Gamma^{(1)}_n)_{p} =
  \bigg(\!(s^*-1)^\frac q{q-p} \Big (\sum\limits_{k=1}^{s^*}
  \widetilde{\psi}_k^{-q}
  \Big)^ {-\frac p{q-p}}\!\!
 +\!\!\!
 \sum\limits_{k=s^*n+1}^\infty  \widetilde{\psi}_k^\frac
 {pq}{q-p}\bigg)^\frac{q-p}{pq}\!\!,
\]
де
 \[
\widetilde{\psi}_k=
\Big(\sum\limits^{k\,n}_{i={(k-1)n+1}}|\psi_i|^{\frac{p\,q}{q-p}}
\Big)^{\frac{q-p}{pq}},\quad k=1,2,\ldots,
  \]
 Число $s$ вибране з умови
 \[
 \widetilde{\psi}_{s^*}^{-q}\le
 \frac 1{s^*-1}\sum\limits_{k=1}^{s^*}\widetilde{\psi}_k^{-q}
 <\widetilde{\psi}_{s^*+1}^{-q}.
   \]
Таке число $s$ завжди існує і єдине.
 \end{theorem}

%%%%%%%%%%%%%%%%%%%%%%%%%%%%%%%%%%%%%%%%%%%%%%%%%%%%%%%%%%%%%%%%%%%%%%%%%
Для величин $e_n(\psi U_\varphi^q;\Gamma_n^{(2)})_{p}$ при $0<p<q$,  взагалі кажучи, має місце  нерівність
 \[
 e_n(\psi
U_\varphi^q;\Gamma_n^{(2)})_{p} \leq \bigg((s^*-1)^\frac
q{q-p} \Big(\sum\limits_{k=1}^{s^*}
  \widetilde{\psi}_k^{-q}
  \Big)^ {-\frac p{q-p}}+
 \sum\limits_{k=s^*n+1}^\infty  \widetilde{\psi}_k^\frac
 {pq}{q-p}\bigg)^\frac{q-p}{pq}.
\]
Більш детально з цим випадком і  зокрема, з умовами, за яких в останньому співвідношенні має місце рівність, можна ознайомитись у роботі  \cite{Stepanets_UMZh2005_4}.

%%%%%%%%%%%%%%%%%%%%%%%%%%%%%%%%%%%%%%%%%%%%%%%%%%%%%%%%%%%%%%%%%%%%%%%%%%%%%%%%%%%%%%%%%%%%%%%%%%%%%%
\section{Наближення в просторах ${\mathcal S}^p$.}\label{Sp_spaces}

%%%%%%%%%%%%%%%%%%%%%%%%%%%%%%%%%%%%%%%%%%%%%%%%%%%%%%%%%%%%%%%%%%%%%%%%%%%%%
\subsection{Основні означення.}\label{Definitions_Sp} Нехай, як і в п.~\ref{Examples},  $L=L({\mathbb T}^d)$, $d\ge 1$,  -- множина всіх $2\pi$-періодичних за кожною зі
змінних функцій $f({\bf x})=f(x_1,\cdots ,x_d)$, сумовних на кубі періодів ${\mathbb T}^d$ і (\ref{b8}) -- ряд Фур'є функції
$f\in L$ за системою (\ref{b9}). Еквівалентні відносно міри Лебега функції ототожнюються.

Нехай, далі, ${\mathcal S}^p $ -- простір, породжений множиною $L$, системою (\ref{b9}) і деяким числом  $ p\in (0,
\infty )$, зі скалярним добутком (\ref{b10}) і (квазі-)нормою $\|\cdot \|_p=\|\cdot \|_{_{\scriptstyle {\mathcal S}^p}}$, визначеною згідно з (\ref{Def_Sp}):
   \begin{equation}\label{w1}%$$
   \|f\|_{_{\scriptstyle {\mathcal S}^p}} = \Big(\sum _{{\bf k}\in {\mathbb Z}^d}|\widehat f({\bf  k})|^p\Big)^{1/p}. %\eqno (9.1) $$
 \end{equation}

Нехай тепер $ {\psi }=\{\psi ({\bf k})\}_{{\bf k}\in {\mathbb Z}^d}$ -- довільна система комплексних чисел -- кратна
послідовність. Якщо для  функції $f\in L$  з рядом Фур'є (\ref{b8}) ряд
  \begin{equation}\label{w2}%$$
\sum _{{\bf k}\in {\mathbb Z}^d}\psi ({\bf k})\widehat {f} ({\bf
k}){\mathrm e}^{{\mathrm i}({\bf k},{\bf x})}. %\eqno (9.2) $$
 \end{equation}
є рядом Фур'є  деякої функції  $F$ з  $L,$ то  $F$ називають $\psi$-інтегралом функції  $f$ і позначають
$F={\mathcal J}^{\psi } (f).$  При цьому функцію  $f$ називають $\psi $-похідною функції  $F$ і позначають  $f=D^{\psi }(F)=F^{\psi }$.
Множину $\psi$-інтегралів всіх функцій   $f\in L$ позначають через $L^{\psi}.$ Якщо  ${\mathfrak {N}}$ -- деяка підмножина з $L,$ то
через $L^{\psi }{\mathfrak {N}}$ позначають  множину  $\psi$-інтегралів всіх функцій з ${\mathfrak {N}}.$
Зрозуміло, що коли $f\in L^{\psi },$  коефіцієнти Фур'є функцій $f$ и $f^{\psi }$ пов'язані співвідношеннями
 \begin{equation}\label{w3}%$$
  \widehat f({\bf k})=\psi ({\bf k})\widehat f^{\psi }({\bf k}), \ \ {\bf k}\in {\mathbb Z}^d. %\eqno (9.3) $$
 \end{equation}
В ролі ${\mathfrak {N}}$ можна обрати одиничну кулю в  $U^p$ в просторі ${\mathcal S}^p $:
\begin{equation}\label{w4}%$$
 U^p=\{f\in {\mathcal S}^p: \quad  \|f\|_p\le 1\}. %\eqno (9.4) $$
 \end{equation}
В такому випадку покладаємо  $L_p^{\psi }:=L_p^{\psi }({\mathbb T}^d)=L^{\psi }U^p.$ Система  $\psi$, як і вище, підпорядкована умові
\begin{equation}\label{w5}%$$
 \lim _{|{\bf k}|\to \infty }\psi ({\bf k})=0. %\eqno (9.5) $$
  \end{equation}
Зазначимо, що коли  $f\in L^{\psi }{\mathcal S}^p$ і  $|\psi ({\bf k})|\le K,$ ${\bf k}\in {\mathbb Z}^d,$  то
$f\in {\mathcal S}^p$. Тому за  умови (\ref{w5}) має місце  вкладення   $L_p^{\psi }\subset {\mathcal S}^p.$

Означимо характеристичні послідовності $\varepsilon(\psi),$ $ g(\psi)$ і $\delta (\psi )$ аналогічно до того як це зроблено у підрозділі \ref{Charteristic}
Через $\varepsilon(\psi )=\varepsilon_1, \varepsilon_2, \ldots $ позначимо множину значень величин $|\psi ({\bf k})|,$ ${\bf k}\in
{\mathbb Z}^d,$ впорядковану за спаданням. Розглянемо також послідовності $g(\psi )=\{g_n\}_{n=1}^{\infty }$ та  $\delta (\psi )=\{\delta _n\}_{n=1}^{\infty },$
 де  $g_n=g_n(\psi)=\{{\bf k}\in {\mathbb Z}^d:  |\psi ({\bf k})|\ge \varepsilon_n\}$ і  $\delta_n=\delta _n^{\psi }=|g_n|$ -- кількість елементів ${\bf k}\in
{\mathbb Z}^d$, що належать множині  $g_n.$

З огляду на умову  (\ref{w5}), в даному  випадку послідовності   $\varepsilon(\psi)$ та $g(\psi )$ означаються
рівностями (\ref{b17i}) з врахуванням того, що ${\bf k}\in {\mathbb Z}^d$. Як і раніше, вважаємо, що $g_0=g_0 (\psi )$ -- порожня множина  і  $\delta _0:=\delta _0 (\psi )=0$.

Зазначимо, що крім природної умови  (\ref{w5}) від системи $\psi $ жодних інших істотних обмежень не вимагатиметься.
Тому ці системи  $\psi $, а з ними і їх характеристичні послідовності $\varepsilon(\psi ),$
$g(\psi )$ та $\delta (\psi )$ в загальному випадку можуть бути різноманітними та достатньо складними.

В багатовимірному випадку, напевно, найбільш простими і природними є
системи $\psi,$ у яких    $\psi({\bf k})$ зображуються добутками
\begin{equation}\label{w35}%$$
 \psi ({\bf {k}})=\psi (k_1, \ldots ,k_d)=\prod _{j=1}^d\psi
 _j({k}_j), \ \ {k_j}\in {\mathbb Z}^1, \ \ j=\overline {1,d}, %\eqno (10.1)$$
  \end{equation}
значень одновимірних послідовностей $\psi _j=\{\psi
_j(k_j)\}_{k_j=1}^{\infty }.$ Якщо при цьому $\psi(-k_j)=\overline {\psi _j(k_j)}$, $j=\overline {1,d}$,
то множини $g_n(\psi)$ будуть симетричними відносно усіх координатних площин і, як неважно переконатися,
 \begin{equation}\label{w36}%$$
 \sum _{{\bf k}\in {\mathbb Z}^d}\psi ({\bf k})e^{i{\bf
 k}{\bf
 t}}=\sum _{{\bf k}\in {\mathbb Z}_{+}^d}2^{d- q({\bf
 k})}\prod _{j=1}^d |\psi _j(k_j)|\cos \Big(k_jt_j - \frac {\beta
 _{k_j}\pi }2\Big), %\eqno (10.2) $$
   \end{equation}
де  $q({\bf k})$ -- кількість координат вектора ${\bf k}$, які дорівнюють  нулю, а
числа $\beta _{k_j}$ означаються рівностями
 $$
 \cos \frac {\beta _{k_j}\pi }2=\frac {{\rm Re} \ \psi _j(k_j)}{|\psi
 _j(k_j)|}, \ \ \ \sin \frac {\beta _{k_j}\pi }2=\frac {{\rm Im }\
 \psi _j(k_j)}{|\psi _j(k_j)|}.
 $$
 В такому випадку множина $L^{\psi}$ $\psi$-інтегралів дійсних функцій
  $\varphi $ з $L({\mathbb T}^d)$ складається з дійсних функцій $f$, і якщо при цьому ряд в
  (\ref{w2}) є рядом Фур'є деякої сумовної  функції ${\mathscr{D}}_{\psi }(t),$ то
  достатньою умовою включення  $f\in L^{\psi }{\mathfrak {N}}$
  є можливість зображення функції  $f$ у вигляді згортки
  $$
  f(x)=(2\pi )^{-d}\int\limits _{{\mathbb T}^{d}}\varphi
  (x-t){\mathscr{D}}_{\psi }(t){\mathrm d}t,
  $$
де $\varphi \in {\mathfrak {N}}$ і майже скрізь  $\varphi ({\bf x})=f^{\psi }({\bf x})$.
Це, зокрема, означає, що множини $L^{\psi}{\mathfrak {N}}$ містять класи функцій, які зображуються у вигляді згорток
з фіксованими сумовними ядрами.

\subsection{Застосування отриманих результатів до задач наближення періодичних функцій багатьох змінних.}\label{Applications}

Наведемо застосування результатів попередніх підрозділів до розв'язання задач теорії наближення функцій багатьох змінних в просторах ${\mathcal S}^p$.

\vskip 2mm
\noindent{\bf \ref{Applications}.1. Найкращі наближення, найкращі $n$-членні наближення та базисні поперечники
класів $L_q^{\psi }$ в просторах ${\mathcal S}^p$.} Нехай $f$ -- довільна функція з простору ${\mathcal S}^p$, $n$ -- будь-яке натуральне число і
$\gamma_n$ -- довільний набір з $n$ векторів ${\bf k}\in {\mathbb Z}^d$. Розглянемо тригонометричні поліноми
 \begin{equation}\label{w14a}%  $$
    P_{\gamma_n} = \sum\limits_{{\bf k}\in \gamma_n}
c_{\bf k} {\mathrm e}^{{\mathrm i}({\bf k},\cdot)}\quad \mbox{\rm та}\quad S_{\gamma_n}(f)=
  \sum\limits_{{\bf k}\in \gamma_n}\widehat{f}({\bf k}){\mathrm e}^{{\mathrm i}({\bf k},\cdot)},
 \end{equation}
де $c_{\bf k}$ -- будь-які комплексні числа, а $\widehat{f}({\bf k})$ -- коефіцієнти Фур'є функції $f$, а також величини
 \begin{equation}\label{w14}%  $$
  E_{\gamma_n}(f)_{_{\scriptstyle {\mathcal S}^p}} \!=\!\!\inf_{c_{\bf k}\in {\mathbb C}}
  \|f - P_{\gamma_n}\|_{_{\scriptstyle {\mathcal S}^p}}\quad \mbox{\rm і}\quad
  {\mathscr E}_{\gamma_n}(f)_{_{\scriptstyle {\mathcal S}^p}}\!=\|f -\! S_{\gamma_n}(f) \|_{_{\scriptstyle {\mathcal S}^p}}.
 \end{equation}
Якщо ${\mathfrak N}$ -- деяка підмножина з ${\mathcal S}^p$, то покладаємо
  \begin{equation}\label{new15}%  $$
   E_{\gamma_n}({\mathfrak N})_{_{\scriptstyle {\mathcal S}^p}} = \sup\limits_{f\in {\mathfrak N}}
   E_{\gamma_n}(f)_{_{\scriptstyle {\mathcal S}^p}}\quad \mbox{\rm і}\quad
   {\mathscr E}_{\gamma_n}({\mathfrak N})_{_{\scriptstyle {\mathcal S}^p}} = \sup\limits_{f\in {\mathfrak N}}
   {\mathscr E}_{\gamma_n}(f)_{_{\scriptstyle {\mathcal S}^p}}
 \end{equation}
а також
\begin{equation}\label{new16}%  $$
   {\mathscr D}_n({\mathfrak N})_{_{\scriptstyle {\mathcal S}^p}} = \inf\limits_{\gamma_n}
  E_{\gamma_n}({\mathfrak N})_{_{\scriptstyle {\mathcal S}^p}} \quad \mbox{\rm і}\quad
   {\mathscr D}_n^\perp({\mathfrak N})_{_{\scriptstyle {\mathcal S}^p}} = \inf\limits_{\gamma_n}
  {\mathscr E}_{\gamma_n}({\mathfrak N})_{_{\scriptstyle {\mathcal S}^p}}.
 \end{equation}
Через $e_n({\mathfrak N})_{_{\scriptstyle {\mathcal S}^p}}$ позначаємо  найкраще $n$-членне тригонометричне наближення
множини ${\mathfrak N}$ в просторі ${\mathcal S}^p$, тобто, величину
 \begin{equation}\label{w12}%$$
    e_n({\mathfrak N})_{_{\scriptstyle {\mathcal S}^p}}=\sup\limits _{f\in {\mathfrak N}}
    \inf_{\gamma _n} E_{\gamma_n}(f)_{_{\scriptstyle {\mathcal S}^p}}.
       \end{equation}
В ролі множини ${\mathfrak N}$ розглядаємо  множину  $L_q^{\psi }$   $\psi$-інтегралів всіх функцій з
одиничної кулі простору ${\mathcal S}^q$, $0<p,q<\infty$, за умов, що гарантують   вкладення   $L_q^{\psi }\subset {\mathcal S}^p$.

За таких позначень  мають місце такі твердження -- наслідки відповідних теорем підрозділів 3-6.

\begin{proposition}[{\cite[Гл.~11]{Stepanets_M2002_2}, \cite{Stepanets_UMZh2006_1}}]\label{Applic_q<p}
    Нехай   $0<q\le p$, і $ \psi =\{\psi ({\bf k})\}_{{\bf k}\in {\mathbb Z}^d} $ -- система чисел, підпорядкована умовам
    (\ref{w5})      і
     \begin{equation}\label{w17}%  $$
               \psi ({\bf k})\not =0\quad \forall {\bf k}\in  {\mathbb Z}^d.
     \end{equation}
 Тоді для будь-якого  $n \in {\mathbb N}$ і
    для довільного набору  $\gamma_n$ із $n$ різних натуральных чисел справджуються рівності
 \begin{equation}\label{w30}%  $$
        E_{\gamma_n}(L_q^{\psi })_{_{\scriptstyle {\mathcal S}^p}}=
        {\mathscr E}_{\gamma_n}(L_q^{\psi })_{_{\scriptstyle {\mathcal S}^p}}
        =\widetilde{\psi}_{\gamma_n}(1),
    \end{equation}
 де $\widetilde{\psi}_{\gamma_n}(1)$ -- перший  член послідовності
 $\widetilde{\psi}_{\gamma_n}=\{\widetilde{\psi}_{\gamma_n}(k) \}^{\infty}_{k=1}$, яка є спадною перестановкою
 системи чисел $\{|\psi_{\gamma_n}({\bf k})|\}_{{\bf k}\in {\mathbb Z}^d}$,
\begin{equation}\label{w31}%$$
 \psi_{\gamma_n}({\bf k})=\left\{\begin{matrix}0,\quad\hfill & {\bf k}\in\gamma_n,\\  \psi({\bf k}),\quad\hfill & {\bf k}\overline{\in}
 \gamma_n\end{matrix};\right.
 %\eqno (9.25) $$
 \end{equation}
при всіх $n \in {\mathbb N}$ виконуються рівності
\begin{equation}\label{w32}%$$
   {\mathscr D}_n(L_q^{\psi })_{_{\scriptstyle {\mathcal S}^p}} =
   {\mathscr D}_n^\perp(L_q^{\psi })_{_{\scriptstyle {\mathcal S}^p}}=
   \widetilde{\psi}_{n+1},
 \end{equation}
\begin{equation}\label{w19}%  $$
 e_n^p(L_q^{\psi})_{_{\scriptstyle {\mathcal S}^p}}=\mathop {\sup }\limits
 _{s>n}(s-n)(\sum _{k=1}^s\bar \psi _k^{-q})^{-\frac
 p{q}}=(s^*-n)(\sum _{k=1}^{s^*}\bar \psi_k^{-q})^{-\frac p{q}},
 %\eqno (9.15) $$
  \end{equation}
 в яких $ \widetilde{\psi} =\{\widetilde{\psi} _k\}_{k=1}^{\infty } $ -- спадна перестановка системи чисел $|\psi_{\bf k}|$, ${{\bf k}\in {\mathbb Z}^d}$, і  $s^*$ -- деяке натуральне число.
\end{proposition}

\begin{proposition}[{\cite[Гл.~11]{Stepanets_M2002_2}, \cite{Stepanets_UMZh2006_1}}]\label{Applic_q>p}
    Нехай   $0<p<q$, і $ \psi =\{\psi ({\bf k})\}_{{\bf k}\in {\mathbb Z}^d} $ -- система чисел,
    яка задовольняє умови  (\ref{w5}) та
    \begin{equation}\label{w21}%  $$
      \sum _{{\bf k}\in {\mathbb Z}^d}|\psi ({\bf k})|^{\frac {pq}{q-p}}<\infty.
    \end{equation}
Тоді для будь-якого  $n \in {\mathbb N}$ і  для довільного набору  $\gamma_n$ із $n$ різних натуральных чисел справджуються рівності
 \begin{equation}\label{w33}%  $$
   E_{\gamma_n}(L_q^{\psi })_{_{\scriptstyle {\mathcal S}^p}}=
  {\mathscr E}_{\gamma_n}(L_q^{\psi })_{_{\scriptstyle {\mathcal S}^p}} =\Big( \sum_{k=1}^{\infty} \,(\bar
 {\psi}_{\gamma_n}(k))^{\frac{p\,q}{q-p}}\Big)^{\frac{q-p}{p\,q}},
 %   \eqno (4.22) $$
 \end{equation}
де  $\widetilde{\psi}_{\gamma_n}=\{\widetilde{\psi}_{\gamma_n}(k) \}^{\infty}_{k=1}$ -- спадна перестановка
 системи чисел $\{|\psi_{\gamma_n}({\bf k})|\}_{{\bf k}\in {\mathbb Z}^d}$;  при всіх $n \in {\mathbb N}$ виконуються рівності
\begin{equation}\label{w34}%$$
   {\mathscr D}_n(L_q^{\psi })_{_{\scriptstyle {\mathcal S}^p}} =
   {\mathscr D}_n^\perp(L_q^{\psi })_{_{\scriptstyle {\mathcal S}^p}}=
   \Big( \sum_{k=n+1}^{\infty} \widetilde{\psi}_k^{\frac{p\,q}{q-p}}
  \Big)^{\frac{q-p}{p\,q}},
 \end{equation}
\begin{equation}\label{w23}%  $$
e_n^p(L_q^{\psi })_{_{\scriptstyle {\mathcal S}^p}}\!=\!\bigg((s^*-n)^\frac q{q-p} \Big(\sum\limits_{k=1}^{s^*} \widetilde{\psi}_k^{-q}
  \Big)^ {-\frac p{q-p}}\!\!+\!\!
 \sum\limits_{k=s^*+1}^\infty  \widetilde{\psi}_k^\frac
 {pq}{q-p}\bigg)^\frac{q-p}{pq}, %\eqno (9.19) $$,
  \end{equation}
 в яких $ \widetilde{\psi} =\{\widetilde{\psi} _k\}_{k=1}^{\infty } $ -- спадна перестановка системи чисел $|\psi_{\bf k}|$, ${{\bf k}\in {\mathbb Z}^d}$ і  $s^*$ -- деяке натуральне число.

\end{proposition}

\vskip 2mm
\noindent{\bf \ref{Applications}.2. Наближення  поліномами, побудованими по областях $g_n(\psi)$, в просторах ${\mathcal S}^p$.} Нехай  $\psi =\{\psi ({\bf k})\}_{{\bf k}\in {\mathbb Z}^d} $ -- система чисел, підпорядкована умовам  (\ref{w5}) та (\ref{w17}),
і функція $f$ належить множині $L^{\psi }{\mathcal S}^p$. Розглянемо випадок,
апроксимуючі поліноми будуються по наборах $\gamma_n$, які визначаються через  елементи $g_n(\psi)$
характеристичної послідовності  $g(\psi )$  системи $ \psi$.  Тоді поліноми (\ref{w14a}) набуватимуть вигляду $ P_{g_n(\psi)}=\sum_{{\bf k}\in g_n(\psi)}
c_{\bf k} {\mathrm e}^{{\mathrm i}({\bf k},\cdot)}$ і
\begin{equation}\label{w6}%$$
 S_n(f)_{\psi}=S_{g_n(\psi)}(f) =\sum_{{\bf k}\in g_n(\psi)}\widehat f({\bf k}) {\mathrm e}^{{\mathrm i}({\bf k},{\bf \cdot})},
  \end{equation}
$S_0(f)_{\psi}:=0$, а величини (\ref{w14}) --
 \[%
 E_n(f)_{_{\scriptstyle \psi ,{\mathcal S}^p}}=\inf\limits _{a_{\bf k}\in {\mathbb C}}\Big\|f- P_{g_n(\psi)}\Big\|_{_{\scriptstyle {\mathcal S}^p}}
 \quad \mbox{\rm і}\quad
 \mathscr {E}_n(f)_{_{\scriptstyle \psi ,{\mathcal S}^p}}=\|f- S_{n-1}(f)_{\psi}\|_{_{\scriptstyle {\mathcal S}^p}}.
  %\eqno (9.9) $$
  \]
Якщо ${\mathfrak N}$ -- деяка підмножина з $L^{\psi }{\mathcal S}^p$, то  покладаємо
 \[
   E_{n}({\mathfrak N})_{_{\scriptstyle \psi ,{\mathcal S}^p}} = \sup\limits_{f\in {\mathfrak N}}
   E_{\gamma_n}(f)_{_{\scriptstyle \psi ,{\mathcal S}^p}}\quad \mbox{\rm і}\quad
   {\mathscr E}_{\gamma_n}({\mathfrak N})_{_{\scriptstyle \psi ,{\mathcal S}^p}} = \sup\limits_{f\in {\mathfrak N}}
   {\mathscr E}_{\gamma_n}(f)_{_{\scriptstyle \psi ,{\mathcal S}^p}}.
 \]

\begin{proposition}[{\cite[Гл.~11]{Stepanets_M2002_2}, \cite{Stepanets_UMZh2006_1}}]\label{Applic_Psi_q<p}
    Нехай  $0<q\le p$  і $ \psi =\{\psi ({\bf k})\}_{{\bf k}\in {\mathbb Z}^d} $ -- система чисел, підпорядкована умовам  (\ref{w5})
    та (\ref{w17}).  Тоді для будь-якого  $n \in {\mathbb N}$ справджуються рівності
\begin{equation}\label{w18}%  $$
 E_n(L_q^{\psi })_{_{\scriptstyle \psi ,{\mathcal S}^p}}={\mathscr{E}}_n(L_q^{\psi })_{_{\scriptstyle \psi ,{\mathcal S}^p}}=\varepsilon_n,
 %\eqno(9.14) $$
  \end{equation}
де  $\varepsilon_n$  -- члени характеристичної послідовності системи $\varepsilon(\psi)$.
\end{proposition}

\begin{proposition}[{\cite[Гл.~11]{Stepanets_M2002_2}, \cite{Stepanets_UMZh2006_1}}]\label{Applic_Psi_q>p}
    Нехай   $0<p<q$  і $ \psi =\{\psi ({\bf k})\}_{{\bf k}\in {\mathbb Z}^d} $ -- система чисел,
    яка задовольняє умови  (\ref{w5}) та (\ref{w21}). Тоді для будь-якого  $n \in {\mathbb N}$ справджуються рівності
\begin{equation}\label{w22}%  $$
E_n(L_q^{\psi})_{_{\scriptstyle \psi ,{\mathcal S}^p}}={\mathscr{E}}_n(L_q^{\psi })_{_{\scriptstyle \psi ,{\mathcal S}^p}}=
\bigg(\sum\limits_{k=\delta_{n-1}+1}^\infty
\widetilde{\psi}_k^{\frac{pq}{q-p}}\bigg)^{\frac{q-p}{pq}},
\end{equation}
в яких  $ \widetilde{\psi} =\{\widetilde{\psi} _k\}_{k=1}^{\infty } $ -- спадна перестановка системи чисел $|\psi_{\bf k}|$, ${{\bf k}\in {\mathbb Z}^d}$, а
 $\delta_n$ -- члени характеристичної послідовності $\delta(\psi)$.

\end{proposition}

Аналоги теорем \ref{Direct_Th} та \ref{In_Th} у просторах ${\mathcal S}^p$ формулюються в такий спосіб.

\begin{proposition}[{\cite[Гл.~11]{Stepanets_M2002_2}, \cite{Stepanets_UMZh2006_1}}]\label{Applic_Direct_Th}  Нехай $ f\in L_p^{\psi },$ $
p>0,$ і $ \psi =\{\psi (\bf k)\}_{{\bf k}\in {\mathbb Z}^d} $
-- система чисел, яка задовольняє умову  (\ref{w5}). Тоді ряд
 $$
 \sum _{k=1}^{\infty }(\varepsilon_k^p -
 \varepsilon_{k-1}^p)E_k^p(f^{\psi })_{_{\scriptstyle \psi ,{\mathcal S}^p}}
 $$
 збігається, і при кожному $ n\in {\mathbb N}$ має місце рівність
 $$
 E_n^p(f)_{_{\scriptstyle \psi ,{\mathcal S}^p}}=\varepsilon_n^pE_n^p(f^{\psi })_{\psi
 ,p} + \sum _{k=n+1}^{\infty
 }(\varepsilon_k^p-\varepsilon_{k-1}^p)E_k^p(f^{\psi })_{_{\scriptstyle \psi ,{\mathcal S}^p}},
 $$
 де  $ \varepsilon_k $ -- елементи характеристичної послідовності  $ \varepsilon(\psi )$.

 \end{proposition}

\begin{proposition}[{\cite[Гл.~11]{Stepanets_M2002_2}, \cite{Stepanets_UMZh2006_1}}]\label{Applic_In_Th}
 Нехай $ f\in {\mathcal S}^p\cap L^\psi$, $ p>0$, система $ \psi =\{\psi (\bf k)\}_{{\bf k}\in {\mathbb Z}^d}$
 задовольняє умову (\ref{w5})  і
 $$
 \lim _{k\to \infty }\varepsilon_k^{-1}E_k(f)_{_{\scriptstyle \psi ,{\mathcal S}^p}}=0.
 $$
Тоді для того, щоб виконувалося включення $ f\in L_p^{\psi }$  необхідно і достатньо, щоб збігався ряд
 $$
 \sum
 _{k=1}^{\infty }(\varepsilon_k^{-p} -
 \varepsilon_{k-1}^{-p})E_k^p(f)_{_{\scriptstyle \psi ,{\mathcal S}^p}}.
 $$
 Якщо цей ряд збігається, то при довільному  $n\in {\mathbb N},$ справджується рівність
 $$
 E_n^p(f^{\psi })_{\psi
 ,p}=\varepsilon_n^{-p}E_n^p(f)_{_{\scriptstyle \psi ,{\mathcal S}^p}} + \sum _{k=n+1}^{\infty
 }(\varepsilon_k^{-p}-\varepsilon_{k-1}^p)E_{k}^p(f)_{_{\scriptstyle \psi ,{\mathcal S}^p}},
 $$
 де  $ \varepsilon_k $ -- елементи характеристичної послідовності  $ \varepsilon(\psi )$.

 \end{proposition}

\vskip 2mm
\noindent{\bf \ref{Applications}.3. Поперечники за Колмогоровим класів $L_p^{\psi }$.}  Нехай
${\mathscr G}_n$ -- множина всіх $n$-вимірних підпросторів $G_n$  в ${\mathcal S}^p$, $n\in {\mathbb N}$, і
 \[
   d_n(L_p^{\psi })_{_{\scriptstyle {\mathcal S}^p}}= \inf\limits_{G_n\in\, {\mathscr G}_n} \sup\limits_{f\in L_p^{\psi }}
      \inf\limits _{u\in G_n}\|f-u \|_{_{\scriptstyle {\mathcal S}^p}}, \ \
      d_0(L_p^{\psi })_{_{\scriptstyle {\mathcal S}^p}}:=\sup \limits _{f\in L_p^{\psi }}\|f\|_{_{\scriptstyle {\mathcal S}^p}}, %\eqno (9.11)$$
 \]
   -- поперечники за  Колмогоровим класів $L_p^{\psi }$  в просторі ${\mathcal S}^p$.

\begin{proposition}[{\cite[Гл.~11]{Stepanets_M2002_2}, \cite{Stepanets_UMZh2006_1}}]\label{Applic_Kolmogorov_W} Нехай $ p\in [1,
\infty )$ і $ \psi =\{\psi ({\bf k})\}_{{\bf k}\in {\mathbb Z}^d} $ -- система чисел, підпорядкована умовам  (\ref{w5}) та (\ref{w17}). Тоді при довільних   $ n\in {\mathbb N}$ виконуються рівності
              $$
              d_{\delta _{n-1}}(L_p^{\psi })_{_{\scriptstyle {\mathcal S}^p}}=  d_{\delta _{{n-1}+1}}(L_p^{\psi })_{_{\scriptstyle {\mathcal S}^p}}=\ldots
              $$
         \begin{equation}\label{w29}%$$
          = d_{\delta _n-1}(L_p^{\psi })_{_{\scriptstyle {\mathcal S}^p}}= {\mathscr{E}}_n(L_p^{\psi })_{_{\scriptstyle \psi ,{\mathcal S}^p}}=\varepsilon_n,
         \end{equation}
 де $ \varepsilon_n$ та $ \delta _n $ -- члени характеристичних  послідовностей     $ \varepsilon(\psi )$ та $
\delta (\psi ),$ відповідно.
 \end{proposition}

\vskip 2mm
\noindent{\bf \ref{Applications}.4. Деякі наслідки для просторів $L_p({\mathbb T}^d)$.} Нехай $L_p=L_p({\mathbb T}^d),$ $p\in [1, \infty ),$ -- простір функцій $f\in L$ зі скінченною нормою  (\ref{Lp_norm}). Зв'язок між просторами  $L_p$ та ${\mathcal S}^p$ встановлює відома теорема Гаусдорфа--Юнга (див., наприклад, \cite{Temlyakov_B1993}),
з якої випливає, що при  $p\in (1,2]$ і $\frac 1p+\frac 1{p'}=1$ мають місце формули
\begin{equation}\label{w38}%$$
 L_p\subset {\mathcal S}^{p'}\ \ \ \ \mbox {і} \ \ \ \ \|f\|_{{\mathcal S}^{p'}}\le \|f\|_{L_{p}}, %\eqno (10.5) $$
  \end{equation}
\begin{equation}\label{w39}%$$
 {\mathcal S}^p\subset L_{p'}\ \ \ \ \mbox {і} \ \ \ \ \|f\|_{L_{p'}}\le \|f\|_{_{\scriptstyle {\mathcal S}^p}}. %\eqno (10.6) $$
  \end{equation}
Зокрема, при  $p=p'=2 $ виконуються рівності
 \begin{equation}\label{w40}%$$
 L_2={\mathcal S}^2 \ \ \ \ \mbox{і} \ \ \ \ \|\cdot\|_{L_2}=\|\cdot
 \|_{{\mathcal S}^2}. %\eqno (10.7)$$
   \end{equation}
Отже,  теореми, доведені для просторів ${\mathcal S}^p,$ містять певну інформацію і
для просторів $L_p,$ яка є найбільш повною внаслідок  (\ref{w40}), у випадку, коли $p=2.$

У останньому випадку, з теорем \ref{Applic_Psi_q<p} та \ref{Applic_Kolmogorov_W} отримуємо   наслідок

\begin{corollary}[{\cite[Гл.~11]{Stepanets_M2002_2}, \cite{Stepanets_UMZh2006_1}}]\label{Applic_Psi_p=2}
    Нехай  $ \psi =\{\psi ({\bf k})\}_{{\bf k}\in {\mathbb Z}^d} $ -- система чисел, підпорядкована умовам  (\ref{w5})
    та (\ref{w17}).  Тоді для будь-якого  $n \in {\mathbb N}$ справджуються рівності
             $$
              d_{\delta _{n-1}}(L_2^{\psi })_{_{\scriptstyle L_2}}=
              d_{\delta _{{n-1}+1}}(L_2^{\psi })_{_{\scriptstyle L_2}}=\ldots
              = d_{\delta _n-1}(L_2^{\psi })_{_{\scriptstyle L_2}}
              $$
         \begin{equation}\label{p=2}%$$
          =E_n(L_2^{\psi })_{_{\scriptstyle \psi ,L_2}}=
          {\mathscr{E}}_n(L_2^{\psi })_{_{\scriptstyle \psi ,L_2}}=\varepsilon_n,
         \end{equation}
 де $ \varepsilon_n$ та $ \delta _n $ -- члени характеристичних  послідовностей     $ \varepsilon(\psi )$ та $
\delta (\psi ),$ відповідно.
\end{corollary}

Рівності  (\ref{p=2}) в одновимірному випадку ($d=1$) для класів Соболєва $W^r_2$ (при $\psi(k)=k^{-r}$, $r\in {\mathbb N}$)
отримав у 1936 році  А.\,М.~Колмогоров \cite{Kolmogoroff_1936}, який започаткував новий напрям в теорії наближень, пов'язаний з
дослідженням поперечників різних функціональних класів.

Як випливає з (\ref{p=2}),  у просторі $L_2$ поперечники множин $L_2^{\psi }$ реалізують
 суми Фур'є (\ref{w6}).

 Зазначимо, що відомі класи диференційовних функцій Соболєва отримуються з  $L^{\psi } U_{_{\scriptstyle L_p}}$, якщо
 $\psi(\bf k)$ має вигляд   (\ref{w35}) з
  \begin{equation}\label{w45}%$$
 \psi
 _j(k_j)=\left \{ \begin{matrix} 1, \hfill & k_j=0, \hfill \\  (ik_j)^{-r_j}, \hfill & k_j\not =0,
 \hfill\end{matrix}\right.\quad j=\overline {1,d},\
 r_j\in {\mathbb R}. %\eqno (10.12)$$
     \end{equation}

Нехай $d=2$, а послідовності $\psi _1(k_1)$ та $\psi _2(k_2)$ означені рівностями  (\ref{w45}) за умови  $r_1=r_2=r>0.$
У цьому випадку. класи $L^{\psi } U_{_{\scriptstyle L_p}}$  з точки зору знаходження поперечників
вперше розглядалися К.\,І.~Бабенком в  \cite{Babenko_1960_5, Babenko_1960_2}, який цьому випадку фактично отримав співвідношення (\ref{p=2}).

В цій ситуації  характеристична послідовність $\varepsilon(\psi )$ складається
з елементів $\varepsilon_n=n^{-r},$ $n\in {\mathbb N}, $ множини $g_n(\psi)$ -- множини
векторів ${\bf k}=(k_1,k_2)\in {\mathbb Z}^2$, які задовольняють умову
 $$
 k_1'k_2'\le n,
 $$
 де
 $$
 k_j'=\left \{ \begin{matrix} 1,
 \hfill & k_j=0, \hfill \\ |k_j|, \hfill & k_j\not =0 , \ \ j=1,2
 .\hfill \end{matrix}\right.
 $$
 Такі множини  вперше з'явились у згаданих роботах К.\,І.~Бабенка  \cite{Babenko_1960_5, Babenko_1960_2} і зараз їх
 прийнято називати
 гіперболічними хрестами.

\subsection{Класи ${\mathcal F}_{q,r}^{\psi}$ та їх апроксимативні характеристики}\label{F^psi results}

\noindent{\bf \ref{F^psi results}.1. Означення класів ${\mathcal F}_{q,r}^{\psi}$.} Позначимо через $l_p^d$, $0<p\le \infty$, простір всіх послідовностей ${\bf x}=\{x_k\}_{k=1}^d\in {\mathbb R}^d$ зі стандартною $l_p$-нормою (квазі-нормою)
 \[
 |{\bf x}|_p:=\|{\bf x}\|_{_{\scriptstyle l_p^d}}=\left\{\begin{matrix}
 (\sum_{k=1}^d |x_k|^p)^{1/p},\quad 0<p<\infty,\\ \sup_{1\le k\le d}|x_k|,\quad\quad\quad p=\infty.\end{matrix}\right.
 \]
Розглянемо наступні функціональні класи:
 \[
  {\mathcal F}_{q,r}^{\psi}={\mathcal F}_{q,r}^{\psi}({\mathbb T}^d):=\bigg\{f\in L({\mathbb T}^d):\  \|\{\widehat{f} ({\bf k})/\psi(|{\bf k}|_{r})\}_{{\bf k}\in {\mathbb Z}^d}\|_{l_q({\mathbb Z}^d)}\le 1\bigg\},
 \]
де $\psi=\psi(t)$, $t\ge 1$, -- деяка  фіксована додатна спадна
функція, $\psi(0):=\psi(1)$ і  $0<q,r\le \infty$.

Зазначимо, що коли  $\psi(t)=t^{-s}$, $s\in {\mathbb N}$ і $q=1$,  класи \label{F_qr^s}  ${\mathcal F}_{q,\infty}^{\psi}=:{\mathcal F}_{q,\infty}^{s}$ є множинами функцій, у яких частинні похідні порядку $s$ мають абсолютно збіжні ряди Фур'є. Якщо ж $q=2$, то класи  ${\mathcal F}_{2,\infty}^{s}$ збігаються з класами  Соболєва $W^s_2$. Апроксимативні характеристики класів ${\mathcal F}_{q,r}^{\psi}$  для різних $r\in (0,\infty]$ і різноманітних функцій $\psi$ досліджувались в роботах \cite{DeVore_Temlyakov_1995, Temlyakov_Greedy_1998, Li_2010, Shidlich_Zb2011, Shidlich_Zb2013, Shidlich_Zb2014, Serdyuk_Stepanyuk_2015, Shidlich_2016} та інших. Отримані для цих класів результати знаходять своє застосування до дослідження поведінки апроксимативних характеристик  функціональних класів у просторах Лебега $L_p({\mathbb T}^d)$.

 \vskip 3mm\noindent {\bf \ref{F^psi results}.2. Порядкові оцінки найкращих $n$-членних тригонометричних наближень в просторах $S^{p}$.} Класи ${\mathcal F}_{q,r}^{\psi}$ збігаються з множинами $L_q^{\psi^* }$ у випадку, коли система $\psi^*=\{\psi^*({\bf k})\}_{{\bf k}\in {\mathbb Z}^d}$ задовольняє рівності $\psi^*({\bf k})=\psi(|{\bf k}|_{r})$, ${\bf k}\in {\mathbb Z}^d$.
 Тому для них справджуються наведені вище твердження \ref{Applic_q<p}, \ref{Applic_q>p}, \ref{Applic_Psi_q<p}, \ref{Applic_Psi_q>p} та \ref{Applic_Kolmogorov_W}. З огляду на можливі застосування результатів до задач наближення  у просторах Лебега $L_p({\mathbb T}^d)$ також є корисними  твердження цього підрозділу, у яких, зокрема, встановлено точні порядкові оцінки апроксимативних величин класів
 ${\mathcal F}_{q,r}^{\psi}$.  Для їх формулювання та доведення використовувалися згадані вище твердження, а також наведена у підрозділі \ref{Order_Estimates}1  класифікація  Степанця опуклих функцій.

 Нехай $\psi=\psi(t)$, $t\ge 1$, --   фіксована додатна спадна до нуля функція. Тоді  спадну перестановку $ \widetilde{\psi} =\widetilde{\psi}(j)$, $j=1,2,\ldots$,  системи чисел $\psi(|{\bf k}|_r)$ можна визначити рівністю
\begin{equation}\label{d11}
 \widetilde{\psi}(l)=\psi(m),\quad l\in (V_{m-1},V_{m}],\quad m=1,2,\ldots,
 \end{equation}
 де $V_m:=|\widetilde{\Delta}_{m,r}^d|$ --- кількість елементів множини
\begin{equation}\label{f300}% $$
\widetilde{\Delta}_{m,r}^d:=\{{\bf k}\in {\mathbb Z}^d:|{\bf k}|_{r}\le  m,\ \ \ m=0,1,\ldots\}.
 \end{equation}

Далі, при формулюванні результатів важливо, щоб  при всіх достатньо великих $m$ (більших, ніж деяке додатне число $k_0$) виконувалось співвідношення
 \begin{equation}\label{a2.212}%$$
{M_r(m-c_1)^d}<V_m= |\widetilde{\Delta}_{m,r}^d|\le {M_r(m+c_2)^d},
\end{equation}
де $M_r$, $c_1$ та  $c_2$ --- деякі додатні сталі.

Зрозуміло, що у випадку, коли $r=\infty$, співвідношення (\ref{a2.212}) виконується і $M_\infty={\rm vol}\{{\bf k}\in {\mathbb R}^d:|{\bf k}|_{\infty}\le1\}=2^d$, якщо ж $r=1$, то $M_1={\rm vol}\{{\bf k}\in {\mathbb R}^d:|{\bf k}|_{1}\le  1\}=2^d/d!$.

\begin{proposition}[\cite{Shidlich_Zb2011,Shidlich_2016}]\label{sigma_n B} Нехай $d\ge 1$, $0<r\le\infty$, $0<p,\,q<\infty$,  виконуються умова (\ref{a2.212}),  і   функція $\psi^p$ належить  множині $B$, а при $0<p<q$, крім того, при всіх $t$, більших ніж
деяке число $t_0$, є опуклою  та задовольняє умову (\ref{1.2c121}) з $\beta=d(1/p-1/q)$. Тоді
\begin{equation}\label{ss1}%$$
{e}_n({\mathcal F}_{q,r}^{\psi})_{_{\scriptstyle {\mathcal S}^{p}}}\asymp
 {\psi(n^\frac 1d)}{\,n^{\frac 1p\,-\frac 1q}}.
\end{equation}
\end{proposition}

Враховуючи вигляд оцінки (\ref{ss1}) і те, що умова (\ref{a2.212}) виконується, зокрема, при $r=1$ та $r=\infty$, з даного твердження легко отримати такий наслідок.

\begin{corollary}\label{sigma_n B col} Нехай  $d\ge 1$, $0<p,\,q<\infty$  і   функція $\psi$ задовольняє умови твердження \ref{sigma_n B}. Тоді для довільного $r\in [1,\infty]$ має місце оцінка (\ref{ss1}).
\end{corollary}

Дійсно, для довільних чисел $r\in [1,\infty]$, $0<q<\infty$ і будь-якої додатної спадної функції $\psi$ мають місце вкладення
 \begin{equation}\label{er1}%$$%
{\mathcal F}_{q,1}^{\psi}\subset {\mathcal F}_{q,r}^{\psi}\subset {\mathcal F}_{q,\infty}^{\psi}.
\end{equation}
Тому якщо виконуються умови наслідку \ref{sigma_n B col}, то для $r\in [1,\infty]$
 $$
 \frac{\psi(n^\frac 1d)}{\,n^{\frac 1q-\frac 1p}}\ll {e}_n({\mathcal F}_{q,1}^{\psi})_{_{\scriptstyle {\mathcal S}^{p}}}\le {e}_n({\mathcal F}_{q,r}^{\psi})_{_{\scriptstyle {\mathcal S}^{p}}}\le
{e}_n({\mathcal F}_{q,\infty}^{\psi})_{_{\scriptstyle {\mathcal S}^{p}}}\ll  \frac{\psi(n^\frac 1d)}{\,n^{\frac 1q\,-\frac 1p}}.
 $$

%%%%%%%%%%%%%%%%%%%%%%%%%%%%%%%%%%%%%%%%%%%%%%%%%%%%%%%%%%%%%%%%%%%%%%%%%%%%%%%%%%%%%%%%%%%%%%%%%%%%%%

\begin{proposition}[\cite{Shidlich_Zb2014}]\label{sigma_n M'} Нехай $d\ge 1$,  $0<r\le\infty$, $0<p,\,q<\infty$,  виконуються умова (\ref{a2.212}), а функція $\psi^p$ належить множині ${\mathfrak M}_{\infty}'$ або ${\mathfrak M}_{\infty}^c$
Тоді
 \[
 {e}_n({\mathcal F}_{q,r}^{\psi})_{_{\scriptstyle {\mathcal S}^{p}}}\asymp
 {\psi(m_n)}{(n \alpha(\psi^p,
 m_n))^{\frac 1p-\frac 1q}}\asymp
 {\psi(m_n)}{(n \alpha(\psi,
 m_n))^{\frac 1p-\frac 1q}},
 \]
де
\begin{equation}\label{f*1}%$$
m_n=(n/M_r)^\frac 1d.
\end{equation}
\end{proposition}

У випадку, коли  $d=1$, класи ${\mathcal F}_{q,r}^{\psi}=:{\mathcal F}_{q}^{\psi}$  не залежать від $r$, умова (\ref{a2.212}) виконується з константою $M_r=2$. Тому для довільної функції $\psi^p\in {\mathfrak M}_{\infty}'\cup{\mathfrak M}_{\infty}^c$ та будь-яких $0<p,\,q<\infty$
$$
{e}_n({\mathcal F}_{q}^{\psi})_{_{\scriptstyle {\mathcal S}^{p}({\mathbb T}^1)}}\asymp
 {\psi(n/2)}{(n \alpha(\psi,n/2))^{\frac 1p-\frac 1q}}.
 \eqno(\ref{f41}')
$$

\begin{proposition}[\cite{Shidlich_Zb2013}]\label{sigma_n M'' p<q} Нехай $d\ge 1$,  $0<r\le\infty$, $0<p<q<\infty$, $m\in {\mathbb N}$, $\psi^p\in {\mathfrak M}_{\infty}''$, виконується умова (\ref{a2.212})  та умова
 \begin{equation}\label{f2*}%$$
 k^{(d-1)/\alpha}\psi(k+1)/\psi(k)\to 0,\quad k\to\infty, %\eqno (P_\alpha)
 \end{equation}
з $\alpha=p$. Тоді при всіх $n\in [V_{m-1},V_{m})$
\begin{equation}\label{f41}%$$
{e}_n({\mathcal F}_{q,r}^{\psi})_{_{\scriptstyle {\mathcal S}^{p}}}\asymp
 \psi(m) {(V_{m}-n)^{\frac 1p}}\,{n^{\frac {1-d}{dq}}}.
\end{equation}
\end{proposition}

%\vskip

Якщо виконуються умови твердження \ref{sigma_n M'' p<q} і $n\in [V_{m-1},V_{m})$, то
\begin{equation}\label{f4'p}%$$
{\psi(m)}{n^{\frac{1-d}{qd}}}\ll {e}_n({\mathcal F}_{q,r}^{\psi})_{_{\scriptstyle {\mathcal S}^{p}}}\ll {\psi(m)}{n^{\frac{d-1}{d}(\frac 1p-\frac 1q)}}.
\end{equation}

\begin{proposition}[\cite{Shidlich_Zb2013}]\label{sigma_n M'' p>q} Нехай  $d\ge 1$,  $0<r\le\infty$, $0<q\le p<\infty$, $m\in {\mathbb N}$, $\psi^p\in {\mathfrak M}_{\infty}''$, виконуються умови (\ref{a2.212})  та (\ref{f2*}) при $\alpha=q$.
Тоді

1)  при $n=V_{m-1}$ має місце оцінка
\begin{equation}\label{f41'}%$$
{e}_n({\mathcal F}_{q,r}^{\psi})_{_{\scriptstyle {\mathcal S}^{p}}}\asymp
 \psi(m);
\end{equation}

2) для всіх $\ n\in (V_{m-1}, V_m)$  таких, що
\begin{equation}\label{f42}
q(V_{m}-V_{m-1})\ge p\,(V_{m}-n);
\end{equation}
 справджується оцінка $(\ref{f41})$;

3) для всіх $\ n\in\!\! (V_{m-1}, V_m)$, що не задовольняють умову (\ref{f42}),
\begin{equation}\label{f43}%$$
{e}_n({\mathcal F}_{q,r}^{\psi})_{_{\scriptstyle {\mathcal S}^{p}}}\asymp
{\psi(m)}{(n-V_{m-1})^{\frac 1p-\frac 1q}}.
\end{equation}
\end{proposition}

Зазначимо, що у випадку, коли $r\,{=}\,\infty$, для довільного $m\,{\in}\, {\mathbb N}$ маємо
$V_m=|\widetilde{\Delta}_{m,\infty}^d|=(2m+1)^d$. Тому якщо $n\in [V_{m-1}, V_m)$, то
число $m$ визначається рівністю $m=[{(n+1)^{1/d}}/2]$.

У випадку  $d=1$,  якщо $n\in [V_{m-1}, V_m)$, то $n=V_{m-1}{=}V_{m}-1$, а $m=[(n+1)/2]$. Тому для довільної функції $\psi^p\in {\mathfrak M}_{\infty}''$ та будь-яких $0<p,\,q<\infty$
\begin{equation}\label{f2413q}%$$
{e}_n({\mathcal F}_{q}^{\psi})_{_{\scriptstyle {\mathcal S}^{p}({\mathbb T}^1)}}\asymp   \psi([(n+1)/2]).
\end{equation}

Як бачимо, при $d>1$ у випадку $0<q\le p<\infty$ отримані оцінки істотно залежать від розміщення числа $n$ на півсегменті $[V_{m-1}, V_{m})$. Розглядаючи у твердженнях \ref{sigma_n M'' p<q} та \ref{sigma_n M'' p>q} деякі конкретні підпослідовності $n(m)$ отримаємо такий наслідок.

\begin{corollary}\label{sigma_n M'' col}  Нехай  $d\ge 1$, $m \in \, {\mathbb N}$, $n\in [V_{m-1},V_{m})$, $0<r\le\infty$, $0<p,\,q<\infty$, $\psi^p\in {\mathfrak M}_{\infty}''$, виконуються умови (\ref{a2.212}) та  (\ref{f2*}) при $\alpha=\min\{p,q\}$. Тоді

1) якщо $n=n(m)=V_{m}-c_m$, $c_m\in {\mathbb N}$, $c_m\le K$, $m=1,2,\ldots,$ то для довільних $0<p,\,q<\infty$
\begin{equation}\label{f2411}%$$
{e}_n({\mathcal F}_{q,r}^{\psi})_{_{\scriptstyle {\mathcal S}^{p}}}\asymp
{\psi(m)}{n^{\frac {1-d}{dq}}};
\end{equation}

2) якщо $n=n(m)=V_{m-1}+c_m$, $c_m\in {\mathbb N}$, $c_m\le K$, $m=1,2,\ldots,$ то для довільних $0<p<q<\infty$
\begin{equation}\label{f2412}%$$
{e}_n({\mathcal F}_{q,r}^{\psi})_{_{\scriptstyle {\mathcal S}^{p}}}\asymp
{\psi(m)}\,{n^{\frac{d-1}d(\frac 1p-\frac 1q)}};
\end{equation}
а для довільних $0<q< p<\infty$
\begin{equation}\label{f2413}%$$
{e}_n({\mathcal F}_{q,r}^{\psi})_{_{\scriptstyle {\mathcal S}^{p}}}\asymp
 \psi(m);
\end{equation}

3) якщо ж підпослідовність $n=n(m)$ така, що
\begin{equation}\label{f24131}%$$
(V_{m}-V_{m-1})\asymp (V_{m}-n),
\end{equation}
то для довільних $0<p<q<\infty$  справджується оцінка (\ref{f2412}), а при $0<q\le p<\infty$ оцінка (\ref{f2412}) справджується за умови, що $n=n(m)\not=V_{m-1}$, $m=1,2,\ldots$.
\end{corollary}

%%%%%%%%%%%%%%%%%%%%%%%%%%%%%%%%%%%%%%%%%%%%%%%%%%%%%%%%%%%%%%%%%%%%%%%%%%%%%%%%%%%%%%%%%%%%%%%%%%%%%%%%%%%%%%%%%%%%%%%%%%%%%%%%%%%%%%%%%%%%%%%%%%%%%%%%%

%%%%%%%%%%%%%%%%%%%%%%%%%%%%%%%%%%%%%%%%%%%%%%%%%%%%%%%%%%%%%%%%%%%%%%%%%%%%%%%%%%%%%%%%%%%%%%%%%%%%%%
%%%%%%%%%%%%%%%%%%%%%%%%%%%%%%%%%%%%%%%%%%%%%%%%%%%%%%%%%%%%%%%%%%%%%%%%%%%%%%%%%%%%%%%%%%%%%%%%%%%%%%%%%%%%%%%%%%%%%%%%%%%%%%%%%%%%%%%%%%%%%%%%%%%%%%%%%

\vskip 3mm\noindent {\bf \ref{F^psi results}.3. Порядкові оцінки величин ${\mathscr D}_n^\perp({\mathcal F}_{q,r}^{\psi})_{_{\scriptstyle {\mathcal S}^{p}}}$ та ${\mathscr D}_n({\mathcal F}_{q,r}^{\psi})_{_{\scriptstyle {\mathcal S}^{p}}}$} містяться в наступних  твердженнях.

\begin{proposition}[\cite{Shidlich_Zb2013}]\label{D_n M''} Нехай   $d\ge 1$,  $0<r\le\infty$, $0<p, \,q<\infty$ і $m\in {\mathbb N}$. Тоді

1) якщо $0<q\le p<\infty$, то для довільної додатної спадної до нуля функції $\psi(t)$, $t\ge 1$, при кожному $ n\in [V_{m-1},V_{m})$ справджується рівність
$$
 {\mathscr D}_n({\mathcal F}_{q,r}^{\psi})_{_{\scriptstyle {\mathcal S}^{p}}}={\mathscr D}_n^\perp({\mathcal F}_{q,r}^{\psi})_{_{\scriptstyle {\mathcal S}^{p}}}=\psi(m);
$$

2) якщо ж $0<p<q<\infty$,  виконується умова (\ref{a2.212}), $\psi^p\in {\mathfrak M}_{\infty}''$ і при $\alpha=\frac{pq}{q-p}$ має місце співвідношення (\ref{f2*}), то при кожному $ n\in [V_{m-1},V_{m})$ справджується оцінка
\begin{equation}\label{d18}
{\mathscr D}_n({\mathcal F}_{q,r}^{\psi})_{_{\scriptstyle {\mathcal S}^{p}}}={\mathscr D}_n^\perp({\mathcal F}_{q,r}^{\psi})_{_{\scriptstyle {\mathcal S}^{p}}}\asymp \psi(m)(V_{m}-n)^{\frac 1p-\frac 1q}.
 \end{equation}
\end{proposition}

У випадку $d=1$  при  $0<q\le p<\infty$ з твердження \ref{D_n M''} випливає, що для довільної додатної спадної до нуля функції $\psi(t)$, $t\ge 1$,
$$
 {\mathscr D}_n({\mathcal F}_{q}^{\psi})_{_{\scriptstyle {\mathcal S}^{p}({\mathbb T}^1)}}={\mathscr D}_n^\perp({\mathcal F}_{q}^{\psi})_{_{\scriptstyle {\mathcal S}^{p}({\mathbb T}^1)}}=\psi([(n+1)/2]),
$$
а при   $0<p<q<\infty$ для  довільної функції $\psi^p\in {\mathfrak M}_{\infty}''$
$$
{\mathscr D}_n({\mathcal F}_{q}^{\psi})_{_{\scriptstyle {\mathcal S}^{p}({\mathbb T}^1)}}={\mathscr D}_n^\perp({\mathcal F}_{q}^{\psi})_{_{\scriptstyle {\mathcal S}^{p}({\mathbb T}^1)}}\asymp  \psi([(n+1)/2]).
$$

\begin{proposition}[\cite{Shidlich_Zb2011,Shidlich_2016}]\label{D_n B} Нехай    $d\ge 1$, $0<r\le\infty$, $0<p, \,q<\infty$ і виконується умова (\ref{a2.212}). Тоді

1) якщо $0<q\le p<\infty$, то для довільної функції $\psi\in B$ має місце оцінка
    \begin{equation}\label{f311}
{\mathscr D}_n({\mathcal F}_{q,r}^{\psi})_{_{\scriptstyle {\mathcal S}^{p}}}={\mathscr D}_n^\perp({\mathcal F}_{q,r}^{\psi})_{_{\scriptstyle {\mathcal S}^{p}}}\asymp \psi(m_n)\asymp \psi(n^{1/d}),
\end{equation}

2) якщо ж $0<p<q<\infty$, а функція $\psi^p$ належить  $B$ і при всіх $t$, більших ніж
деяке число $t_0$, є опуклою  та задовольняє умову (\ref{1.2c121}) при $\beta=d(1/p-1/q)$, то
    \begin{equation}\label{f31122}
{\mathscr D}_n({\mathcal F}_{q,r}^{\psi})_{_{\scriptstyle {\mathcal S}^{p}}}={\mathscr D}_n^\perp({\mathcal F}_{q,r}^{\psi})_{_{\scriptstyle {\mathcal S}^{p}}}\asymp
\psi(n^\frac 1d)n^{\frac 1p-\frac 1q}.
\end{equation}
\end{proposition}

 \begin{proposition}[\cite{Shidlich_Zb2014}]\label{D_n M'} Нехай   $d\ge 1$, $0<r\le\infty$, $0<p,\,q<\infty$,  виконуються умова (\ref{a2.212}), а функція $\psi^p$ належить ${\mathfrak M}_{\infty}'$ або ${\mathfrak M}_{\infty}^c$. Тоді для довільних $0<q\le p<\infty$ має місце співвідношення
 \begin{equation}\label{f312}
{\mathscr D}_n({\mathcal F}_{q,r}^{\psi})_{_{\scriptstyle {\mathcal S}^{p}}}={\mathscr D}_n^\perp({\mathcal F}_{q,r}^{\psi})_{_{\scriptstyle {\mathcal S}^{p}}}\asymp \psi(m_n),
\end{equation}
 а для довільних $0<p<q<\infty$ --- співвідношення
    \begin{equation}\label{f311223}
{\mathscr D}_n({\mathcal F}_{q,r}^{\psi})_{_{\scriptstyle {\mathcal S}^{p}}}={\mathscr D}_n^\perp({\mathcal F}_{q,r}^{\psi})_{_{\scriptstyle {\mathcal S}^{p}}}\asymp
\psi(m_n)(n \alpha(\psi,  m_n))^{\frac 1p-\frac 1q}.
\end{equation} \end{proposition}

На підставі твердження \ref{D_n M'} та означення множини   ${\mathfrak M}_\infty^c$ бачимо, що   для довільної функції $\psi^p\in {\mathfrak M}_{\infty}^c$ та будь-яких $0<p,\,q<\infty$ справджується оцінка
\begin{equation}\label{fa1}
{\mathscr D}_n({\mathcal F}_{q}^{\psi})_{_{\scriptstyle {\mathcal S}^{p}({\mathbb T}^1)}} ={\mathscr D}_n^\perp({\mathcal F}_{q}^{\psi})_{_{\scriptstyle {\mathcal S}^{p}({\mathbb T}^1)}} \asymp {e}_n({\mathcal F}_{q}^{\psi})_{_{\scriptstyle {\mathcal S}^{p}({\mathbb T}^1)}}  \asymp   \psi(n/2).
\end{equation}

Співставляючи оцінки  величин ${e}_n({\mathcal F}_{q,r}^{\psi})_{_{\scriptstyle {\mathcal S}^{p}}}$,  ${\mathscr D}_n({\mathcal F}_{q,r}^{\psi})_{_{\scriptstyle {\mathcal S}^{p}}}$ та ${\mathscr D}_n^\perp({\mathcal F}_{q,r}^{\psi})_{_{\scriptstyle {\mathcal S}^{p}}}$, робимо висновок, що у випадку, коли $d=1$ і $\psi^p\in {\mathfrak M}''_{\infty}\cup {\mathfrak M}^c_{\infty}$, для будь-яких чисел $0<p,\,q<\infty$
$$%\begin{equation}\label{d61}%$$
{e}_n({\mathcal F}_{q}^{\psi})_{_{\scriptstyle {\mathcal S}^{p}({\mathbb T}^1)}}\asymp {\mathscr D}_n({\mathcal F}_{q}^{\psi})_{_{\scriptstyle {\mathcal S}^{p}({\mathbb T}^1)}}={\mathscr D}_n^\perp({\mathcal F}_{q}^{\psi})_{_{\scriptstyle {\mathcal S}^{p}({\mathbb T}^1)}}.
$$%\end{equation}
Якщо $d>1$, то аналогічне співвідношення
 \begin{equation}\label{d6}%$$
{e}_n({\mathcal F}_{q,r}^{\psi})_{_{\scriptstyle {\mathcal S}^{p}}}\asymp {\mathscr D}_n({\mathcal F}_{q,r}^{\psi})_{_{\scriptstyle {\mathcal S}^{p}}}={\mathscr D}_n^\perp({\mathcal F}_{q,r}^{\psi})_{_{\scriptstyle {\mathcal S}^{p}}}
\end{equation}
справджується коли  $0<p<q$ і функція $\psi$ задовольняє умови тверджень \ref{sigma_n B} та \ref{sigma_n M'}. Якщо ж $d>1$, а $0<q\le p$, то для довільної функції, яка задовольняє умови тверджень \ref{sigma_n B} та \ref{sigma_n M'}
 \begin{equation}\label{d318}
{e}_n({\mathcal F}_{q,r}^{\psi})_{_{\scriptstyle {\mathcal S}^{p}}}=o\Big({\mathscr D}_n({\mathcal F}_{q,r}^{\psi})_{_{\scriptstyle {\mathcal S}^{p}}}\Big),\ n\to\infty.
 \end{equation}

 Коли $\psi^p$ належить множині ${\mathfrak M}_{\infty}''$ і  $d>1$, то при   $0<q<p$ співвідношення (\ref{d6}) виконується (за відповідних додаткових умов тверджень \ref{sigma_n M'' p<q} та \ref{sigma_n M'' p>q}) для підпослідовності вигляду $n=n(m)=V_{m-1}+c_m$, $c_m\in {\mathbb N}$, $c_m\le K$, $m=1,2,\ldots$, а при   $0<p=q$ --- для підпослідовності $n=n(m)=V_{m-1}$, $m=1,2,\ldots$ У випадку, коли  $0<p<q$   співвідношення (\ref{d6}) виконується для підпослідовностей $n=n(m)$, що задовольняють умову (\ref{f24131}); якщо ж дана підпослідовність \mbox{$n=n(m)$} така, що \mbox{$(V_{m}-n)=o(V_{m}-V_{m-1}),$} $m\to \infty$, то справджується співвідношення (\ref{d318}).

%%%%%%%%%%%%%%%%%%%%%%%%%%%%%%%%%%%%%%%%%%%%%%%%%%%%%%%%%%%%%%%%%%%%%%%%%%%%%%%%%%%%%%%%%%%%%%%%%%%%%%%%%%%%%%%%%%%%%%%%%%%%%%%%%%%%%%%%%%%%%%%%%%%%%%%%%

\subsection{Прямі та обернені теореми наближення в просторах ${\mathcal S}^p$.}\label{DI_TH_Sp}

 \noindent {\bf \ref{DI_TH_Sp}.1. Попередні позначення.} Нехай  $f$ -- довільна функція з простору $L=L({\mathbb T}^d)$  з рядом Фур'є  вигляду (\ref{b8}). Для будь-якого  $\nu\in {\mathbb N}_0:=\{0,1,2,\ldots\}$ покладемо
 \begin{equation}\label{H_nu}
H_{\nu}(f)({\bf x}):=\sum_{|{\bf k}|_1=\nu}\widehat f({\bf
k}){\rm e}^{{\mathrm i}({\bf k,x})},~|{\bf
k}|_1:=\sum_{j=1}^d|k_j|.
 \end{equation}
Тоді ряд Фур'є функції  $f$ можна записати у вигляді
 \begin{equation}\label{Fourier series}
    S[f]({\bf x}):=\sum_{{\bf k}\in {\mathbb Z}^d} \widehat f({\bf k}){\rm e}^{{\mathrm i}({\bf x,k})}=\sum_{\nu=0}^\infty H_{\nu}(f)({\bf x}).
 \end{equation}

 Розглянемо множину  ${\mathscr T}_n^\vartriangle$, $n\in {\mathbb N}_0$, всіх поліномів вигляду
  $
 \tau_n({\bf x}):=%\sum_{\nu=0}^n
 \sum_{|{\bf k}|_1\le n} a_{\bf k}{\rm e}^{{\mathrm i}({\bf x,k})},
 $
де $a_{\bf k}$ -- довільні комплексні числа.

 Величину
\begin{equation}\label{E_n_def}
E_n^\vartriangle(f)_{_{\scriptstyle {\mathcal S}^p}}=\inf\limits_{
\tau_{n-1}\in {\mathscr T}_{n-1}^\vartriangle}\|f-\tau_{n-1}\|_{_{\scriptstyle {\mathcal S}^p}}
\end{equation}
називають найкращим наближенням функції  $f \in   {\mathcal S}^p$ порядку ${n-1}$ поліномами, побудованим за трикутними областями.

Модулем гладкості функції $f\in {\mathcal S}^p$ порядку  $\alpha>0$ називають величину
\[
\omega_{\alpha}^{\vartriangle} (f,t)_{_{\scriptstyle {\mathcal S}^p}}:=\sup\limits_{|h|\le t}\|\Delta_h^\alpha f\|_{_{\scriptstyle {\mathcal S}^p}}=\sup\limits_{|h|\le t}\Big\|\sum\limits_{j=0}^\infty (-1)^j {\alpha \choose j} f(\cdot-jh)\Big\|_{_{\scriptstyle {\mathcal S}^p}}.
\]
Функції  $\omega_{\alpha}^{\vartriangle} (f,t)_{_{\scriptstyle
{\mathcal S}^p}}$ володіють усіма звичайними властивостями класичних модулів гладкості.  Зокрема, має місце таке твердження.

\begin{lemma}\label{Lemma_2} Нехай   $f, g\in {\mathcal S}^p$, $\alpha\ge \beta>0$, $t, t_1, t_2\ge 0$. Тоді

{\rm (i)} $\omega_{\alpha}^{\vartriangle}(f,t)_{_{\scriptstyle
{\mathcal S}^p}}$, $t\in (0,\infty)$, є невід'ємною неперервною зростаючою функцією   і $\lim\limits_{t\to 0+}
\omega_{\alpha}^{\vartriangle}(f,t)_{_{\scriptstyle {\mathcal
S}^p}}=0$;

{\rm (ii)} $\omega_{\alpha}^{\vartriangle}(f,t)_{_{\scriptstyle
{\mathcal S}^p}}\le 2^{\{\alpha-\beta\}}
\omega_{\beta}^{\vartriangle}(f,t)_{_{\scriptstyle {\mathcal
S}^p}}$, де $\{\alpha\}=\inf\{k\in {\mathbb N_0}:k\ge \alpha\}$;

{\rm (iii)} $\omega_{\alpha}^{\vartriangle}(f+g,t)_{_{\scriptstyle
{\mathcal S}^p}}\le
\omega_{\alpha}^{\vartriangle}(f,t)_{_{\scriptstyle {\mathcal
S}^p}}+\omega_{\alpha}^{\vartriangle}(g,t)_{_{\scriptstyle {\mathcal
S}^p}}$;

{\rm (iv)}  $\omega _1^{\vartriangle}(f,t_1+t_2)_{_{\scriptstyle
{\mathcal S}^p}}\le  \omega _1^{\vartriangle}(f,t_1)_{_{\scriptstyle
{\mathcal S}^p}}+ \omega _1^{\vartriangle}(f,t_2)_{_{\scriptstyle
{\mathcal S}^p}}$;

{\rm (v)}   $\omega _\alpha(f,t)_{_{\scriptstyle {\mathcal S}^p}}\le
2^{\{\alpha\}}\|f\|_{_{\scriptstyle {\mathcal S}^p}}$.
\end{lemma}

 \vskip 3mm\noindent {\bf \ref{DI_TH_Sp}.2.  Прямі теореми наближення.}

\begin{proposition}[\cite{Abdullayev_Ozkartepe_Savchuk_Shidlich_2019}]\label{Proposition 1}  Нехай $\psi=\{\psi({\bf k})\}_{{\bf k}\in\mathbb Z^d}$
-- система чисел, підпорядкована умовам  (\ref{w5})  та (\ref{w17}). Якщо  для функції
 $f\in {\mathcal S}^{p}$, $1\le p<\infty$, існує похідна $f^{\psi}$ з простору ${\mathcal S}^{p}$,
 тоді при   $n\in {\mathbb N}$,
\[
{E_{n}^\vartriangle} (f)_{_{\scriptstyle {\mathcal S}^p}} \le
\varepsilon_n {E_{n}^\vartriangle} (f^{\psi})_{_{\scriptstyle
{\mathcal S}^p}},\quad \mbox{де}\quad
    \varepsilon_n=\max\limits_{|{\bf k}|_1 \ge n} |\psi({\bf k})|.
\]
\end{proposition}

 Сформулюємо тепер прямі теореми наближення в просторах  ${\mathcal S}^p$ в термінах найкращих наближень та модулів гладкості функцій. Зокрема, наведемо  нерівності типу Джексона  вигляду
 \[
 E_n^\vartriangle(f)_{_{\scriptstyle {\mathcal S}^p}}\le
 K(\tau )\omega_{\alpha}^{\vartriangle}\Big(f, \frac {\tau }n\Big)_{_{\scriptstyle {\mathcal S}^p}}, \ \ \tau >0,
\]
і розглянемо питання про точні константи в цих нерівностях при фіксованих   $n$, $\alpha$, $\tau$ та $p$. Для цього розглянемо величину
\begin{equation}\label{(6.6)}
 \!K_{n,\alpha,p}(\tau )=
  \sup \left \{\frac {E_n^\vartriangle(f)_{_{\scriptstyle {\mathcal S}^p}}}
  {\omega_{\alpha}^{\vartriangle}(f, \frac {\tau}n)_{_{\scriptstyle {\mathcal S}^p}}}: f\in {\mathcal S}^{p} \cap  L_{1,Y},
\ f\not \equiv {\rm const }\right \},
 \end{equation}
де   $Y:=\mathbb Z^d_+\cup\mathbb Z^d_-$,
$\mathbb Z^d_+$ та $\mathbb Z^d_-$ -- підмножини векторів ${\bf z}\in \mathbb Z^d$, всі координати яких відповідно невід'мні або від'ємні,
 \begin{equation}\label{L_1,Y}
             L_{1,Y}:=L_{1,Y}({\mathbb T}^d)=\{f\in  L({\mathbb T}^d):\, \widehat f({\bf k})=0\  \forall{\bf k}\in\mathbb Z^d\setminus Y \}.
  \end{equation}

Через $M(\tau)$, $\tau>0$, позначимо множину обмежених неспадних функцій  $\mu$, відмінних від константи на  $[0, \tau].$

\begin{theorem}[\cite{Stepanets_Serdyuk_UMZh2002, Abdullayev_Ozkartepe_Savchuk_Shidlich_2019}]\label{Theorem_2.1} Нехай $f\in {\mathcal S}^{p} \cap  L_{1,Y} $, $1\le p<\infty$. Тоді для довільних $\tau >0$, $n\in {\mathbb N}$ та $\alpha>0$ справджується нерівність
\begin{equation}\label{(6.7)}
E_n^\vartriangle(f)_{_{\scriptstyle {\mathcal S}^p}}\le
C_{n,\alpha,p}(\tau ) \omega_{\alpha}^{\vartriangle}\Big(f, \frac
{\tau }n\Big)_{_{\scriptstyle {\mathcal S}^p}},
  \end{equation}
де
 \begin{equation}\label{(6.8)}
 C_{n,\alpha,p}(\tau ):=\left(\inf\limits _{\mu \in  M(\tau )} \frac {\mu  (\tau ) - \mu  (0)}{2^{\frac {\alpha p}2}I_n(\tau ,\mu
 )}\right)^{1/p},
  \end{equation}
і
 \begin{equation}\label{(6.9)}
 I_n(\tau ,\mu  )=I_{n,\alpha,p}(\tau ,\mu  )=
 \inf\limits _{\nu \in {\mathbb N}:\nu \ge n} \int\limits _0^{\tau }\Big(1-\cos \frac {\nu }nt\Big)^{\frac {\alpha p}2}d\mu  (t).
 \end{equation}
Крім цього, існує функція $\mu _*\in M(\tau )$, яка реалізує точну нижню межу в  (\ref{(6.8)}). Нерівність
(\ref{(6.7)}) є непокращуваною на множині всіх функцій  $f\in {\mathcal S}^{p} \cap  L_{1,Y}$, $f\not \equiv {\rm const}$, в тому
сенсі, що для довільних   $\alpha>0$ та $n\in {\mathbb N}$  маємо
  \begin{equation}\label{(6.10)}
  C_{n,\alpha,p}(\tau )= K_{n,\alpha,p}(\tau ).
  \end{equation}
\end{theorem}

Зазначимо, що у просторах  $L_2({\mathbb T}^1)$ при $\alpha = 1$ дане твердження доведено О.\,Г.~Бабенком \cite{Babenko_1986}.
У просторах  ${\mathcal S}^p$   цей та інші результати цього підрозділу отримано для функцій однієї та багатьох змінної відповідно в роботах   \cite{Stepanets_Serdyuk_UMZh2002} та \cite{Abdullayev_Ozkartepe_Savchuk_Shidlich_2019}.

Зазначимо також, що для  $f\in {\mathcal S}^p$ умова $\widehat f({\bf k})=0$, ${\bf k}\in\mathbb Z^d\setminus \mathbb
Z^d_\pm$ в теоремі \ref{Theorem_2.1} взагалі кажучи є необхідною. Наприклад, розглянемо функцію $f({\bf x})={\rm
e}^{{\mathrm i}({\bf k^*}, {\bf x})}$, де $k^*=(l,-l,0,\ldots)$, $l\in
{\mathbb N}$. Тоді при всіх  $n<2l$ маємо $E_n(f^*)_{_{\scriptstyle {\mathcal S}^p}}=1$, однак $\omega(f^*,t)\equiv 0$.

\begin{corollary}[\cite{Stepanets_Serdyuk_UMZh2002, Abdullayev_Ozkartepe_Savchuk_Shidlich_2019}]\label{Theorem 2.2.} Нехай  $f\in {\mathcal S}^{p} \cap  L_{1,Y}$, $1\le p<\infty$. Тоді для довільних
 $n\in {\mathbb N}$ та $\alpha>0$ справджується нерівність
 \begin{equation}\label{(6.11)}
E_n^\vartriangle(f)^p_{_{\scriptstyle {\mathcal S}^p}}\le \frac
1{2^{\frac {\alpha p }2}I_n(\frac {\alpha p }2)}\int\limits_0^{\pi
}\omega_{\alpha}^{\vartriangle}\Big(f, \frac
t{n}\Big)_{_{\scriptstyle {\mathcal S}^p}}^p \sin t{\mathrm d}t,
  \end{equation}
де
 \begin{equation}\label{(6.12)}
I_n(\lambda ):=\inf\limits _{\nu \in {\mathbb N}:\nu \ge
n}\int\limits _0^{\pi } \Big(1-\cos \frac {\nu }nt\Big)^{\lambda
}\sin t{\mathrm d}t, \ \ \lambda >0, \ \ n\in {\mathbb N}.
   \end{equation}
Якщо при цьому  $\frac {\alpha p }2\in {\mathbb N},$ то
 \begin{equation}\label{(6.13)}
 I_n\Big(\frac {\alpha p }2\Big)=\frac {2^{\frac {\alpha p }2+1}}{\frac {\alpha p }2+1},
   \end{equation}
і нерівність  (\ref{(6.11)}) не може бути покращена  для  $n\in {\mathbb N}$ в тому сенсі, що для кожного
 $n\in {\mathbb N}$ існує функція  $f^*\in {\mathcal S}^{p} \cap  L_{1,Y} $ така, що
  \begin{equation}\label{(6.53)}
E_n^\vartriangle(f^*)^p_{_{\scriptstyle {\mathcal S}^p}} =\frac
{\frac {\alpha p }2+1}{2^{\alpha p+1}}\int\limits _0^{\pi }
\omega_{\alpha}^{\vartriangle}\Big(f^*, \frac
t{n}\Big)_{_{\scriptstyle {\mathcal S}^p}}^p  \sin t{\mathrm d}t.
  \end{equation}
\end{corollary}

В наступному твердженні містяться оцінки зверху для констант  $K_{n,\alpha,p}(\tau )$ при $\tau=\pi$, які не залежать від  $n$
і є непокращуваними у низці важливих випадків.

\begin{corollary}[\cite{Stepanets_Serdyuk_UMZh2002, Abdullayev_Ozkartepe_Savchuk_Shidlich_2019}]\label{Theorem 2.3.}
   Для довільних  $1\le p<\infty$, $\alpha>0$ та  $n\in {\mathbb N}$ справджуються нерівності
  \begin{equation}\label{(6.14)}
     K_{n,\alpha,p}^p(\pi )\le \frac 1{2^{\frac {\alpha p }2-1}I_n(\frac
     {\alpha p }2)}\le \frac {\frac {\alpha p }2+1}{2^{\alpha p}+2^{\frac
     {\alpha p }2-1}(\frac {\alpha p}2+1)\sigma (\frac {\alpha p }2)},
 \end{equation}
 де величина $I_n(\lambda )$ означена рівністю $(\ref{(6.12)})$ і
 \[
    \sigma (\lambda ):=-\sum _{m=[\frac {\lambda }2]+1}^{\infty } {\lambda \choose 2m}
    \frac 1{2^{2m-1}}\bigg(\frac {1-(-1)^{[\lambda ]}}2 {2m \choose m}
 \]
 \[
    - \sum _{j=0}^{m-1} {2m \choose j}\frac 2{2(m-j)^2-1}\bigg),\ \ \ \lambda >0.
  \]
Якщо при цьому  $\frac {\alpha p }2\in {\mathbb N},$ то величина  $\sigma (\frac {\alpha p }2)=0$ і
\begin{equation}\label{6.14)'}
  K_{n,\alpha,p}^p(\pi )\le \frac {\frac {\alpha p }2+1}{2^{\alpha p}}.
 \end{equation}
\end{corollary}

Наступне твердження встановлює оцінку величин  $K_{n,\alpha,p}(\pi )$  рівномірно обмежену відносно усіх параметрів  $\alpha>0$, $n\in {\mathbb N}$ та $1\le p<\infty$.

\begin{corollary}[\cite{Stepanets_Serdyuk_UMZh2002, Abdullayev_Ozkartepe_Savchuk_Shidlich_2019}] \label{Theorem 2.4.}  Нехай
$f\in {\mathcal S}^{p} \cap  L_{1,Y} $, $1\le p<\infty$,  ${f\not \equiv {\rm const}}.$ Тоді для довільних $\alpha>0$ та $n\in {\mathbb N}$
  \begin{equation}\label{(6.16)}
  E_n^\vartriangle(f)_{_{\scriptstyle {\mathcal S}^p}} < \frac {4}{3\cdot 2^{\alpha/2}}\omega_\alpha
  \Big(f, \frac {\pi }n\Big)_{_{\scriptstyle {\mathcal S}^p}}.
   \end{equation}
\end{corollary}

У просторах  $L_2({\mathbb T}^1)$ при  $\alpha =1$ нерівність  (\ref{(6.11)})
доведено М.\,І.~Чернихом   \cite{Chernykh_1967, Chernykh_1967_MZ}. Нерівності такого типу, а також суміжні питання, пов'язані із обчисленнням значень поперечників класів функцій, що задаються мажорантами їх модулів неперервності,  досліджувалися в роботах
\cite{Taikov_1976, Taikov_1979, Voicexivskij_UMZh2003, Serdyuk_2003, Vakarchuk_2004, Gorbachuk_Grushka_Torba_2005, Vakarchuk_Shchitov_2006, Timan_M2009,
Vakarchuk_2016, Babenko_Konareva_2019, Abdullayev_Chaichenko_Shidlich_2020, Abdullayev_Serdyuk_Shidlich_2020} та ін.

 \vskip 3mm\noindent {\bf \ref{DI_TH_Sp}.3. Обернені теореми наближення.}

Перед формулюванням оберненої теореми наближення наведемо також нерівність Бернштейна, у якій норма узагальненої похідної тригонометричного полінома оцінюється через норму самого полінома (див., наприклад,  \cite[Гл.~4]{Timan_M2009}).

\begin{proposition}[\cite{Stepanets_Serdyuk_UMZh2002, Abdullayev_Ozkartepe_Savchuk_Shidlich_2019}]\label{Proposition 2}
Нехай  $\psi =\{\psi ({\bf k})\}_{{\bf k}\in {\mathbb Z}^d} $ -- система чисел, які задовольняють умову  (\ref{w17}). Тоді для довільного   $\tau_n\in {\mathscr T}_{n}$, $n\in \mathbb{N}$, справджується нерівність
\[
    \|\tau^\psi_n \|_{_{\scriptstyle  {\mathcal S}^p}}\le \frac 1{\epsilon_n}\|\tau_n\|_{_{\scriptstyle  {\mathcal S}^p}}, \quad \mbox{де}\ \
    \epsilon_n:=\min_{0<| {\bf k}|_1 \le n}|\psi({\bf k})|.
\]
\end{proposition}

\begin{corollary}\label{Corollary 1.1}  Нехай $\psi({\bf k})=\nu^{-r}$,
$|{\bf k}|_1=\nu$, $\nu=0,1,\ldots$, $r\ge 0$. Тоді для довільного полінома
$\tau_n\in {\mathscr T}_{n}$, $n\in \mathbb{N}$
$$
    \|\tau^\psi_n \|_{_{\scriptstyle  {\mathcal S}^p}} = \|\tau^{(r)}_n \|_{_{\scriptstyle  {\mathcal S}^p}} \le
    n^r \|\tau_n\|_{_{\scriptstyle  {\mathcal S}^p}}.
$$
\end{corollary}

Обернена апроксимаційна теорема в просторі ${\mathcal S}^p$ має такий вигляд.

\begin{theorem}[\cite{Stepanets_Serdyuk_UMZh2002, Abdullayev_Ozkartepe_Savchuk_Shidlich_2019}]\label{Theorem_2}
Нехай  $f\in {\mathcal S}^p$, $1\le p<\infty$. Тоді для довільних  $n\in {\mathbb N}$ та $\alpha>0$
%справджується нерівність
\begin{equation}\label{(6.71)}
 \omega_{\alpha}^{\vartriangle}\Big(f, \frac {\pi }n\Big)_{_{\scriptstyle {\mathcal S}^p}}\le
\frac {\pi ^\alpha }{n^\alpha }\bigg(\sum _{\nu =1}^n(\nu ^{\alpha p
}-(\nu -1)^{\alpha p })E_{\nu}^\vartriangle(f)_{_{\scriptstyle
{\mathcal S}^p}}^p\bigg)^{1/p}.
    \end{equation}
\end{theorem}

Зазначимо, що в   (\ref{(6.71)}) сталу   $\pi^\alpha$, взагалі кажучи, зменшити не можна, оскільки  для
довільного  числа  $\varepsilon>0$ знайдеться функція  $f^*\in {\mathcal S}^p$ така, що при всіх
$n$, більших деякого номера  $n_0$ виконується протилежна нерівність
 \[
\omega_{\alpha}^{\vartriangle}\Big(f^*, \frac {\pi
}n\Big)_{_{\scriptstyle {\mathcal S}^p}}> \frac {\pi
^\alpha-\varepsilon }{n^\alpha }\bigg(\sum _{\nu =1}^n(\nu ^{\alpha
p }-(\nu -1)^{\alpha p })E_{\nu}^\vartriangle(f^*)_{_{\scriptstyle
{\mathcal S}^p}}^p\bigg)^{1/p}.
 \]

Оскільки $\nu ^{\alpha p}-(\nu -1)^{\alpha p}\le \alpha p \nu ^{\alpha p-1},$ то з нерівності
(\ref{(6.71)}) випливає, що
 \begin{equation}\label{(6.71')}
 \omega_{\alpha}^{\vartriangle}\Big(f, \frac {\pi }n\Big)_{_{\scriptstyle {\mathcal S}^p}}\le
 \frac {\pi ^\alpha (\alpha p
 )^{1/p}}{n^\alpha }\bigg(\sum _{\nu =1}^n\nu ^{\alpha p -1}E_{\nu}^\vartriangle(f)^p
 _{_{\scriptstyle {\mathcal S}^p}}\bigg)^{1/p}.
 \end{equation}
Звідси, зокрема, отримуємо такий наслідок.

\begin{corollary}[\cite{Stepanets_Serdyuk_UMZh2002, Abdullayev_Ozkartepe_Savchuk_Shidlich_2019}]\label{Corollary 1}   Нехай  $f\in {\mathcal
S}^p$, $1\le p<\infty$,  послідовність найкращих наближень ${E_{n}^\vartriangle} (f)_{_{\scriptstyle {\mathcal S}^p}}$
функції  $f$ задовольняє співвідношення
 $
 {E_{n}^\vartriangle} (f)_{_{\scriptstyle {\mathcal S}^p}}= O(n^{-\beta })
 $
при деякому  $\beta >0$. Тоді для всіх  $\alpha>0$,
 \[
 \omega_{\alpha}^{\vartriangle}(f, t)_{_{\scriptstyle {\mathcal S}^p}}=\left \{ \begin{matrix}  O(t^{\beta }) & \hfill \mbox{при}
\ \ \beta <\alpha, \hfill \cr O(t^\alpha|\ln t|^{1/p}) & \hfill \mbox{при}\ \ \beta =\alpha, \hfill \cr O(t^\alpha) & \hfill
\mbox{при}
\ \ \beta >\alpha.\hfill \end{matrix} \right.
\]
\end{corollary}

У просторах  ${\mathcal S}^p({\mathbb T}^1)$  $2\pi$-періодичних функцій однієї змінної нерівності
$(\ref{(6.71')})$ були отримані в   \cite{Sterlin_1972} та \cite{Stepanets_Serdyuk_UMZh2002}. У просторах ${\mathcal S}^p({\mathbb T}^d)$ функцій багатьох змінних ці нерівності отримано в  \cite{Abdullayev_Ozkartepe_Savchuk_Shidlich_2019}. У просторах
 $L_p({\mathbb T}^d)$ нерівності типу $(\ref{(6.71')})$ доведено М.\,П.~Тіманом (див.
 \cite{Timan_M_MS1958, Timan_M_DAN1958} та \cite[Гл. 2]{Timan_M2009}).

 Прямі та обернені теореми наближення
функцій, заданих на сфері, у просторах   $S^{p,q}(\sigma^d)$, $d\ge 3$, отримано в роботах \cite{Lasuriya_2007, Lasuriya_2015}.

 \vskip 3mm\noindent {\bf \ref{DI_TH_Sp}.4. Конструктивні характеристики класів функцій, визначених
  їх модулями гладкості}

Нехай  $\omega$ -- довільна мажоранта, визначена на відрізку  $[0,1]$. Для фіксованого
$\alpha>0$ покладемо
\begin{equation} \label{omega-class}
    {\mathcal S}^{p}H^{\omega}_{\alpha} =
    \Big\{f\in {\mathcal S}^{p}:  \quad \omega_\alpha^{\vartriangle}(f; \delta)_{_{\scriptstyle {\mathcal S}^p}}=
    {\mathcal O}  (\omega(\delta)),\quad  \delta\to 0+\Big\}.
\end{equation}
Далі, розглядаємо мажоранти  $\omega(t)$, $t\in [0,1]$, які задовольняють наступні умови: \textbf{1)} $\omega(\delta)$ неперервна на $[0,1]$;\
\textbf{  2)} $\omega(\delta)\uparrow$;\  \textbf{  3)}
$\omega(\delta)\not=0$ для $\delta\in (0,1]$;\  \textbf{  4)}
$\omega(\delta)\to 0$ при $\delta\to 0$, а також відомі умови Барі $({\mathscr B}_\alpha)$ та $({\mathscr B})$  (див., наприклад, \cite{Bari_Stechkin_1956}):
\[
({\mathscr B}_\alpha),\ \alpha>0:\quad  \sum_{v=1}^n v^{\alpha-1}\omega\Big({\frac {1}{v}}\Big) =
{\mathcal O}  \Big[n^\alpha \omega \Big( {\frac {1}{n}}\Big)\Big],
\]
\[
({\mathscr B}): \quad  \sum\limits_{v=n+1}^\infty \frac 1v\,\omega\bigg(\frac 1v\bigg)={\mathcal O} \bigg[\omega\bigg(\frac 1n\bigg)\bigg].
\]

\begin{theorem}[\cite{Stepanets_Serdyuk_UMZh2002, Abdullayev_Ozkartepe_Savchuk_Shidlich_2019}]\label{Theorem 6.1}  Нехай  $\alpha>0$, $\omega$ -- довільна функція, яка задовольняє умови   1)--4) та умову
$({\mathscr B}_\alpha)$. Для того, щоб функція  $f\in {\mathcal S}^{p} \cap  L_{1,Y}$ належала множині
 ${\mathcal S}^{p} H^{\omega}_{\alpha}\cap  L_{1,Y}$, необхідно та достатньо, щоб
 \[
    E_n^\vartriangle(f)_{_{\scriptstyle {\mathcal S}^p}}={\mathcal O} \Big[ \omega \Big({\frac {1}{n}} \Big) \Big].
\]
\end{theorem}

Функція  $t^r$, $r \le \alpha$, задовольняє умови   $1)$--\,$4)$ та
$({\mathscr B}_\alpha)$.  Тому позначаючи через
${\mathcal S}^{p} H_{\alpha}^r$  множину  ${\mathcal
S}^{p} H^{\omega}_{\alpha}$ при $\omega(t)=t^r$, $0<r\le \alpha,$  отримаємо таке твердження.

\begin{corollary}\label{corollary 6.1.} Нехай $\alpha >0$, $0<r\le \alpha.$
Для того, щоб функція  $f\in {\mathcal S}^{p}\cap  L_{1,Y}$ належала множині ${\mathcal S}^{p} H_{\alpha}^r\cap  L_{1,Y}$,
 необхідно та достатньо, щоб
$$
    E_n^\vartriangle(f)_{_{\scriptstyle {\mathcal S}^p}}={\mathcal O}   ({n^{-r}} ).
$$
\end{corollary}

%%%%%%%%%%%%%%%%%%%%%%%%%%%%%%%%%%%%%%%%%%%%%%%%%%%%%%%%%%%%%

\subsection{Наближення лінійними методами функцій з просторів  ${\mathcal S}^p $}\label{Sp_ linear_chap}

У низці робіт (див., наприклад, \cite{Voicexivskij_2002, Stepanets_Shydlich_UMZh2003, Voicexivskij_Serdyuk_2005, Savchuk_Shidlich_UMZh2007, Savchuk_Shidlich_2014, Shydlich_Zb2003, Shydlich_UMZh2004, Shidlich_UMZh2008} та ін.)  досліджувалися різні  проблеми, пов'язані з апроксимацією   лінійними методами підсумовування рядів Фур'є у просторах ${\mathcal S}^p_\varphi$ та просторах ${\mathcal S}^p$ зокрема.  Так, в роботах
\cite{Shydlich_UMZh2004, Shidlich_UMZh2008} розглядались загальні питанні теорії лінійних методів підсумовування рядів Фур'є (регулярність, насичення) у просторах ${\mathcal S}^p_\varphi$. В  \cite{Stepanets_Shydlich_UMZh2003} та \cite{Shydlich_Zb2003} встановлено точні значення найкращих $n$-членних наближень такими методами $q$-еліпсоїдів у цих просторах. В роботі
 \cite{Voicexivskij_2002} отримано нерівності типу Джексона наближення сумами Зигмунда в просторах ${\mathcal S}^p$. Наведемо результати робіт \cite{Savchuk_Shidlich_UMZh2007, Savchuk_Shidlich_2014}, у яких встановлено прямі та обернені теореми наближення функцій середніми Тейлора-Абеля-Пуассона, і в термінах похибок наближення цими середніми в просторі ${\mathcal S}^p$ отримано конструктивну характеристику  класів функцій, узагальнені
похідні яких належать множинам ${\mathcal S}^pH_\omega$.

\vskip 3mm\noindent {\bf \ref{Sp_ linear_chap}.1. Позначення та постановка задачі.} Нехай $f$ --- довільна функція з простору $L({\mathbb T}^d)$.  Виходячи зі співвідношення (\ref{Fourier series}), розглянемо лінійні оператори
$S^\vartriangle_{n},$ $\sigma^\vartriangle_n$,
$P^\vartriangle_{\varrho,s}$ і $A^\vartriangle_{\varrho,r}$,
визначені на $L({\mathbb T}^d)$ відповідно рівностями
\[
S^\vartriangle_{n}(f)({\bf x})=\sum_{\nu=0}^{n}H_{\nu}(f)({\bf
x}),~n=0,1,\ldots,
\]
\[
\sigma^\vartriangle_n(f)({\bf
x})=\frac{1}{n+1}\sum_{\nu=0}^{n}S^\vartriangle_{\nu}(f)({\bf
x})
=
\sum_{\nu=0}^{n}\left(1-\frac{\nu}{n+1}\right)H_{\nu}(f)({\bf
x}),~n\in\mathbb N,
\]
\[
P^\vartriangle_{\varrho,s}(f)({\bf
x})=\sum_{\nu=1}^{\infty}\varrho^{\nu^s}H_{\nu}(f)({\bf x}),\quad
s>0,~\varrho\in[0,1),
\]
і
\begin{equation}\label{def A}
A^\vartriangle_{\varrho,r}(f)({\bf x})=S^\vartriangle_{r-1}(f)({\bf
x})+\sum_{\nu=r}^{\infty}\lambda_{\nu,r}H_{\nu}(f)({\bf x}),
\end{equation}
де при $r\in\mathbb N$ та $\varrho\in[0,1)$
\[
\lambda_{\nu,r}:=\lambda_{\nu,r}(\varrho):=\sum_{k=0}^{r-1}{\nu\choose k}(1-\varrho)^{k}\varrho^{\nu-k}
=\sum_{k=0}^{r-1}\frac{(1-\varrho)^{k}}{k!}~\frac{d^k}{d\varrho^k}\varrho^{\nu}.
\]

Вирази $S^\vartriangle_n(f)({\bf x})$,
$\sigma^\vartriangle_n(f)({\bf x})$ і
$P^\vartriangle_{\varrho,s}(f)({\bf x})$ називають відповідно
трикутною частинною сумою ряду Фур'є,  трикутною сумою Фейєра і
узагальненою трикутною сумою Абеля--Пуассона функції $f$. Вираз
$A^\vartriangle_{\varrho,r}(f)({\bf x})$ називають трикутною сумою
Тейлора--Абеля--Пуассона функції $f$.

Оператори $P^\vartriangle_{\varrho,s}$ в загальному випадку, як
агрегати наближення функцій функцій однієї змінної, мабуть, вперше
розглядалися в \cite{Bugrov_1963, Bugrov_1972}.
Оператори $A^\vartriangle_{\varrho,r}$ введені  в \cite{Savchuk_UMZh2007}, де в їх термінах дано конструктивну характеристику класів Гарді--Ліпшиця $H^r_p\mathop{\rm Lip}\alpha$ функцій однієї
змінної, голоморфних в одиничному крузі комплексної площини. Апроксимацій властивості цих операторів
вивчалися в роботах
\cite{Savchuk_UMZh2007, Savchuk_Shidlich_UMZh2007, Savchuk_Shidlich_2014, Prestin_Savchuk_Shidlich_2017, Prestin_Savchuk_Shidlich_2019, Chaichenko_Savchuk_Shidlich_2020} та ін. В
частинному  випадку, коли $r=s=1$ оператори
$A^\vartriangle_{\varrho,1}$ та $P^\vartriangle_{\varrho,1}$
збігаються між собою і породжують класичний метод Абеля--Пуассона
підсумовування кратних рядів Фур'є по трикутних областях.

Нагадаємо, що інтегралом Пуассона функції $f\in L({\mathbb T}^d)$ називається
функція $P(f)$, визначена в $[0,1)^d\times\mathbb R^d$ рівністю
\[
f(\mbox{\boldmath${\bf \varrho}$},{\bf x})=\int_{\mathbb T^d}f({\bf x+t})P(\mbox{\boldmath${\bf \varrho}$},{\bf t})d{\bf t},
\]
де
\[
P(\mbox{\boldmath${\bf \varrho}$} ,{\bf t}):=\prod_{j=1}^d\frac{1-\varrho_j^2}{1-2\varrho_j\cos
t_j+\varrho_j^2},~\varrho_j\in[0,1),
\]
--- кратне ядро Пуассона і ${\bf x+t}:=(x_1+t_1,\ldots,x_d+t_d).$

Надалі домовимось під виразом $f(\varrho,{\bf x})$ розуміти
інтеграл Пуассона, в якому  {\boldmath${\bf \varrho}$} --- це вектор з
однаковими координатами, тобто {\boldmath${\bf \varrho}$}$\,=(\varrho,\ldots,\varrho).$

В даному підрозділі вивчаються оператори
$A^\vartriangle_{\varrho,r}$ і $P^\vartriangle_{\varrho,s}$ як
лінійних методів наближення функцій в просторах ${\mathcal S}^p.$ При цьому
основна увага звертається на зв'язок апроксимативних властивостей
сум $A^\vartriangle_{\varrho,r}(f)$ і
$P^\vartriangle_{\varrho,s}(f)$ із диференціальними властивостями
функції $f$, а саме, властивостями похідних, означених в такий спосіб.

Нехай $\psi=\{\psi({\bf k})\}_{{\bf k}\in\mathbb
Z^d}$ -- довільна система комплексних чисел і
\[
\mathscr Z(\psi):=\mathscr Z^d(\psi):=\left\{{\bf k}\in\mathbb Z^d :
\psi({\bf k})=0\right\}.
\]

Надалі вважаємо, що множина $\mathscr Z(\psi)$ має скінченну
кількість елементів.

Якщо для даної функції $f\in L({\mathbb T}^d)$ знайдеться функція $g\in L({\mathbb T}^d)$
така, що
\begin{equation}\label{series for derivative}
S[f]({\bf x})=\sum_{{\bf k}\in\mathscr Z(\psi)}\widehat{f}({\bf
k})e^{i({\bf k,x})}+\sum_{{\bf k}\in\mathbb Z^d}\psi({\bf
k})\widehat{g}({\bf k})e^{i({\bf k,x})},
\end{equation}
то кажуть, що у функції $f$ існує $\psi$-похідна $g$, для якої
використовують позначення $g=f^{\psi}.$ При цьому, якщо $\mathscr
Z(\psi)=\varnothing$, то перша сума в (\ref{series for derivative})
покладається рівною нулеві.

Зрозуміло, що $\psi$-похідна для функцій з простору ${\mathcal S}^p$ є єдиною з
точністю до суми $\sum_{{\bf k}\in\mathscr Z(\psi)}a_{{\bf
k}}e^{i({\bf k, x})}$, де $a_{{\bf k}}$ --- будь-які числа. Дане  означення
$\psi$-похідної пристосоване для потреб досліджень, викладених у цьому підрозділі,
і за суттю  не відрізняється від
поняття $\psi$-похідної О.\,І.~Степанця, наведеного в  підрозділі \ref{Definitions_Sp}

Далі, розглядаються $\psi$-похідні функцій з $L({\mathbb T}^d)$ в таких
двох випадках: $1)\  \psi({\bf k})=\nu^{-r}$ при $|{\bf
k}|_1=\nu,~\nu=0,1,\ldots,~r\ge 0$ і  $2)\ \psi({\bf k})=0$ при $|{\bf k}|_1=0,1,\ldots,r-1$ та
$\psi({\bf k})={(\nu-r)!}/{\nu!}$ при $|{\bf k}|_1=\nu,~\nu\ge r,~r\in\mathbb N$.

При цьому у першому випадку для $\psi$-похідної функції $f$
використовуємо позначення $f^{(r)}$, у другому --- $f^{[r]}$, а при
$r=0$ покладаємо $f^{(0)}=f^{[0]}=f.$
Відмітимо також, що $f^{(1)}=f^{[1]}$.

%\subsection{ Основні результати}\label{{\mathcal S}^p Linear results}

\vskip 3mm\noindent {\bf \ref{Sp_ linear_chap}.2. Прямі та обернені теореми наближення лінійними методами.}
Перейдемо до формулювання основних результатів підрозділу \ref{Sp_ linear_chap}.
 При цьому будемо використовувати
позначення, наведені в підрозділах \ref{DI_TH_Sp}.2 та \ref{DI_TH_Sp}.4.

\begin{proposition}[\cite{Savchuk_Shidlich_UMZh2007, Savchuk_Shidlich_2014}]\label{Linear S^p} Нехай $1\le p<\infty$, $f\in L(\mathbb T^d),~d\in\mathbb N$ і $\omega$ --- довільна функція, яка задовольняє умови  1)--4) та $({\mathscr B})$. Наступні твердження рівносильні:

i) $\left\|S^\vartriangle_n\left(f^{[1]}\right)\right\|_{_{\scriptstyle
{\mathcal S}^p}}=O(n\omega(\frac 1n)),\quad n\to\infty;$

ii) $\left\|f-\sigma^\vartriangle_n\left(f\right)\right\|_{_{\scriptstyle
{\mathcal S}^p}}=O(\omega(\frac 1n)),\quad n\to\infty$.

\noindent Крім цього, якщо виконується одне із тверджень i)--iii), то

iii) $f\in {\mathcal S}^pH_\omega^1$.

\noindent Якщо ж $f\in L_{1,Y}(\mathbb T^d)$, то всі твердження 1)--4) є еквівалентними.

\end{proposition}

Зазначимо, що імплікація $ii)\Rightarrow iii)$ є твердженням типу прямих та обернених теорем для методу Фейєра \cite{Butzer_Nessel_B1971}.

В наступній теоремі даються прямі та обернені теореми наближення функцій оператором $A^\vartriangle_{\varrho,r}$ в просторі ${\mathcal S}^p$ в термінах мажорант $\omega$.

\begin{theorem}[\cite{Savchuk_Shidlich_UMZh2007, Savchuk_Shidlich_2014}]\label{Linear S^p 1} Нехай $1\le p<\infty$, $r\in\mathbb N$, $f\in L(\mathbb T^d),~d\in\mathbb N$ і $\omega$ --- довільна функція, яка задовольняє умови  1)--4) та $({\mathscr B})$. Наступні твердження рівносильні:

i) $\|f-A^\vartriangle_{\varrho,r}(f)\|_{_{\scriptstyle
{\mathcal S}^p}}=O((1-\varrho)^{r-1}\omega(1-\varrho)),\quad\varrho\to 1-;$

ii) $\left\|P(f^{[r]})({\varrho},\cdot)\right\|_{_{\scriptstyle
{\mathcal S}^p}}=O(\frac {\omega(1-\varrho)}{1-\varrho}),\quad\varrho\to 1-;$

\noindent Крім цього, якщо виконується одне із тверджень i) чи ii), то

iii) $f^{[r-1]}\in {\mathcal S}^pH^1_\omega$.

\noindent Якщо ж $f\in L_{1,Y}(\mathbb T^d)$, то всі твердження i)--iii) є еквівалентними.
\end{theorem}

Зазначимо, що імплікація $ii)\Rightarrow iii)$ є твердженням типу теорем Гарді--Літтлвуда \cite{Hardy_Littlewood_1932}.

Наведемо також апроксимаційні властивості сум
$P^\vartriangle_{\varrho,s}(f)$ в просторі ${\mathcal S}^p$. Застосування теореми  \ref{Linear S^p 1} до функції $f=g^{(s-1)}$ зі значенням
параметра $r=1$ і врахування співвідношення
\begin{equation}\label{asymp equality}
\|f-P^\vartriangle_{\varrho,s}(f)\|_{_{\scriptstyle {\mathcal S}^p}}\sim
\|f^{(s-1)}-P^\vartriangle_{\varrho,1}(f^{(s-1)})\|_{_{\scriptstyle
{\mathcal S}^p}},\quad\varrho\to 1-,
\end{equation}
дозволяє записати таке твердження.

\begin{theorem}[\cite{Savchuk_Shidlich_UMZh2007, Savchuk_Shidlich_2014}]\label{Linear S^p 2} Нехай $1\le p<\infty$,  $s\in\mathbb N$,  $f\in L(\mathbb T^d),~d\in\mathbb N$ і $\omega$ --- довільна функція, яка задовольняє умови  1)--4) та $({\mathscr B})$. Наступні твердження рівносильні:

i) $\|f-P^\vartriangle_{\varrho,s}(f)\|_{_{\scriptstyle
{\mathcal S}^p}}=O(\omega(1-\varrho)),\quad\varrho\to 1-;$

ii) $\left\|P(f^{(s)})({\varrho},\cdot)\right\|_{_{\scriptstyle
{\mathcal S}^p}}=O(\frac {\omega(1-\varrho)}{1-\varrho}),\quad\varrho\to 1-;$

\noindent Крім цього, якщо виконується одне із тверджень i) чи ii), то

iii) $f^{(s-1)}\in {\mathcal S}^pH_\omega$.

\noindent Якщо ж $f\in L_{1,Y}(\mathbb T^d)$, то всі твердження i)--iii) є еквівалентними.
\end{theorem}

При $d=1$ простір $L_{1,Y}(\mathbb T^1)$
збігається з простором $L_{1}(\mathbb T^1)$ і тому твердження i)--iii) в твердженні \ref{Linear S^p} і теоремах \ref{Linear S^p 1} та \ref{Linear S^p 2} є рівносильними без жодних застережень.

Зазначимо, що в \cite{Chaichenko_Savchuk_Shidlich_2020} результати твердження \ref{Linear S^p}, а також теорем \ref{Linear S^p 1} та  \ref{Linear S^p 2}, зокрема, розповсюджено на простори типу Орлича ${\mathcal S}_M$.

%%%%%%%%%%%%%%%%%%%%%%%%%%%%%%%%%%%%%%%%%%%%%%%%%%%%%%%%%%%%%%%%%%%%%%%%%%%%%%%%%%%%%%%%%%%%%%%%%%%%%%

%%%%%%%%%%%%%%%%%%%%%%%%%%%%%%%%%%%%%%%%%%%%%%%%%%%%%%%%%%%%%%%%%%%%%%%%%%%%%%%%%%%%%%%%%%%%%%%%%%%%%%%%%%%%%%%%

\begin{thebibliography}{99}



%%%%%%%%%%%%%%%%%%%%%%%%%%%%%%%%%%%%%%%%%%%%%%%%%%%%%%%%%%%%%%%%%%%%%%%%%%%%%%%%%%%%%%%%%%%%%%%%%%%%

\bibitem{Abdullayev_Chaichenko_Shidlich_2020}
                   F.~Abdullayev, S.~Chaichenko, A.~Shidlich,
                  \emph{Direct and inverse approximation theorems of functions in the Musielak-Orlicz type spaces},
                  Math. Inequal. Appl. (accepted for publication, see also arXiv: 2004.09807).


\bibitem{Abdullayev_Ozkartepe_Savchuk_Shidlich_2019}
                   F.\,G.~Abdullayev,  P.~\"{O}zkartepe,  V.\,V.~Savchuk,  A.\,L.~Shidlich,
                  \emph{Exact constants in direct and inverse approximation theorems for functions
                   of several variables in the spaces ${\mathcal S}^p$},
                  Filomat. \textbf{33},  (5) (2019), 1471-1484.

 \bibitem{Abdullayev_Serdyuk_Shidlich_2020}
                   Ф.~Абдуллаєв, А.~Сердюк, А.~Шидліч,
                  \emph{Поперечники функціональних класів, визначених мажорантами узагальнених модулів гладкості в просторах  ${\mathcal S}^{p}$},
                  {Укр. мат. журн.}  (прийнята до друку, див. також arXiv: 2005.05597).

\bibitem{Babenko_1986}
                  А.\,Г.~Бабенко,
                  \emph{О точной константе в неравенстве Джексона в $L^2$},
                  Матем. заметки, \textbf{39}, (5) (1986), 651-664.

\bibitem{Babenko_1960_2}
                  К.\,И.~Бабенко,
                  \emph{О приближении периодических функций многих переменных тригонометрическими полиномами},
                  ДАН СССР. \textbf{132}, (2) (1960), 247-250.

\bibitem{Babenko_1960_5}
                  К.\,И.~Бабенко,
                  \emph{ О приближении одного класса периодических функций многих переменных тригонометрическими полиномами},
                  ДАН СССР. \textbf{132}, (5) (1960), 982-985.

\bibitem{Babenko_Konareva_2019}%17
                В.\,Ф.~Бабенко, С.\,В.~Конарева,
                \emph{Неравенства типа Джексона-Стечкина для аппроксимации элементов гильбертова пространства},
                {Укр. мат. журн.} \textbf{70}, (9) (2018), 1155-1165.

\bibitem{Bari_Stechkin_1956}
                  Н.\,K.~Бари, С.\,Б.~Стечкин,
                  \emph{Наилучшие приближения и дифференциальные свойства двух сопряженных функций},
                   Тр. московского мат. об-ва.  \textbf{5}, (1956), 483-522.

\bibitem{Bugrov_1963}
                  Я.\,С.~Бугров,
                  \emph{Неравенства типа Бернштейна и их применение к исследованию дифференциальных свойств решений дифференциальных уравнений высшего порядка},
                  Mathematica (Cluj).  \textbf{5}, (28) (1963), 5-25.

\bibitem{Bugrov_1972}
                  Я.\,С.~Бугров,
                  \emph{Свойства решений дифференциальных уравнений
                  высшего порядка в терминах весовых классов},
                  Труды Мат. ин-та АН СССР. \textbf{117}, (1972), 47-61.

\bibitem{Butzer_Nessel_B1971}
                  P.~Butzer, R.~Nessel,
                  \emph{Fourier Analysis and Approximation. Volume 1:  One--Dimentional Theory.}
                  Basel--New York, 1971. % --- 554~p.

\bibitem{Vakarchuk_2004}%22
                  С.\,Б.~Вакарчук,
                  \emph{Неравенства типа Джексона и точные значения поперечников классов функций в пространствах $S^p$, $1\leq p< \infty$},
                  Укр. мат. журн. \textbf{56}, (5) (2004), 595-605.

\bibitem{Vakarchuk_2016}%22
                  С.\,Б.~Вакарчук,
                  \emph{Неравенства типа {Д}жексона с обобщенным модулем непрерывности и точные значения $n$-поперечников классов
       $(\psi,\beta)$-дифференцируемых функций в ${L}_2$. {I}},
                  Укр. мат. журн. \textbf{68}, (6) (2016), 723-745.

\bibitem{Vakarchuk_Shchitov_2006}%22
                  С.\,Б.~Вакарчук, А.\,Н.~Щитов,
                  \emph{О некоторых экстремальных задачах теории  аппроксимации функций в пространствах  $S^p$, $1\leq p< \infty $},
                  Укр. мат. журн. \textbf{58}, (3) (2006),  303-316.

\bibitem{Voicexivskij_2002}
                   В.\,Р.~Войцехівський,
                  \emph{Нерівності типу Джексона при наближенні функцій з простору $S^p$ сумами Зігмунда},
                  Теорія наближення функцій та суміжні питання: Праці ІМ НАН України, \textbf{35}, (2002), 33-46.

\bibitem{Voicexivskij_2003}
                  В.\,Р.~Войцехівський,
                  \emph{Поперечники деяких класів простору $S^p$},
                   Екстремальні задачі теорії функцій та суміжні питання:
                   Праці Ін-ту математики НАН Украини, \textbf{46}, (2003),  17-26.

\bibitem{Voicexivskij_UMZh2003}%22
                  В.\,Р.~Войцехівський,
                  \emph{Нерівності типу Джексона в просторі $S^p$},
                  Укр. мат. журн. \textbf{55}, (9) (2003),  1167-1177.

\bibitem{Voicexivskij_Serdyuk_2005}%22
                  В.\,Р.~Войцехівський, А.\,С.~Сердюк,
                  \emph{Нерівності типу Джексона при наближенні функцій з простору ${S}^p$ методом {В}ороного},
                   Проблеми теорії наближення функцій та суміжні питання: Зб. праць Ін-ту математики НАН України,
                   \textbf{2},  (2) (2005),  43--53.

\bibitem{Fuchang_Gao_2010}
                  F.~Gao,
                  \emph{ Exact value of the $n$-term approximation of a diagonal operator},
                   J. Approx. Theory.  \textbf{162} (4) (2010), 646-652.

\bibitem{Gorbachuk_Grushka_Torba_2005}%22
                  М.\,Л.~Горбачук, Я.\,І.~Грушка, С.\,М.~Торба,
                  \emph{Прямі й обернені теореми в теорії наближень методом Рітца},
                  Укр. мат. журн. \textbf{57}, (5) (2005),  633-643.

\bibitem{DeVore_Temlyakov_1995}
                  R.\,A.~DeVore, V.\,N.~Temlyakov,
                  \emph{Nonlinear approximation by trigonometric sums},
                  J. Fourier Anal. Appl. \textbf{2}, (1) (1995), 29-48.



\bibitem{Dung_Temlyakov_Ullrich_2018}
                  D.~Dung, V.~Temlyakov, T.~Ullrich,
                  \emph{ Hyperbolic cross approximation}.
                  %Edited and with a foreword by Sergey Tikhonov.
                  Advanced Courses in Mathematics. CRM Barcelona. Birkh\"{a}user/Springer, Cham, 2018. %xi+218 pp.

\bibitem{Kahan_M1976}
                  Ж.-П.~Кахан,
                  \emph{Абсолютно сходящиеся ряды Фурье}.
                  М.: Мир, 1976. %--- 204~c.

\bibitem{Kolmogoroff_1936}
                  A.~Kolmogoroff,
                  \emph{ \"{U}ber die beste Ann\"{a}herung von Funktionen einer gegebenen Funktionenklasse.},
                  Ann. of Math., Second series, \textbf {37} (1) (1936), 107-110.

\bibitem{Lasuriya_2007}
                 Р.\,А.~Ласурия,
                  \emph{Прямые и обратные теоремы приближения функций, заданных на сфере, в пространстве ${S}^{p,q}(\sigma^m)$},
                  Укр. мат. журн. \textbf{59}, (7) (2007), 901-911.

\bibitem{Lasuriya_2015}
                  Р.\,А.~Ласурия,
                  \emph{Прямые и обратные теоремы приближения функций суммами {Ф}урье–{Л}апласа в пространствах ${S}^{p,q}(\sigma^{m-1})$},
                 Матем. заметки. \textbf{98}, (4) (2015),  530-543.


\bibitem{Li_2010}
                 R.\,S.~Li, Y.\,P.~Liu,
                  \emph{Best $m$-term one-sided trigonometric approximation of some function classes defined by a kind of multipliers},
                  Acta Mathematica Sinica, English Series.  \textbf{26} (5) (2010), 975--984.


\bibitem{Pinkus_1985}
                 A.~Pinkus,
                  \emph{$n$-widths in approximation theory}.
                  Springer-Verlag, 1985. %- 291~p.


\bibitem{Prestin_Savchuk_Shidlich_2017}
                 Ю.~Престін, В.\,В.~Савчук, А.\,Л.~Шидліч,
                  \emph{Прямі та обернені теореми наближення $2\pi$-періодичних функцій
                  середніми Тейлора-Абеля-Пуассона},
                  Укр. мат. журн. \textbf{69}, (5) (2017), 657-669.

\bibitem{Prestin_Savchuk_Shidlich_2019}
                  J.~Prestin, V.\,V.~Savchuk, A.\,L.~Shidlich,
                  \emph{Approximation theorems for multivariate Taylor-Abel-Poisson means},
                  Stud. Univ. Babe\c{s}-Bolyai Math. \textbf{64}, (3) (2019), 313-329.

\bibitem{Romanyuk_2012}
                 А.\,С.~Романюк,
                 \emph{  Аппроксимативные характеристики классов периодических функций многих переменных},
                  Праці Інституту математики НАН України. \textbf{40} (2012).%, 352~p.

\bibitem{Rukasov_UMZh2003}
                  В.\,И.~Рукасов,
                  \emph{ Наилучшие $n$-членные приближения в пространствах с несимметричной метрикой},
                  Укр. мат. журн. \textbf{55}, (4) (2003), 500-509.

\bibitem{Savchuk_UMZh2007}
                  В.\,В.~Савчук,
                  \emph{Наближення голоморфних функцій середніми Тейлора-Абеля-Пуассона}
                  Укр. мат. журн. \textbf{59}, (9) (2007), 1253-1260.

\bibitem{Savchuk_Shidlich_UMZh2007}
                  В.\,В.~Савчук,  А.\,Л.~Шидліч,
                  \emph{ Наближення функцій багатьох змінних лінійними методами в просторах $S^p$},
                   Проблеми теорії наближення функцій та суміжні питання:
                   Зб. праць Ін-ту математики НАН України, \textbf{4}, (1) (2007),  302-317.

\bibitem{Savchuk_Shidlich_2014}
                  V.\,V.~Savchuk, A.\,L.~Shidlich,
                  \emph{ Approximation of functions of several variables by linear methods in the space $S^p$},
                  Acta Sci. Math. (Szeged). \textbf{80}, (3-4) (2014),  477-–489.

\bibitem{Serdyuk_2003}
                 А.\,С.~Сердюк,
                  \emph{Поперечники в просторі $S^p$ класів функцій, що означаються модулями  неперервності їх $\psi$-похідних}
                   Екстремальні задачі теорії функцій та суміжні питання:
                   Праці Ін-ту математики НАН Украини, \textbf{46}, (2003), 229-248.

\bibitem{Serdyuk_Stepanyuk_2015}
                  А.\,С.~Сердюк, Т.\,А.~Степанюк,
                  \emph{ Оцінки найкращих $m$-членних тригонометричних наближень класів аналітичних функцій},
                   Допов. НАН України, \textbf{2},  (2015),  32-37.

\bibitem{Stepanets_Preprint2001}
                  А.\,И.~Степанец,
                  \emph{   Аппроксимационные характеристики пространств $S^p_\varphi$}.
                   Киев, 2001 %. -85 с. -
                   (Препр./ НАН Украины, Ин-т математики; 2001.2).

\bibitem{Stepanets_UMZh2001_3}
                  А.\,И.~Степанец,
                  \emph{   Аппроксимационные характеристики пространств $S^p_\varphi$},
                   Укр. мат. журн. \textbf{53}, (3) (2001),  392-416.

\bibitem{Stepanets_UMZh2001_8}
                  А.\,И.~Степанец,
                  \emph{   Аппроксимационные характеристики  пространств $S^p_\varphi$ в разных метриках},
                   Укр. мат. журн. \textbf{53}, (8) (2001),   1121-1146.


\bibitem{Stepanets_M2002_1}
                 А.\,И.~Степанец,
                 \emph{Методы теории приближений: В 2 ч.}
                 Праці Ін-ту математики НАН Украини.  Математика та її застосування, \textbf{40}, (I) (2002). % 427~с.


\bibitem{Stepanets_M2002_2}
                 А.\,И.~Степанец,
                 \emph{Методы теории приближений: В 2 ч.}
                 Праці Ін-ту математики НАН Украини.  Математика та її застосування, \textbf{40}, (II) (2002). % 468~с.

\bibitem{Stepanets_UMZh2003_10}
                  А.\,И.~Степанец,
                  \emph{  Экстремальные задачи теории приближений в линейных пространствах},
                   Укр. мат. журн. \textbf{55}, (10) (2003),   1392-1423.

\bibitem{Stepanets_UMZh2004_10}
                  А.\,И.~Степанец,
                  \emph{Наилучшие приближения $q$-эллипсоидов в просторанствах ${S}^{p,\mu}_\varphi$},
                   Укр. мат. журн. \textbf{56}, (10) (2004),   1378-1383.

\bibitem{Stepanets_UMZh2005_4}
                  А.\,И.~Степанец,
                  \emph{  Наилучшие $n$-членные приближения с ограничениями},
                   Укр. мат. журн. \textbf{57}, (4) (2005),    533-553.

\bibitem{Stepanets_UMZh2006_1}
                  А.\,И.~Степанец,
                  \emph{  Задачи теории приближений в линейных пространствах},
                   Укр. мат. журн. \textbf{58}, (1) (2006),   47-92.

\bibitem{Stepanets_B2006}
                  A.\,I.~Stepanets,
                  \emph{Methods of approximation theory}. VSP, Leiden, 2005. %xviii+919 pp.

\bibitem{Stepanets_Rukasov_UMZh2003_2}
                  А.\,И.~Степанец, В.\,И.~Рукасов,
                  \emph{  Пространства $\ S^p$ с несимметричной метрикой},
                  Укр. мат. журн. \textbf{55}, (2) (2003),  264-277.

\bibitem{Stepanets_Rukasov_UMZh2003_5}
                  А.\,И.~Степанец, В.\,И.~Рукасов,
                  \emph{Наилучшие "cплошные" $n$-членные приближения в пространствах ${S}^p_\varphi$},
                  Укр. мат. журн. \textbf{55}, (5) (2003),  801-811.

\bibitem{Stepanets_Serdyuk_UMZh2002}
                  А.\,И.~Степанец, А.\,С.~Сердюк,
                  \emph{Прямые и обратные теоремы приближения функций в пространстве $S^p$},
                  Укр. мат. журн. \textbf{54}, (1) (2002),  106-124.

\bibitem{Stepanets_Shydlich_UMZh2003}
                  О.\,І.~Степанець, А.\,Л.~Шидліч,
                  \emph{  Найкращі $n$-членні наближення $\Lambda$-методами в просторах $S_\varphi^p$},
                  Укр. мат. журн. \textbf{55}, (8) (2003),  1107-1126.

\bibitem{Stepanets_Shidlich_Pr_2007}
                  А.\,И.~Степанец, А.\,Л.~Шидлич,
                  \emph{  Экстремальные задачи для интегралов от неотрицательных функций}.
                   Киев, 2007 %. -103 с. -
                  (Препринт / НАН Украины. Ин-т математики; 2007.2).


\bibitem{Stepanets_Shidlich_NK2007}
                  А.\,И.~Степанец, А.\,Л.~Шидлич,
                  \emph{  О порядках наилучших приближений интегралов функций
                   при помощи интегралов ранга $\sigma$},
                   Нелін. колив. \textbf{10}, (4) (2007), 528--559.


\bibitem{Stepanets_Shidlich_JAT_2010}
                 A.\,I.~Stepanets, A.\,L.~Shidlich,
                  \emph{ Best approximations of integrals by integrals of finite rank},
                  J. Approx. Theory \textbf{162}, (2) (2010), 323-348.

\bibitem{Stepanets_Shidlich_IZV_2010}
                  А.\,И.~Степанец, А.\,Л.~Шидлич,
                  \emph{  Экстремальные задачи для интегралов от неотрицательных функций},
                  Изв. РАН. Сер. матем. \textbf{74}, (3) (2010), 169–224.

\bibitem{Sterlin_1972}%21
                 М.\,Д.~Стерлин,
                  \emph{Точные постоянные в обратных теоремах теории приближений},
                  Докл. АН СССР. \textbf{202}, (3) (1972), 545-547.


\bibitem{Stechkin_1955_DAN}
                  С.\,Б.~Стечкин,
                  \emph{ Об абсолютной сходимости ортогональных рядов},
                  Докл. АН СССР. \textbf{102}, (1) (1955),   37-40.



\bibitem{Taikov_1976}
                  Л.\,В.~Тайков,
                  \emph{ Неравенства, содержащие наилучшие приближения и модуль непрерывности функций из $L_2$},
                  Матем. заметки. \textbf{20}, (3) (1976), 433-438.

\bibitem{Taikov_1979}
                  Л.\,В.~Тайков,
                  \emph{Структурные и конструктивные характеристики функций из $L_2$},
                 Матем. заметки. \textbf{25}, (2) (1979),  217-223.

\bibitem{Temlyakov_1986}
                  В.\,Н.~Темляков,
                  \emph{ Приближение функций с ограниченной смешанной производной}.
                  Тр. МИАН СССР. \textbf{178} (1986), 113~c.

\bibitem{Temlyakov_B1993}
                  V.\,N.~Temlyakov,
                  \emph{ Approximation of periodic functions}.
                  Computational Mathematics and Analysis Series. Commack, New York: Nova Science Publ.,  1993. %-- 419~p.

\bibitem{Temlyakov_Greedy_1998}
                  V.\,N.~Temlyakov,
                  \emph{ Greedy algorithm and $m$-term trigonometric approximation},
                  Constr. Approx. \textbf{14}, (4)  (1998), 569--587.

\bibitem{Temlyakov_B2011_Greedy}
                  V.\,N.~Temlyakov,
                  \emph{Greedy approximation}.
                  Cambridge Monographs on Applied and Computational Mathematics, 20.
                  Cambridge: Cambridge University Press, 2011. %--- 418~p.

\bibitem{Temlyakov_B2015}
                  V.\,N.~Temlyakov,
                  \emph{Sparse approximation with bases}.
                  Advanced Courses in Mathematics. %Edited by S.~Tikhonov.
                  CRM Barcelona. Birkhauser / Springer, Basel, 2015. %--- 261~p.

\bibitem{Timan_M_MS1958}
            М.\,Ф.~Тиман,
            Обратные теоремы конструктивной теории функций в пространствах $L_p$, $(1\le p\le \infty)$,
            Матем. сб., \textbf{46(88)}, (1) (1958),  125-132.

\bibitem{Timan_M_DAN1958}
            М.\,Ф.~Тиман,
            Обратные теоремы конструктивной теории функций многих переменных
            Докл. АН СССР, \textbf{120}, (6) (1958),  1207-1209.

\bibitem{Timan_M2009}
                М.\,Ф.~Тиман,
                  \emph{Аппроксимация и свойства периодических функций}.
                  Киев: Наук. думка, 2009.


\bibitem{Hardy_Littlewood_1932}
                G.\,H.~Hardy, J.\,E.~Littlewood,
                \emph{ Some properties of fractional integrals. II},
                Math. Zeitschr. \textbf{34}, (1) (1932), 403-439.

\bibitem{Hardy_B1948}
                 Г.\,Г.~Харди, Д.\,Е.~Литтльвуд, Г.~Полиа,
                 \emph{Неравенства}. Москва: Изд-во иностр. лит., 1948. %-- 456~c.

\bibitem{Chaichenko_Savchuk_Shidlich_2020}
                  S.~Chaichenko, V.~Savchuk, A.~Shidlich,
                  \emph{Approximation of functions by linear summation methods in the Orlicz-type spaces},
                   Укр. матем. вiсник. \textbf{17}, (2) (2020), 152-170.

\bibitem{Chaichenko_Shidlich_2018}
                  С.\,О.~Чайченко, А.\,Л.~Шидліч,
                  \emph{Апроксимативнi характеристики модулярних просторiв Орлича},
                  Укр. матем. вiсник. \textbf{15}, (2) (2018), 194-209.

\bibitem{Chernykh_1967}
                  Н.\,И.~Черных,
                  \emph{О неравенстве Джексона в  $L_2$},
                  Тр. МИАН СССР, \textbf{88}, (1967),  71-74.

\bibitem{Chernykh_1967_MZ}
                  Н.\,И.~Черных,
                  \emph{О наилучшем приближении периодических функций тригонометрическими полиномами в  $L_2$},
                  Матем. заметки, \textbf{2}, (5) (1967),  513-522.

\bibitem{Shydlich_Zb2003}
                  А.\,Л.~Шидліч,
                  \emph{  Найкращі $n$-членні наближення $\Lambda$-методами в просторах $S_\varphi^p$},
                   Екстремальні задачі теорії функцій та суміжні питання:
                   Праці Ін-ту математики НАН Украини, \textbf{46}, (2003), 283-306.

\bibitem{Shydlich_UMZh2004}
                  А.\,Л.~Шидліч,
                  \emph{ Про насичення лінійних методів підсумовування рядів Фур'є у просторах $S_\varphi^p$},
                  Укр. мат. журн. \textbf{56}, (1) (2004),  133-138.

\bibitem{Shidlich_UMZh2008}
                  А.\,Л.~Шидліч,
                  \emph{ Насичення лінійних методів підсумовування  рядів Фур'є в просторах $S_\varphi^p$},
                  Укр. мат. журн. \textbf{60}, (6) (2008),  815-828.


\bibitem{Shydlich_Zb_2008}
                  А.\,Л.~Шидліч,
                  \emph{Апроксимативні характеристики просторів $S_\Phi^p$},
                   Теорія наближення функцій та суміжні питання:
                   Зб. праць Ін-ту математики НАН Украини, \textbf{5}, (1) (2008), 404-430.


\bibitem{Shydlich_UMZh2009}
                  А.\,Л.~Шидліч,
                  \emph{Порядкові рівності для деяких функціоналів та їх застосування
                  до оцінок найкращих $n$-членних наближень та поперечників},
                  Укр. мат. журн. \textbf{61}, (10) (2009), 1403-1423.

\bibitem{Shidlich_Zb2011}
                  А.\,Л.~Шидліч,
                  \emph{Порядкові оцінки найкращих $n$-членних ортогональних тригонометричних наближень класів функцій ${\mathcal F}_{q,\infty}^{\psi}$ в просторах $L_p(\mathbb T^d)$},
                    Проблеми теорії наближення функцій та суміжні питання:
                   Зб. праць Ін-ту математики НАН Украини, \textbf{8}, (1) (2011), 302-317.


\bibitem{Shidlich_Zb2013}
                  А.\,Л.~Шидліч,
                  \emph{Порядкові оцінки для деяких апроксимаційних характеристик},
                   Теорія наближення функцій та суміжні питання:
                   Зб. праць Ін-ту математики НАН Украини, \textbf{10}, (1) (2013), 304-337.


\bibitem{Shidlich_Zb2014}
                  А.\,Л.~Шидліч,
                  \emph{Порядкові оцінки функціоналів, в термінах яких виражаються найкращі $n$-членні наближення класів ${\mathcal F}_{q,r}^{\psi}$},
                   Теорія наближення функцій та суміжні питання:
                   Зб. праць Ін-ту математики НАН Украини, \textbf{11}, (3) (2014),  287-314.


\bibitem{Shidlich_2016}
                  A.\,L.~Shidlich,
                  \emph{Nonlinear approximation of the classes ${\mathcal F}_{q,r}^{\psi}$ of
                  functions of several variables  in the integral metrics},
                  Математичні проблеми механіки та обчислювальної математики:
                  Зб. праць Ін-ту математики НАН України, \textbf{13}, (3) (2016), 256-274.

\bibitem{Shidlich_Chaichenko_LM_2014}
                  А.\,Л.~Шидліч, С.\,О.~Чайченко,
                  \emph{Деякі екстремальні задачі в просторах Орліча},
                  Матем. студії. \textbf{42}, (1) (2014),   21-32.

\bibitem{Shidlich_Chaichenko_lp_2014}
                  А.\,Л.~Шидліч, С.\,О.~Чайченко,
                  \emph{Апроксимаційні характеристики діагональних операторів в просторах $l_{\bf p}$},
                  Математичні проблеми механіки та обчислювальної математики:
                  Зб. праць Ін-ту математики НАН України, \textbf{11}, (2) (2014),  399-412.


\bibitem{Shidlich_Chaichenko_Orlicz_lM_2015}
                   A.\,L.~Shidlich,  S.\,O.~Chaichenko,
                  \emph{Approximative properties of diagonal operators in Orlicz spaces},
                   Numer. Funct. Anal. Optim.  \textbf{36}, (10) (2015),  1339–1352.



\bibitem{Schmidt_1906}
                    E.~Schmidt,
                  \emph{Zur Theorie der linearen und nichtlinearen Integralgleichungen. I},
                   Math. Annalen.  \textbf{63}, (1906),   433-476.




\end{thebibliography}
\end{document}